\documentclass[onefignum,onetabnum]{siamart250211}

\usepackage{paralist}
\usepackage{amsmath,amsfonts}
\usepackage{graphicx}
\usepackage{xcolor}
\usepackage{rotating}
\usepackage{colortbl}
\usepackage{enumitem}
\usepackage{soul}
\usepackage{makecell,multirow,multicol,booktabs}
\usepackage[margin=1in]{geometry}

\usepackage{cleveref}
\definecolor{red}{RGB}{200,16,46}
\hypersetup{
  colorlinks=true,
  linkcolor=red,
  citecolor=red,
  urlcolor=red
}

\newsiamremark{remark}{Remark}

\def\ititle{Neural Networks for Bayesian Inverse Problems Governed by a Nonlinear ODE}
\headers{NNs for Bayesian Inverse Problems}{G. Villalobos, J. Rudi, A. Mang}
\title{\ititle\thanks{\funding{This work was partly supported by the National Science Foundation ({\bf NSF}) under the awards DMS-2012825 and DMS-2145845. Any opinions, findings, and conclusions or recommendations expressed herein are those of the authors and do not necessarily reflect the views of NSF. This work was completed in part with resources provided by the Research Computing Data Core at the University of Houston. Access to the Neocortex Supercomputer at the Pittsburgh Supercomputing Center was awarded through the Neocortex Spring 2023 Call.}}}

\author{German Villalobos\thanks{Department of Mathematics, University of Houston, Houston, TX, USA (\email{andreas@math.uh.edu}).}
\and Johann Rudi\thanks{Department of Mathematics, Virginia Tech, Blacksburg, VA, USA (\email{jrudi@vt.edu}).}
\and Andreas Mang\footnotemark[2]}

\newcommand{\iacr}[1]{{\bf #1}}
\newcommand\ipoint[1]{\textbf{#1}}
\newcommand\fhn{FitzHugh--Nagumo}

\newcommand{\defeq}{\ensuremath{\mathrel{\mathop:}=}}

\newcommand\vect[1]{\ensuremath{\boldsymbol{#1}}}
\newcommand\mat[1]{\ensuremath{\boldsymbol{#1}}}
\newcommand\ns[1]{\ensuremath{\mathbb{#1}}}
\renewcommand\d{\ensuremath{\mathrm{d}}}
\newcommand{\diag}{\ensuremath{\operatorname{diag}}}
\newcommand\spd[1]{\ensuremath{\operatorname{Sym}^+(#1)}}
\newcommand\sym[1]{\ensuremath{\operatorname{Sym}(#1)}}

\DeclareMathOperator*{\minopt}{minimize}

\newcommand\nfnf{$(0\mid0)$}
\newcommand\nn{$(1\mid1)$}
\newcommand\nfn{$(0\mid1)$}
\newcommand\ts{TS} % time series
\newcommand\fc{FC} % fourier coefficients
\newcommand\tsfc{TS\&FC} % time series and fourier coefficients
\newcommand{\ms}{\phantom{$-$}}
\newcommand{\adjtab}{\footnotesize\renewcommand{\arraystretch}{0.5}}

\usepackage{siunitx}
\newcommand{\decc}[1]{\num[round-mode=places,round-precision=3]{#1}}
\newcommand{\decs}[1]{\num[round-mode=places,round-precision=1]{#1}}

\newcommand\sci[1]{\num[
scientific-notation=true,
round-precision=1,
fixed-exponent=1,
detect-weight=true,
detect-family=true,
round-mode=places,
retain-explicit-plus=true,
mode=text,fixed-exponent=0,
retain-explicit-plus=true,
output-exponent-marker=\text{e}]{#1}}

\newcommand{\inum}[1]{\num[group-minimum-digits=3,group-separator={,}]{#1}}

\ifpdf\hypersetup{
  pdftitle={\ititle},
  pdfauthor={G. Villalobos, J. Rudi, A. Mang}
}\fi

\begin{document}

\maketitle

\begin{abstract}
We investigate the use of neural networks (NNs) for the estimation of hidden model parameters and uncertainty quantification from noisy observational data for inverse parameter estimation problems. We formulate the parameter estimation as a Bayesian inverse problem. We consider a parametrized system of nonlinear ordinary differential equations (ODEs), which is the \fhn\ model. The considered problem exhibits significant mathematical and computational challenges for classical parameter estimation methods, including strong nonlinearities, nonconvexity, and sharp gradients. We explore how NNs overcome these challenges by approximating reconstruction maps for parameter estimation from observational data. The considered data are time series of the spiking membrane potential of a biological neuron. We infer parameters controlling the dynamics of the model, noise parameters of autocorrelated additive noise and noise modelled via stochastic differential equations, as well as the covariance matrix of the posterior distribution to expose parameter uncertainties---all with just one forward evaluation of an appropriate NN. We report results for different NN architectures and study the influence of noise on prediction accuracy. We also report timing results for training NNs on dedicated hardware. Our results demonstrate that NNs are a versatile tool to estimate parameters of the dynamical system, stochastic processes, as well as uncertainties as they propagate though the governing ODE.
\end{abstract}

\begin{keywords}
Neural Networks, Deep Learning for Inverse Problems, Bayesian Inverse Problems, \fhn, Nonlinear ODEs
\end{keywords}

\begin{MSCcodes}
68T07, 62F15, 34A55
\end{MSCcodes}

\section{Introduction}\label{s:intro}
We develop methodology for and test the use of neural networks (\iacr{NN}s) for the estimation of
\begin{inparaenum}[\it (i)]
\item hidden parameters of physical models,
\item hidden parameters of stochastic processes, and
\item model and data uncertainties from noisy observations.
\end{inparaenum}
The present work extends our contributions described in~\cite{rudi2022:parameter} in categories \begin{inparaenum}[\it (i)] \item and \item\!\!\!; and it introduces \item (see \Cref{s:contributions} for specific details)\end{inparaenum}. We assume that the observational data can be explained through the solution of a dynamical system. The inverse problem to estimate parameters is typically formulated as a variational problem governed by the underlying system~\cite{tarantola2005:inverse, Troltzsch:2010a,Mang:2019a,Mang:2017a,Mang:2017b}. In a deterministic setting, we seek point estimates for the hidden model parameters and compute these with numerical optimization~\cite{Nocedal:2006a}. In a statistical setting, we seek samples from a posterior distribution~\cite{Stuart:2010a,Calvetti:2023a, tenorio2017:introduction,Kaipio:2006a}, which in practice is often a probability density function of the hidden model parameters conditioned on the data.

\begin{figure}
\centering
\includegraphics[width=0.95\textwidth]{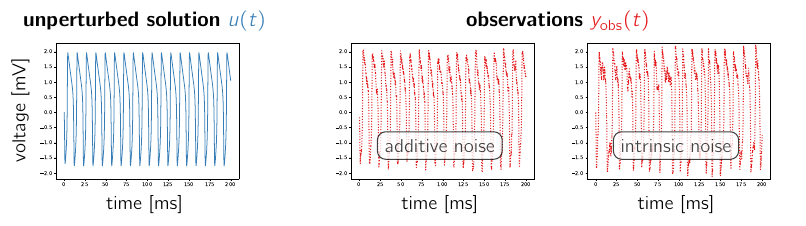}
\caption{Illustration of the model dynamics and noise perturbations. We show time series of the membrane potential $u$ for $\tau = 200\,\text{ms}$. Left: Solution of the \fhn\ model (membrane potential $u$) for representative model parameters $\vect{\theta}_{\text{dyn}} \in \mathbb{R}^3$. Right: Plots of the observational data $y_{\text{obs}}$ for an additive and an intrinsic noise model. We limit the observational data associated with the \fhn\ model to the membrane potential $u$. Our goal is to learn a map, $g_{\Xi}$, from the observational data $y_{\text{obs}}$ to the hidden parameters, $\vect{\theta}_{\text{pred}} = g(y_{\text{obs}})$, such that $\vect{\theta}_{\text{pred}} \approx \vect{\theta}_{\text{true}}$.\label{f:fwdsim}}
\end{figure}

\begin{figure}
\centering
\includegraphics[width=0.75\textwidth]{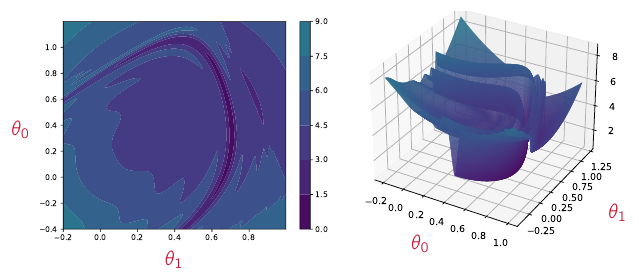}
\caption{Illustration of the parameter estimation problem. We show a contour plot and the associated surface plot of the objective that enters our variational optimization problem as a function of the hidden parameters $\theta_0$ and $\theta_1$ of $\vect{\theta}_{\text{dyn}} = (\theta_0, \theta_1, \theta_2) \in \mathbb{R}^3$ that control the model output. We can observe that the optimization problem is nonconvex with a volatile change in objective values that features a narrow, flat optimality zone (dark blue region in plot on the left) with sharp gradients (surface plot on the right).\label{f:lscp}}
\end{figure}

\begin{figure}
\centering
\includegraphics[width=\textwidth]{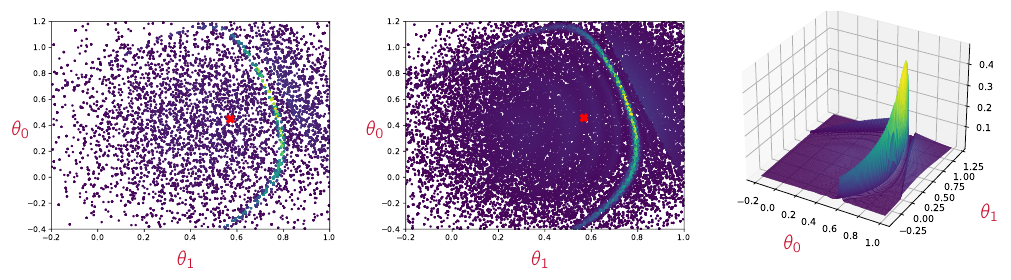}
\caption{Illustration of the posterior distribution. We show samples drawn from the posterior distribution as a function of the hidden parameters $\theta_0$ and $\theta_1$ of $\vect{\theta}_{\text{dyn}} = (\theta_0, \theta_1, \theta_2) \in \mathbb{R}^3$ that control the model output. We can observe that the problem is nonconvex with a volatile change in posterior values that features a narrow optimality zone (bright region) with sharp transitions. We note that the plot in \Cref{f:lscp} (the objective function) corresponds to the negative log posterior of the posterior distribution. The samples in the visualization on the left and middle are generated using a Metropolis Hastings Markov chain Monte Carlo sampling algorithm (left: 20\,000 samples; middle: 100\,000 samples).\label{f:posterior}}
\end{figure}

The considered inverse problem is governed by the \fhn\ model~\cite{Fitzhugh:1961a, nagumo1962:active, cebrian2024:six}, a nonlinear system of ordinary differential equations (\iacr{ODE}s). It simulates electrical voltage spikes in biological neurons that are generated in response to current stimuli. The \fhn\ model is a simplification of neural dynamics with a system of two ODEs; while more complex dynamical systems utilize Hodgkin--Huxley-type models based on~\cite{HodgkinHuxley:1952}. We consider an inverse parameter estimation problem that at first may appear computationally simple to solve, but it poses significant mathematical and computational challenges.
Key \emph{mathematical challenges} associated with the inverse parameter estimation problem are
\begin{inparaenum}[\it (i)]
\item nonlinearity of the forward model,
\item a multiplicative coupling and strong nonlinear dependencies between model parameters,
\item strong nonlinearities in the inverse problem,
\item nonconvexity of the data mismatch term with sharp gradients~\cite{ramsay2007:parameter},
\item and weak a priori information on adequate model parameter values~\cite{Gutenkunst:2007a,Prinz:2004a}.
\end{inparaenum}
These challenges render the inverse problem as a global optimization problem with multiple local minima~\cite{Wang:2022a}. Therefore, most classical methods for the solution of this inverse problem become ineffective and motivate our work. In the present work, we explore if the use of NNs to learn the reconstruction map for hidden parameters given noisy observations allows us to reliably generate parameter predictions under the stated mathematical challenges. We illustrate the data associated with the considered parameter estimation problem in \Cref{f:fwdsim}. We show a representative optimization landscape in \Cref{f:lscp} to highlight some of the mathematical challenges. We visualize a representative posterior distribution in \Cref{f:posterior}.

We note that the low-dimensional character (only three scalar parameters control the state of the dynamical system) make for an excellent testbed to explore algorithms without facing the challenge of a computationally intractable forward operator.

Estimating the parameters of the dynamical system, $\vect{\theta}_{\text{dyn}}$, is challenging on its own. It becomes even more challenging when the goal is to also infer parameters of a noise model, $\vect{\theta}_{\text{noise}}$, from noise that pollutes the observed state of the dynamical system. Estimating both simultaneously is rarely attempted but important for inverse problems on noisy measurements. In addition, we provide methodology to expose how uncertainties propagate throughout the dynamical system. These uncertainties can arise from various sources. One source of error is the discrepancy between the mathematical model subject to inference and the ``true'' physical system and the choice of the dynamical system class (``model form uncertainty''). These discrepancies alongside inaccuracies in specifying initial or boundary conditions of the dynamical system can introduce significant modeling bias. Measurement errors, stemming from sensor noise, finite resolution of measuring equipment, or sampling errors, are another source of error. Furthermore, numerical discretization and limited computational precision contribute to the uncertainty. There is also uncertainty in the parameters that control the dynamical system due to a lack of information or intrinsic variability. Uncertainty quantification offers a rigorous mathematical framework for exposing these errors. For effective uncertainty quantification in inverse problems governed by dynamical systems, it is essential to employ a careful, often hierarchical approach to treating sources of error. Each source affects the overall uncertainty in the inferred parameters or states differently. Modern methods for inference go beyond parameter estimation and strive to quantify and manage these uncertainties, providing robust and credible inferences. Numerous approaches exist to address individual sources of error, including modeling unknowns as random variables, model averaging, and deriving a posteriori error estimates. Our work specifically aims to provide parameter estimates that control error perturbations in the measured signals, while also exposing uncertainties by estimating the covariance matrix associated with the posterior distribution.

\subsection{Outline of the Method}\label{s:method-outline}

Our work aims at designing a NN surrogate model that learns the mapping from noise observations to various quantities of interest. Let $\vect{\theta}_{\text{dyn}} \in \Theta_{\text{dyn}}$ denote the parameters of a considered dynamical system (i.e., the parameters of interest), where $\Theta_{\text{dyn}}$ is the set of feasible parameters. Let $f : \Theta_{\text{dyn}} \to \mathcal{Y}$ denote the \emph{parameter-to-observation map}, where $\mathcal{Y}$ is the set of all realizable observations. The parameter-to-observation map is defined as the composition $f = o \circ e$. The two maps comprising $f$ are $e : \Theta_{\text{dyn}} \to \mathcal{U}$ that maps parameters to realizable states, $\mathcal{U}$, and $o : \mathcal{U} \to \mathcal{Y}$ that restricts and/or transforms states to what can be observed.  We denote $e(\vect{\theta}_{\text{dyn}})$ as the \emph{forward map} (i.e., the solution operator of a dynamical system), $u = e(\vect{\theta}_{\text{dyn}}) \in
\mathcal{U}$ as the state variable, $o(u)$ as the \emph{observation operator}, and $y_{\text{obs}} \in \mathcal{Y}$ as observed data. We relate the parameters of the dynamical system $\vect{\theta}_{\text{dyn}}$ to the observed data $y_{\text{obs}} \in \mathcal{Y}$ by
\[
y_{\text{obs}} = f(\vect{\theta}_{\text{dyn}}) + \eta
= o \circ e(\vect{\theta}_{\text{dyn}}) + \eta
= o(u) + \eta,
\]

\noindent where $\eta$ represents additive observation noise. In the \emph{inverse parameter-estimation problem} we seek the (hidden) parameters $\vect{\theta}_{\text{dyn}}$ of the dynamical system that best explain $y_{\text{obs}}$, i.e., we try to find predictions $u_{\text{pred}} = e(\vect{\theta}_{\text{pred}})$ so that $f(\vect{\theta}_{\text{pred}}) \approx y_{\text{obs}}$. Informally, we can cast this problem as identifying the inverse map $g_{\text{dyn}} : \mathcal{Y} \to  \Theta_{\text{dyn}}$,
\[
\vect{\theta}_{\text{pred}}
= g_{\text{dyn}}(y_{\text{obs}})
= f^{-1}(y_{\text{obs}})
= (o \circ e)^{-1}(y_{\text{obs}}),
\]

\noindent that relates measurements $y_{\text{obs}}$ (NN input) to parameters $\vect{\theta}_{\text{pred}}$ (NN output). We note that this mapping is not well-defined due to the ill-posedness of the inverse problem.

In the present work, we first propose and test an end-to-end machine learning framework to generate a NN approximation to $g_{\text{dyn}}$.  That is, we learn a NN surrogate model for the \emph{reconstruction map}. Second, we extend the NN surrogate for $g_{\text{dyn}}$ to a reconstruction map that maps from $\mathcal{Y}$ to $\Theta$, where $\Theta$ does not only comprise the parameters $\vect{\theta}_{\text{dyn}}$ of the dynamical system, but it also includes the parameters $\vect{\theta}_{\text{noise}}$ of different noise models and matrices $\mat{\Gamma}_{\text{post}}$. The matrices $\mat{\Gamma}_{\text{post}} \succ 0$ are positive semi-definite covariance matrices that form a Laplace approximation of the (generally non-Gaussian) posterior at the maximum a-posteriori ({\bf MAP}) point; they are important to expose uncertainties in our Bayesian inverse problem. In summary, our target is computationally learning the NN surrogate $g_{\Xi} : \mathcal{Y} \to \Theta$ by optimizing the NN parameters $\Xi$. As such, the dimension of the NN prediction/output space $\Theta$ significantly increases compared to the dimension of the space $\Theta_{\text{dyn}}$.

To generate the data for training $g_{\Xi}$, we use forward simulations for samples $\vect{\theta}_{\text{dyn}}^i$, $i = 1,\ldots,n_{\text{train}}$, and $\vect{\theta}_{\text{noise}}^i$ to generate $y_{\text{obs}}^i$ and $\mat{\Gamma}_{\text{post}}^i$. Subsequently, we train a NN to learn the mapping from the NN input $y_{\text{obs}}$ to the NN output $\vect{\theta}_{\text{pred}}$ comprised of $\vect{\theta}_{\text{dyn}}$, $\vect{\theta}_{\text{noise}}$, and $\mat{\Gamma}_{\text{post}}$.

We note that our framework does not require to draw samples from the posterior distribution to construct the training data for $\mat{\Gamma}_{\text{post}}$. We approximate $\mat{\Gamma}_{\text{post}}$ based on the inverse of the Hessian $\mat{H}$ associated with the negative log-posterior at known (true) model parameters $\vect{\theta}_{\text{dyn}}^i$. We use adjoints to approximate the Hessians~\cite{BuiThanh:2013a, petra2014:computational, Martin:2012a,Stuart:2010a} resulting in a scalable method that has a fixed cost independent of the dimensions of the parameter and state spaces, $\Theta_{\text{dyn}}$ and $\mathcal{U}$, respectively.

\subsection{Contributions}\label{s:contributions}

The considered inverse problem \emph{poses significant computational and mathematical challenges} that render traditional methods for the solution of inverse problems inefficient. Our work extends~\cite{rudi2022:parameter} in several ways. In~\cite{rudi2022:parameter}, a NN was used to infer two parameters of the dynamical system (namely, $\theta_0$ and $\theta_1$ in \cref{e:fhnfwd}) as well as noise parameters of an autocorrelated additive noise model. Here, we consider a second noise model based on a stochastic differential equation (\iacr{SDE})~\cite{leon2018:hypoelliptic,samson2025:inference, buckwar2022:splitting}, estimate one additional parameter of the dynamical system, and provide methodology for uncertainty quantification. In addition, we explore the performance of the proposed approach on dedicated hardware. In particular, our contributions are:
\begin{itemize}[leftmargin=*]
%//////////////////////////////////////
\item We computationally learn an observation-to-parameter map $g_{\Xi} : \mathcal{Y} \to \Theta$ that is able to simultaneously estimate the parameters of the dynamical system $\vect{\theta}_{\text{dyn}} \in \mathbb{R}^3$, noise parameters $\vect{\theta}_{\text{noise}} \in \mathbb{R}^m$, $m \in \{1,2,3\}$, as well as the covariance matrix $\mat{\Gamma}_{\text{post}} \in \spd{3}$ of the posterior distribution from noisy observations. In~\cite{rudi2022:parameter} only two model and noise parameters were considered.
%//////////////////////////////////////
\item We consider the estimation of model parameters and noise parameters for different noise models---in particular, autocorrelated additive noise and noise modelled via SDEs, and a combination thereof. We refer to the SDE noise model as ``intrinsic noise.'' In~\cite{rudi2022:parameter} only additive noise was considered.
%//////////////////////////////////////
\item We develop a computational framework for training and estimating covariance matrices of the posterior distribution of the model parameters (along with model and noise parameters) to expose uncertainties in the model and data. This framework is absent from~\cite{rudi2022:parameter}.
%//////////////////////////////////////
\item We consider a log-Euclidean framework to learn $\mat{\Gamma}_{\text{post}}$.
%//////////////////////////////////////
\item We investigate different NN architectures for the map $g_{\Xi}$.
%//////////////////////////////////////
\item We provide a detailed performance analysis of the proposed methodology. In particular, we quantify prediction accuracy as a function of varying network architecture, NN training parameters, input data, and noise perturbations.
%//////////////////////////////////////
\item We report results for computational performance on dedicated hardware, in particular, the Neocortex System at the Pittsburgh Supercomputing Center, a testbed supercomputer consisting of the Cerebras CS-2 AI accelerator hardware.
\end{itemize}

\subsection{Limitations}
The limitations of our work include:
\begin{itemize}[leftmargin=*]
%//////////////////////////////////////
\item The observation-to-parameter map $g$ is not well-defined due to the ill-posedness of the inverse problem. In the absence of noise, it is reasonable to assume that the unpolluted observation lies in the range of $f$. However, since only polluted data is available this observation may not hold.  In addition, the map $f$ will have a non-trivial kernel; multiple instances of $\vect{\theta}_{\text{dyn}}$ will correspond to indistinguishable data $y_{\text{obs}}$. As such, we can interpret $g_{\Xi}$ as a surrogate for a pseudo-inverse $f^\dag$ of $f$.
%//////////////////////////////////////
\item We consider an ``end-to-end'' learning approach~\cite{Ongie:2020a,Lin:2018a, rudi2022:parameter}; we learn a direct map from observables $y_{\text{obs}} \in \mathcal{Y}$ to parameters $\vect{\theta} \in \Theta$. Therefore, our predictions (of, e.g., $\vect{\theta}_{\text{dyn}}$) are no longer informed by a likelihood~\cite{Chung:2024a}. However, experimental evidence (see supplementary material \Cref{s:predicted_state_supp}) seems to suggest that the predictions are for the most part good agreement with the data, especially when an additive noise model is considered. We also do not require imposing a metric on $\mathcal{Y}$, which turns out to be a main source of concern in neural dynamics~\cite{Prinz:2004a}.
%//////////////////////////////////////
\item We use an explicit time integration method for evaluating the forward and adjoint operators. Exploring more sophisticated methods (for example, adaptive time stepping; higher accuracy; etc.) remains subject to future work.
%//////////////////////////////////////
\item The amount of samples for training NNs can increase significantly with the dimension of the parameter space provided sufficiently uninformative priors; and each sample incurs the cost of solving the ODE model. If the dimension of the parameter space becomes ``too large,'' the cost of data generation can become prohibitive.
\end{itemize}

\subsection{Related Work}\label{s:review}

In~\cite{rudi2022:parameter}, we estimated a subset of the parameters of the dynamical system (namely, $\theta_0$ and $\theta_1$ in \cref{e:fhnfwd}) and the two parameters of an autocorrelated, additive noise model. Here, we generalize our approach to include a noise model based on an SDE and we estimate the parameters of the SDE. The SDE introduces stochastic dynamics to the deterministic ODE; and inference for SDE parameters is significantly more challenging than inference for parameters of an additive noise model. Moreover, we extend our framework to enable the prediction of covariance matrices associated with the posterior distribution of the model parameters to expose model and data uncertainties. Overall, we estimate up to 12 parameters associated with our problem. While the parameters of the dynamical system can be estimated based on global (e.g., stochastic) optimization methods, the overall problem tackled in this work is more complex; it also involves estimating parameters for the noise model and entries of the covariance matrix associated with the posterior distribution of the parameters of the dynamical system conditioned on the data. The simultaneous estimation of all of those quantities of interest---ODE parameters, SDE parameters, stochastic noise parameters, and posterior covariance---has not been demonstrated to our knowledge.
Modeling neural activity has a rich mathematical history~\cite{Chialva:2023a}. So does parameter estimation for models of neural dynamics~\cite{VanGeit:2008a}. In the following, we highlight work that we view most pertinent to ours. The dynamical system that governs the parameter estimation problem is the \fhn\ model~\cite{Fitzhugh:1961a, nagumo1962:active}. We describe this model in greater detail in \Cref{s:methods}. It represents a common simplification of the Hodgkin--Huxley model to a nonlinear, coupled system of two ODEs. Traditional approaches for estimating parameters for these types of models are based on heuristics, follow trial and error strategies, and/or global optimization approaches. Examples for such approaches include simulated annealing, differential evolution, genetic algorithms, particle-swarm optimization, and brute-force grid searches~\cite{buhry2008:parameter, buhry2011:automated, VanGeit:2008a,Abdelaty:2022a,Hartoyo:2019a,druckmann2007:novel, gurkiewicz2007:numerical}. While some of these methods are, in general, equipped with guarantees to converge to a global minimizer, their rate of convergence is typically poor. Conversely, methods that exploit local gradient information (e.g., see~\cite{Doi:2002a,Toth:2011a,Meliza:2014a, ramsay2007:parameter}) suffer from strong nonlinearities and the nonconvexity of the problem.

An additional shortcoming of the deterministic methods mentioned above is that they only provide point estimates for the parameter values but no information about uncertainties or error bounds for the inferred parameters. Work on uncertainty quantification for computational models in neuroscience can be found in~\cite{Tennoe:2018a,Wang:2022a}. The work in~\cite{Tennoe:2018a} uses polynomial chaos expansions instead of more traditional Markov chain Monte Carlo (\iacr{MCMC}) sampling approaches (e.g., considered in \cite{Wang:2022a,Kostuk:2012a}). Other probabilistic approaches that have been used include techniques such as state-space modeling, Kalman filters and variants thereof~\cite{Huang:2022a,Ernst:2014a,Campbell:2020a, Vavoulis:2012s}.

Generally, the solution of inverse problems using deep NNs has recently gained attraction~\cite{Adler:2017a,Li:2020a,Ongie:2020a,Mukherjee:2023a,Genzel:2022a}. For a recent overview, see \cite{Bach:2024}.

Our proposed approach falls into the broad class of likelihood-free inference (\iacr{LFI}) or simulation-based inference (\iacr{SBI})~\cite{cranmer2020:frontier, lueckmann2021:benchmarking, bahl2024:advancing} methods in the sense that we learn a surrogate model that maps from observations to hidden parameters and summary statistics without explicitly evaluating, manipulating, and/or forming a posterior or likelihood. This class of algorithms has been developed for scenarios in which the likelihood is intractable or unavailable but it is possible to generate data from the model (i.e., evaluate the forward map $e(\vect{\theta}_{\text{dyn}})$). It subsumes a large class of approaches. In our framework, we directly learn a surrogate map from observations to quantities of interest. This is different from many existing LFI/SBI algorithms, where expensive likelihood evaluations are avoided and informative simulations are selected to increase sampling efficiency~\cite{csillery2010:approximate, sisson2018:handbook, sisson2007:sequential}. Moreover, there are primarily three NN-based approaches to SBI: Neural Posterior Estimation (NPE)~\cite{Papamakarios:2016a, greenberg:2019}, Neural Likelihood Estimation (NLE)~\cite{papamakarios:2019, glockler:2022}, and Neural Ratio Estimation (NRE)~\cite{hermans:2020, durkan:2020}. These approaches are conceptually different from ours. We outline the differences to the potentially closest approach of the three, which is amortized NPE. In amortized NPE, a NN is trained once to learn the mapping from observations to a posterior. The NN employed by NPE typically is a normalizing flow~\cite{Rezende:2015}. The differences between our approach and NPE are that, first, we are free to choose any data misfit function, which takes the role of the negative-log likelihood function in the posterior.  Second, we approximate posteriors with Gaussian distributions with the benefit of a reduction in computational cost, because we obtain the covariance of the Gaussian using the inverse Hessian matrix and adjoint methods; while NPE allows for non-Gaussian posteriors.  Third, we train feed-forward-type NN that have a lower computational cost in practice compared with normalizing flows. The aforementioned NN-based approaches to SBI are primarily either concerned with inference from output of a deterministic model (with some amount of added noise) or from output of a stochastic model; representative benchmark problems can be found in~\cite{lueckmann2021:benchmarking}. We go beyond this deterministic/stochastic separation, because we target the simultaneous inference of both ODE/SDE and noise parameters with just one single NN. The noise models are time-correlated and not merely scaled Gaussian white noise. The levels of noise vary across a range of parameters that go into stochastic models; therefore, our treatment of noise presents additional challenges that are not considered in the literature.

Machine learning methods for inverse problems governed by dynamical systems are discussed in, for example~\cite{Raissi:2019a,Court:2023a,Pakravan:2021a,nguyen2024tnet,lye2021iterative}. In the broader context of inverse problems, end-to-end approaches that directly learn the reconstruction map $g_{\Xi} : \mathcal{Y} \to \Theta$---as done in our work---are presented in~\cite{Rick:2017a,Adler:2017a,Jin:2017a,Wurfl:2016a,Antholzer:2019a,Chlemper:2017a,Lenzi:2023a,Sainsbury:2024a}. Instead of learning the reconstruction map, one can also learn the surrogates for the solution operators for the dynamical systems to speed up the optimization~\cite{lye2021iterative,wang2021fast}. Alternatively, one can use NNs to compute the gradient for updating the hidden parameters~\cite{Chen:2018a}. Other approaches to learn the reconstruction map use two autoencoders---one for representing the space $\mathcal{Y}$ and one for representing the space $\Theta$. Works that consider this type of pairing of autoencoders are~\cite{Chung:2024a,Feng:2023a,Zeng:2015a}. Another strategy is to learn a model representation without considering data using diffusion networks~\cite{Cao:2024a}, generative adversarial networks~\cite{Shah:2018a}, or variational autoencoders~\cite{Goh:2022a}.

Classical (feed-forward) NNs do not allow to generate predictions of uncertainties; they only provide point estimates. Thus using such a NN for parameter prediction results in a point estimate, which is similar to solving an inverse problem with traditional optimization methods. In the present work, we use feed-forward NNs and consequently are only able to generate point estimates. However, we propose to learn the mapping from observational data to covariance matrices approximating the posterior distribution. As an alternative, generative NN architectures map samples (e.g., from a standard normal distribution) to samples from an unknown target distribution (e.g., a posterior). These NNs have been used to approximate the underlying likelihood density or posterior distribution. In practice, the computational cost for training the weights of generative NNs is significantly higher than for classical deep NNs. Hence, the present work proposes a computationally less demanding approach for prediction of parameter uncertainties.

The aforementioned generative NNs for approximating the underlying likelihood density or posterior distribution from samples are normalizing flow architectures~\cite{Rezende:2015,Papamakarios:2016a,Lueckmann:2017a,Papamakarios:2019a,Dasgupta:2024}, which can be applied to neural dynamics~\cite{Gonccalves:2020a}, generative adversarial networks~\cite{Patel:2022a,Ray:2023}, and score-based diffusion NNs~\cite{Batzolis:2021,Dasgupta:2025}. Alternatively, one can learn the forward operator to, for instance, speed up the exploration of the posterior distribution~\cite{Antil:2023a}.

Machine learning has not only been used to learn the solution of the inverse problem given some observational data but has also been successfully used to estimate hyperparameters that appear in the formulation of inverse problems~\cite{Kofler:2023a,Byrne:2023a,Afkham:2021a} or learn/improve regularization operators~\cite{Court:2023a,Alberti:2021a,Boink:2023a,Cho:2021a,Obmann:2020a,Mukherjee:2021a,Schwab:2019a,Li:2020a,Altekruger:2023a}. Theoretical results associated with the use of machine learning for inference in inverse problems can be found in~\cite{Puthawala:2022a,Castro:2022a}. Alternatively, machine learning has been used to accelerate optimization~\cite{Chen:2022a} and accelerate MCMC~\cite{Gabrie:2022}.

\subsection{Outline}

The paper is organized as follows. In \Cref{s:methods} we introduce the methodology. The general problem formulation is discussed in \Cref{s:formulation}. This includes the variational problem formulation in the deterministic setting (see \Cref{s:variational}) and the formulation of the inverse problem in a Bayesian setting (see~\Cref{s:bayesianinference}). This allows us to derive problem specific information used in our NN framework, explore analogies between the NN approach and classical approaches, and motivate our setup of the training for the NN architectures. We present the considered NN architectures in \Cref{s:architectures}. This includes dense NN (\iacr{DNN}) architectures (see \Cref{s:dnn}) and convolutional NNs (\iacr{CNN}s; see \Cref{s:cnn}). A description of the training and testing data can be found in \Cref{s:data}. We present the measures to assess the performance of the proposed approach in \Cref{s:performance}. The results are reported in \Cref{s:results}. We conclude with \Cref{s:conclusions}.

We accompany the manuscript with supplementary materials. These include the derivation of the Hessian (see \Cref{s:curvature_supp}) and additional numerical results that are either complementary to those reported in this study, provide additional experimental evidence for some of the claims we make, or explore variations of our framework (see \Cref{s:supplementary_results}).

\section{Methods and Material}\label{s:methods}

This section presents $(i)$ classical variational and Bayesian problem formulations to provide a theoretical underpinning for the choices and strategies we consider in our NN approach, $(ii)$ the proposed and examined NN architectures for inferring the parameters of the dynamical system, noise parameters and entries of the covariance matrix of the posterior distribution, $(iii)$ the training and testing data for the NN, and $(iv)$ performance measures to assess the NN predictions. We summarize the notation in \Cref{t:notation}.

\begin{table}
\caption{Common notation, symbols, and acronyms.\label{t:notation}}
\centering\adjtab
\renewcommand{\arraystretch}{1.2}
\begin{tabular}{ll}\toprule
\bf symbol/acronym & \bf meaning \\\midrule
$[0,\tau]\subset\ns{R}$ & time horizon\\
$\mathcal{U} \times \mathcal{V}$ & state space\\
$\Theta_{\text{dyn}}$ & control space\\
$\mathcal{Y}$ & observation space\\
$\Theta$ & quantities of interest (NN output space)\\
$\spd{n}$ & space of real $n \times n$ symmetric positive definite matrices \\
$\sym{n}$ & vector space of real $n \times n$ symmetric matrices \\
$u : [0,\tau] \to \ns{R}$ & membrane potential (state variable) \\
$v : [0,\tau] \to \ns{R}$ & recovery variable (state variable) \\
$z \in \ns{R}$ & stimulus \\
$\lambda : [0,\tau] \to \ns{R}$ & dual variable associated with state equation for $u$ \\
$\nu : [0,\tau] \to \ns{R}$ & dual variable associated with state equation for $v$ \\
$\vect{\theta}_{\text{dyn}} \defeq (\theta_0,\theta_1,\theta_2)\in\mathbb{R}^3$ & model parameters (dynamical system) \\
$\vect{\theta}_{\text{noise}} \defeq (\delta,\sigma,\beta)\in\mathbb{R}^3$ & noise parameters \\
$(\delta,\sigma) \in \mathbb{R}^2$ & noise parameters of additive noise model \\
$\beta \in \mathbb{R}$ & noise parameters of intrinsic noise model \\
$p\in\mathbb{N}$ & number of unknowns; $p \in \{4,5,6,10,11,12\}$ \\
$f : \Theta_{\text{dyn}} \to \mathcal{Y}$ & parameter-to-observation map\\
$e : \Theta_{\text{dyn}} \to \mathcal{U}$ & forward map\\
$o : \mathcal{U} \to \mathcal{Y}$ & observation operator\\
$g_{\text{dyn}} : \mathcal{Y} \to \Theta_{\text{dyn}}$ & reconstruction map; ``pseudo-inverse'' of $f$\\
$g_{\Xi} : \mathcal{Y} \to \Theta $ & NN surrogate parametrized by $\Xi$ \\
$y_{\text{obs}} \in \mathcal{Y}$ & observational data (NN input)\\
$\vect{y}_{\text{obs}} = (y_{\text{obs},1}, \ldots, y_{\text{obs},n_t}) \in \ns{R}^{n_t}$ & discrete observational data \\
$\vect{\theta}_{\text{pred}}\in\mathbb{R}^p$ & NN prediction (NN output) \\
$\mat{\Gamma}_{\text{prior}} = \diag(\vect{\sigma}_{\text{prior}}) \in \mathbb{R}^{3,3}$ & prior covariance \\
$\vect{\bar{\theta}}_{\text{prior}}\in\mathbb{R}^3$ & prior mean; $\sigma_{\text{prior}} \in \ns{R}^3$ \\
$\mat{L} = \mat{\Gamma}_{\text{prior}}^{-1} \in \mathbb{R}^{3,3}$ & regularization operator \\
$\mat{\Gamma}_{\text{noise}} = \sigma_{\text{noise}}^2\tau^2 \mat{I} \in \mathbb{R}^{3,3}$ & noise covariance; $\sigma_{\text{noise}} > 0$, $\tau > 0$, $\mat{I} = \diag(1,1,1) \in \mathbb{R}^{3,3}$  \\
$\mat{\Gamma}_{\text{post}} \in \spd{3}$  & posterior covariance with entires $\Gamma_{ij}$, $i=0,1,2$\\
$\mat{H} \in \spd{3}$  & Hessian matrix\\
$n_t \in \ns{N}$ & number of time steps (numerical time integration) \\
$n_{\text{layers}} \in \ns{N}$ & number of NN layers \\
$n_{\text{train}} \in \ns{N}$ & number of training data \\
$n_{\text{test}} \in \ns{N}$ & number of testing data \\
$n_{\text{feat}} \in \ns{N}$ & number of features \\
$n_f \in \ns{N}$ & filter count (CNN) \\
$n_u \in \ns{N}$ & number of units (DNN) \\
\midrule
CDET    & Coefficient of Determination; see \cref{e:cdet}\\
(C)MSE  & (Centered)Mean Squared Error; see \cref{e:cmse}\\
(C/D)NN & (Convolutional/Dense) Neural Network \\
MAP     & Maximum-A-Posteriori (point) \\
MdAPE   & Median of Absolute Percentage Error; see \cref{e:mdape}\\
MCMC    & Markov Chain Monte Carlo (method) \\
ODE     & Ordinary Differential Equation \\
SDE     & Stochastic Differential Equation \\
SQB     & Squared Bias (error); see \cref{e:sqb} \\
\ts     & time series features (NN input)\\
\fc     & Fourier coefficient features (NN input)\\
\tsfc   & time series and Fourier coefficient features (NN input)\\
\bottomrule
\end{tabular}
\end{table}

\subsection{Classical Problem Formulation}\label{s:formulation}

Next, we introduce the formulation of the inverse problem for estimating the parameters of the dynamical system in \cref{e:fhnfwd} from time series data $y_{\text{obs}}$ in a variational setting~\cite{Troltzsch:2010a, tarantola2005:inverse}. Subsequently, we develop the associated Bayesian formulation of the inverse problem~\cite{tenorio2017:introduction,Stuart:2010a,Calvetti:2023a}. We include this information for several reasons: First, we use it to illustrate mathematical challenges. Second, we use analogies to the statistical problem formulation to interpret the results obtained by the NNs. Third, we use curvature information of the variational problem formulation to generate training data for the estimation of the posterior covariances $\mat{\Gamma}_{\text{post}} \in \spd{3}$, with
\[
\spd{n} \defeq \left\{ \mat{M} \in \mathbb{R}^{n,n} : \vect{x}^\mathsf{T} \mat{M} \vect{x} > 0\,\,\text{for all}\,\,\vect{x} \in \mathbb{R}^n \setminus \{\vect{0}\} \right\}.
\]

\subsubsection{Forward Problem}\label{s:forward}

The forward problem is given by the nonlinear system of ODEs
\begin{subequations}\label{e:fhnfwd}
\begin{align}
\d_t u - \theta_2 (u - (u^3/3) + v + z ) & = 0
&& \text{in}\;\;(0,\tau],
\label{e:stateu}
\\
\d_t v + ((u-\theta_{0} + \theta_{1} v)/\theta_2) & = 0
&& \text{in}\;\;(0,\tau],
\label{e:statev}
\end{align}
\end{subequations}

\noindent with model parameters $\vect{\theta}_{\text{dyn}} \defeq (\theta_0,\theta_1,\theta_2) \in \ns{R}^3$, initial conditions $u(t=0) = u_0$, $v(t=0) = v_0$, and state variables $u \in \mathcal{U} \subset \{ w : [0,\tau] \to \mathbb{R}\}$, $v \in \mathcal{V} \subset \{ w : [0,\tau] \to \mathbb{R}\}$, respectively, defined on a given time horizon $[0,\tau]\subset\ns{R}$, $\tau > 0$. The state variable $u$ is related to the \emph{membrane potential} across an axon membrane, $v$ represents the \emph{recovery variable} summarizing outward currents, and $z$ is the \emph{total membrane potential} and corresponds to a stimulus applied to the neuron. We assume that the stimulus $z$ is a known (fixed) constant. We show representative membrane potentials $u$ as a function of $\vect{\theta}_{\text{dyn}}$ in the supplementary materials (see \Cref{f:fwdsim-range} in \Cref{s:predicted_state_supp}).

This model is a phenomenological simplification of more complex neural models to capture qualitative dynamics. As such, the model parameters do not really have a physiological meaning. The parameters $\theta_0$ and $\theta_1$ dominate the dynamics.  The parameter $\theta_0$ denotes the excitability threshold; it determines the threshold above which the system enters an excited state, influencing spiking behavior. The parameter $\theta_1$ represents the recovery scaling; it controls how the recovery variable responds to voltage changes. The parameter $\theta_2$ is a time scale factor; it controls the time scale separation between fast and slow dynamics.

\subsubsection{Variational Problem Formulation}\label{s:variational}

We overview traditional, variational approaches to estimate the model parameters $\vect{\theta}_{\text{dyn}} \in \ns{R}^3$ that enter \cref{e:fhnfwd} given noisy observations $y_{\text{obs}} \in \mathcal{Y}$. The problem formulation is given by\
\begin{equation}\label{e:varopt}
\begin{aligned}
\minopt_{(u,v) \in \mathcal{U} \times \mathcal{V},\,\vect{\theta}_{\text{dyn}} \in \Theta} &
\quad
\frac{\gamma}{2}\int_0^\tau \!\!(u(t) - y_{\text{obs}}(t))^2 \,\d t
+ \frac{1}{2}\|\vect{\theta}_{\text{dyn}} - \vect{\theta}_{\text{ref}}\|^2_{\mat{L}}
\\
\text{subject to}\quad\quad\;\; &
\quad
(u,v) \text{ solve }\cref{e:fhnfwd},
\end{aligned}
\end{equation}

\noindent The first part of the objective function in \cref{e:varopt} measures the proximity between $u$ and $y_{\text{obs}}$, where $\gamma >0$ controls its contribution to the objective. We assume that $y_{\text{obs}}$ is available at all time points $t \in [0,\tau]$. If we assume partial observations $y_{\text{obs}} \in \mathbb{R}^{n_{\text{obs}}}$, $n_{\text{obs}} \in \mathbb{N}$, a distance function that is in accordance with \cref{e:varopt} can, e.g., be chosen as $(\gamma/2) \| o(u) - y_{\text{obs}}\|_2^2$, with observation operator $o : \mathcal{U} \to \mathbb{R}^{n_{\text{obs}}}$. We select $y_{\text{obs}}$ to correspond to the membrane potential $u$ perturbed by noise: $y_{\text{obs}}(t) = u_{\vect{\theta}}(t)$, $\vect{\theta} = (\vect{\theta}_{\text{dyn}},\vect{\theta}_{\text{noise}})$, for all $t \in [0,\tau]$, where $\vect{\theta}_{\text{noise}} \in \mathbb{R}^m$, $m \in \{1,2,3\}$, controls the amount of noise. We describe the considered noise models in \Cref{s:noisemodel}. The second part of the objective function is a Tikhonov-type regularization model that is introduced to stabilize the solution process~\cite{Engl:1996a}, where $\mat{L} \in \spd{3}$ is a regularization operator that is specified below and $\vect{\theta}_{\text{ref}} \in \ns{R}^3$ is a reference vector added to control the kernel of the regularization.

The reduced form of the optimization problem in~\cref{e:varopt} is given by
\begin{equation}\label{e:varoptred}
\minopt_{\vect{\theta}_{\text{dyn}} \in \Theta_{\text{ad}}}
\quad
\frac{\gamma}{2}\int_0^\tau (f(\vect{\theta}_{\text{dyn}}) - y_{\text{obs}}(t))^2 \,\d t
+ \frac{1}{2}\|\vect{\theta}_{\text{dyn}} - \vect{\theta}_{\text{ref}}\|^2_{\mat{L}}.
\end{equation}

\noindent In \cref{e:varoptred} we eliminated the constraint from the problem by replacing the model prediction with the parameter-to-observation map $f : \Theta \to \mathcal{Y}$. This formulation establish a connection to the problem formulation in a statistical setting presented next.

\subsubsection{Bayesian Problem Formulation}\label{s:bayesianinference}

We refer to~\cite{tenorio2017:introduction,Stuart:2010a,Calvetti:2023a} for an introduction to Bayesian inverse problems. With slight abuse of notation, we assume that the problem has been discretized. Consequently, the observations $y_{\text{obs}} \in \mathcal{Y}$ are given by $\vect{y}_{\text{obs}} = (y_{\text{obs},1}, \ldots, y_{\text{obs},n_t}) \in \ns{R}^{n_t}$, where $n_t \in \ns{N}$ corresponds to the number of time steps used to solve the forward problem. We model all unknown variables as random variables.

For the likelihood model we assume that observational uncertainty (i.e., uncertainty in the data related to measurement errors) and the uncertainty in the model (i.e., the uncertainty in $\vect{f}(\vect{\theta}_{\text{dyn}})$ due to modeling errors) are each centered, additive, and Gaussian. We integrate both errors into a single noise model. More precisely, we assume $\vect{y}_{\text{obs}} \sim \mathcal{N}(\vect{f}(\vect{\theta}_{\text{dyn}}), \mat{\Gamma}_{\text{noise}})$, with parameter-to-observation map $\vect{f} : \ns{R}^3 \to \ns{R}^{n_t}$ and covariance $\mat{\Gamma}_{\text{noise}} \in \spd{n_t}$, $\mat{\Gamma}_{\text{noise}} = \sigma_{\text{noise}}^2\tau^2\mat{I}$, $\mat{I} = \operatorname{diag}(1,\ldots,1) \in \ns{R}^{n_t \times n_t}$, $\sigma_{\text{noise}} > 0$. We obtain
\begin{equation}\label{e:yobs}
\vect{y}_{\text{obs}} = \vect{f}(\vect{\theta}_{\text{dyn}}) + \vect{\eta}
\end{equation}

\noindent with $\vect{\eta} \sim \mathcal{N}(\vect{0}, \mat{\Gamma}_{\text{noise}})$. Under the above assumptions, the likelihood probability density function is given by
\begin{equation}
\label{e:like}
\pi_{\text{like}}(\vect{y}_{\text{obs}}\,|\, \vect{\theta}_{\text{dyn}})
\propto
\exp\left(
-\frac{1}{2}
\left\langle
\vect{y}_{\text{obs}} - \vect{f}(\vect{\theta}_{\text{dyn}}), \mat{\Gamma}_{\text{noise}}^{-1}(\vect{y}_{\text{obs}} - \vect{f}(\vect{\theta}_{\text{dyn}}))
\right\rangle
\right).
\end{equation}

We also have to decide on a prior distribution for $\vect{\theta}_{\text{dyn}}$. We model each component of $\vect{\theta}_{\text{dyn}}$ as normal and independent identically distributed random variable; $\vect{\theta}_{\text{dyn}} \sim \mathcal{N}(\vect{\bar{\theta}}_{\text{prior}}, \mat{\Gamma}_{\text{prior}})$, $\vect{\bar{\theta}}_{\text{prior}} \in \ns{R}^3$, $\mat{\Gamma}_{\text{prior}} = \diag(\vect{\sigma}_{\text{prior}}) \in \spd{3}$, with standard deviation $\vect{\sigma}_{\text{prior}} \in \mathbb{R}^3$. The associated probability density function is given by
\begin{equation}
\label{e:prior}
\pi_{\text{prior}}(\vect{\theta}_{\text{dyn}})
\propto
\exp\left(-\frac{1}{2}
\left\langle
\vect{\theta}_{\text{dyn}} - \vect{\bar{\theta}}_{\text{prior}},
\mat{\Gamma}_{\text{prior}}^{-1}
(\vect{\theta}_{\text{dyn}} - \vect{\bar{\theta}}_{\text{prior}})
\right\rangle\right).
\end{equation}

\begin{remark}
We use the prior distribution to draw samples for the training of the NNs. The selection of the hyper-parameters $\mat{\Gamma}_{\text{prior}}^{-1}$ and $\vect{\bar{\theta}}_{\text{prior}}$ of $\pi_{\text{prior}}$ in \cref{e:prior} is a challenge in itself. In a Bayesian setting, one can treat these hyperparameters as unknown random variables and estimate them along with the parameters $\vect{\theta}_{\text{dyn}}$. This leads to so-called \emph{hierarchical methods} (or \emph{hyperprior models})~\cite{Kaipio:2006a, tenorio2017:introduction}. Since we know little about the range of the parameters $\vect{\theta}_{\text{dyn}}$, we use a rather uninformative prior distribution with loose bounds on $\vect{\sigma}_{\text{prior}}$. More details regarding the parameter selection are provided in \Cref{s:parameters}.

We use a Gaussian prior distribution to draw samples for the parameters $\vect{\theta}_{\text{dyn}}$. Existing work that considers a Gaussian prior for the parameters for the \fhn\ model includes~\cite{jensen2012:markov}. The work in~\cite{samson2025:inference} uses uniform distributions. Since we use the samples drawn from the prior distribution solely to generate NN inputs (i.e., evaluate the forward operator $e$), our framework is, in general, amendable to any type of prior distribution.

We assume Gaussian distributions when approximating the posterior at each (known) ``true solution'' $\vect{\theta}_{\text{dyn}}^i$ (a sample chosen from the prior). The mean of the Gaussian is $\vect{\theta}_{\text{dyn}}$ and the covariance is computed from the Hessian of the negative log-posterior.  This quadratic approximation (or linearization) is adequate, since we evaluate the model during data generation at the ``true solution.''
\end{remark}

According to Bayes theorem, the posterior distribution $\pi_{\text{post}}$ of $\vect{\theta}_{\text{dyn}}$ conditioned on the data $\vect{y}_{\text{obs}}$ is given by
\begin{equation}
\label{e:posterior}
\pi_{\text{post}}(\vect{\theta}_{\text{dyn}}\,|\,\vect{y}_{\text{obs}})
\propto
\pi_{\text{like}}(\vect{y}_{\text{obs}}\,|\,\vect{\theta}_{\text{dyn}})\,\,
\pi_{\text{prior}}(\vect{\theta}_{\text{dyn}}).
\end{equation}

The MAP point $\vect{\theta}_{\text{map}} \in \mathbb{R}^3$ can be found by minimizing the negative log posterior
\[
-\log \pi_{\text{post}}(\vect{\theta}_{\text{dyn}}\,|\,\vect{y}_{\text{obs}}) = \phi(\vect{\theta}_{\text{dyn}}) + \text{const},
\]

\noindent with $\phi : \ns{R}^3 \to \ns{R}$,
\begin{equation}\label{e:neglogpost}
\phi(\vect{\theta}_{\text{dyn}}) =
\frac{1}{2}
\|
\vect{y}_{\text{obs}} - \vect{f}(\vect{\theta}_{\text{dyn}})
\|^2_{\mat{\Gamma}_{\text{noise}}^{-1}}
+
\frac{1}{2}
\|
\vect{\theta}_{\text{dyn}} - \vect{\bar{\theta}}_{\text{prior}}
\|_{\mat{\Gamma}_{\text{prior}}^{-1}}^2.
\end{equation}

The function $\phi$ is a discrete version of the objective function in~\cref{e:varoptred} with $\gamma = 1/(\sigma^2_{\text{noise}}\tau^2)$, precision matrix $\mat{L} = \mat{\Gamma}_{\text{prior}}^{-1}$ and $\vect{\theta}_{\text{ref}} = \vect{\bar{\theta}}_{\text{prior}}$. This establishes a direct connection to the deterministic optimization problem in~\cref{e:varoptred}. We use this connection in \Cref{s:data} to derive an approximation of the covariance matrix $\mat{\Gamma}_{\text{post}}\in \spd{3}$ of the posterior distribution in~\cref{e:posterior}.

\begin{remark}
The proposed \emph{NN-based approach} discussed next is not guaranteed to estimate the MAP point $\vect{\theta}_{\text{map}}$ or any other estimates such as the maximum likelihood or conditional mean.

We also note that we cannot \emph{a priori} guarantee that our NN predictor is unbiased. According to classical ML theory, there exists a bias-variance tradeoff that is controlled by the complexity (or capacity) of the ML model~\cite{geman1992:neural}. If the ML model has a strong bias it is not sufficiently complex; we expect to observe poor performance during the training, validation and testing phase (\emph{underfitting}). Conversely, increasing model complexity will increase the variance of the predictor, leading to an improved performance during training but a decrease in performance on the validation and testing data (\emph{overfitting}). Overall, large bias or variance of the ML predictor will lead to a large generalization error. Other sources of generalization error include insufficient or biased training data, a poor validation strategy, or instabilities in the optimization strategy (e.g., poor convergence, bad local minima). Consequently, in the absence of any theoretical guarantees, we have to perform an empirical study to assess the generalization error on validation and testing data; we do so in the experimental part of this manuscript.

We emphasize that there exists recent evidence that modern learning architectures can escape the classical bias-variance trade-off; bias and variance can decrease with increasing network complexity~\cite{belkin2019:reconciling, neal2018:modern, yang2020:rethinking, yu2025:near}.
\end{remark}

\subsection{Neural Network Architectures}\label{s:architectures}

Next, we describe the NN architectures. As we have mentioned in \Cref{s:method-outline}, we want to learn the \emph{reconstruction map} $g_{\Xi} : \mathcal{Y} \to \Theta$. Our approach is supervised. The architectures are in line with our past work~\cite{rudi2022:parameter}. Consequently, we keep their discussion brief.

The inputs (features) $\vect{y}_{\text{obs}}$ to the NNs include, e.g., time series that correspond to the solution of the ODE in~\cref{e:fhnfwd} perturbed by noise. We specify the training and testing data as well as the noise models considered in this work in~\Cref{s:data}. The output of the NNs are predictions $\vect{\theta}_{\text{pred}} \in \ns{R}^p$ of hidden parameters for a given noisy observation $\vect{y}_{\text{obs}} \in \ns{R}^{n_t}$. We define the NN-based reconstruction map by
\begin{equation}\label{e:reconmap}
\vect{\theta}_{\text{pred}} = \vect{y}_n,
\quad \vect{y}_l = \vect{g}_l(\vect{y}_{l-1})
\quad \text{for}\,\, 1 \leq l \leq n_{\text{layers}},
\quad \vect{y}_0 = \vect{y}_{\text{obs}},
\end{equation}

\noindent where $\vect{g}_l$ denotes the $l$th layer of the NN. In this sense, we treat the reconstruction map as a sequence of recursive functions, where each function $\vect{g}_l$ corresponds to the output of the $l$th network layer.

We describe the different architectures considered for the layers $\vect{g}_l$ in~\cref{e:reconmap} in \Cref{s:dnn} and \Cref{s:cnn}, respectively. We summarize the NN parameters in \Cref{t:nnconfig}. Some of the settings that are shared across architectures are the loss function (mean squared error; \iacr{MSE}) and the ADAM optimization algorithm~\cite{Kingma:2014a}. We consider a learning rate of 0.002 (we explored other choices without noticeable benefit) and a SWISH activation function defined by
\begin{equation}\label{e:swish_activation}
\vect{h}_l(\vect{x}) = \vect{x} \oslash (\vect{e}_{n_l} + \vect{\exp}(-\vect{x})),
\quad
\vect{e}_{n_l} = (1,\ldots,1) \in \ns{R}^{n_l},
\end{equation}

\noindent for any $\vect{x} \in \ns{R}^{n_l}$. Here, $\oslash$ denotes the Hadamard division and $\vect{\exp}: \ns{R}^{n_l} \to \ns{R}^{n_l}$ is applied element-wise.

\begin{table}
\caption{DNN and CNN configurations. Entries that are sets denote different choices that we use in the exploration of NN architectures.}
\label{t:nnconfig}
\adjtab\centering
\begin{tabular}{lll}\toprule
\bf architecture & \bf option            & \bf value(s)              \\
\midrule
DNN \& CNN       & loss function         & MSE                       \\
                 & optimizer             & ADAM~\cite{Kingma:2014a}  \\
                 & learning rate         & 0.002                     \\
                 & batch size            & 32                        \\
                 & epochs                & 64                        \\
                 & activation function   & Swish/SiLU                \\
\midrule
DNN              & layers                & \{2,4,8,12,16,20\}        \\
                 & nodes/units $n_u$     & \{4,8,16,32,64,128,256\}  \\
\midrule
CNN              & layers                & \{2,3,4,5,6\}             \\
                 & filters setting $n_f$ & \{2,4,8,16,32,64,128\}    \\
                 & kernel size           & 3                         \\
                 & kernel stride         & 2                         \\
                 & pooling type          & average                   \\
                 & pooling size          & 2                         \\
                 & pooling stride        & 2                         \\
                 & post-flattening layers& 2 dense layers with $n_u=32$\\
\bottomrule
\end{tabular}
\end{table}

\subsubsection{Dense Neural Networks}\label{s:dnn}

This architecture consists of a sequence of fully connected layers. Each layer $l$ consists of $n_u \defeq n_u(l) \in \ns{N}$ units or nodes. Each $\vect{g}_l$ is modelled as the composition of an affine mapping with a nonlinear activation function: $\vect{g}_l(\vect{y}_{l-1}) = \vect{h}_l(\mat{W}_{\!l}\vect{y}_{l-1} + \vect{b}_l)$, where $\mat{W}_{\!l}$ denotes an $n_u(l) \times n_u(l-1)$  matrix of weights, and $\vect{b}_l$ is the $n_u(l) \times 1$ bias vector for the $l$th layer. The activation function $\vect{h}_l : \ns{R}^{n_u(l)} \to \ns{R}^{n_u(l)}$ is defined in~\cref{e:swish_activation}. We consider 4--24 layers and 4--256 units (see \Cref{t:nnconfig}) to cover a somewhat wide range of network sizes. This leads to a total of \inum{8079} up to \inum{2026243} trainable NN parameters. For experiments in \Cref{s:hpc_performance}, we utilize a DNN composed of 12 layers and 128 units per hidden layer, consisting of a total \inum{695179} trainable NN parameters.

\subsubsection{Convolutional Neural Networks}\label{s:cnn}

The second architecture are CNNs. This architecture allows us to exploit temporal locality of the considered time series based on convolutional layers. Each filter is applied to a subinterval of the time series. The filters of one layer are connected to all neighboring layers. This is similar to the units of the DNNs described in \Cref{s:dnn}. The weights of these connections represent the NN parameters to be learned during training. Since the weights are shared across time, significantly less weights have to be estimated. Pooling layers, which are nontrainable, are often used in CNNs for their reduction in spatial dimensionality.

\begin{figure}
\centering
\includegraphics[width=0.75\textwidth]{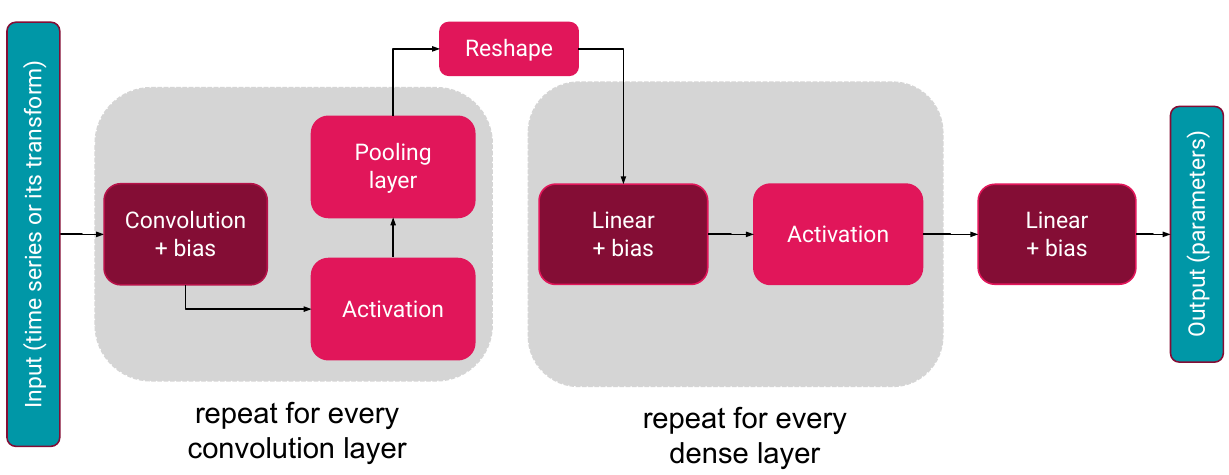}
\caption{Illustration of the CNN architecture. The input is passed into a sequence of convolutional block (each consisting of a convolution, activation, and pooling layers); the output of the last convolutional block is passed to a sequence of dense blocks (each consisting of a dense and activation layers). The last dense layer, without activation, predicts parameters.}
\label{f:cnn_diagram}
\end{figure}

The considered architecture is based on~\cite{Krizhevsky:2017a}. We interleave convolutional layers with pooling layers to aggregate a small block from a convolutional step into a single value. We determine the number of filters for the $l$th layer using a multiplier based on the depth (i.e., the layer number). That is, we select $n_f \times (1, 2, 4, 8, \ldots)$ filters. For example, for a CNN with $n_{\text{layers}}=3$ convolutional layers, we have $n_f \times 1$ filters in the first layer, $n_f \times 2$ filters in the second layer, $n_f \times 4$ filters in the third layer, and so on. Each convolutional layer reduces the size of the output compared with the input by a factor of $1/2$ (via a stride equals to two). The reduction by $1/2$ also holds for the pooling layers and is complemented by an increase in the filter count $n_f$ by a factor of 2. The output of these layers is passed to a sequence of two dense layers with $n_u = 32$ units. The diagram in \Cref{f:cnn_diagram} illustrates our CNN architecture. Across our experiments, the number of trainable parameters for this setup ranges between \inum{17095} and \inum{33662371}. For our main experiments, we utilize a CNN composed of 5 layers with $n_f=8$, reaching \inum{38516} trainable NN parameters.

\subsection{Training and Testing Data}\label{s:data}

Next, we describe how to generate the training and testing data. In analogy to the Bayesian problem formulation in \Cref{s:bayesianinference}, we generate the NN-outputs (labels) for the training and testing data
\begin{equation}\label{e:train-test-data}
\left\{\Theta_{\text{train}}, \Theta_{\text{test}}\right\}
=
\left\{
\{\vect{\theta}_{\text{train}}^i\}_{i=1}^{n_{\text{train}}},
\{\vect{\theta}_{\text{test}}^i\}_{i=1}^{n_{\text{test}}}
\right\}
\subset \mathbb{R}^p
\end{equation}

\noindent by drawing samples of the model parameters $\vect{\theta}_{\text{dyn}} \in \mathbb{R}^3$ from the prior distribution~\cref{e:prior}. The specific choices for $\vect{\sigma}_{\text{prior}}$ and $\vect{\bar{\theta}}_{\text{prior}}$ can be found in \Cref{s:parameters}. We discard any samples for $\vect{\theta}_{\text{dyn}}$ that are outside of predefined parameter bounds identified in~\cref{e:parabounds}. In addition, we also include noise parameters $\vect{\theta}_{\text{noise}}$ and approximations to the posterior covariances $\mat{\Gamma}_{\text{post}} \in \spd{3}$, $\mat{\Gamma}_{\text{post}} = (\Gamma_{ij})_{i,j=0,1,2}$, associated with $\vect{\theta}_{\text{dyn}}$ into the training and testing datasets, respectively. Since $\mat{\Gamma}_{\text{post}} \in \spd{3}$, we only consider the six upper triangular entries of $\mat{\Gamma}_{\text{post}}$ (see \Cref{s:logeuclidean} for more details). Overall, $p$ in \cref{e:train-test-data} ranges from 5 to 12, depending on the experiment. To generate the observations associated with the dataset $\left\{\Theta_{\text{train}}, \Theta_{\text{test}}\right\}$, we solve the forward problem in \cref{e:fhnfwd} and perturb the state variables according to the noise models specified in \Cref{s:noisemodel}. We obtain the NN-inputs (features)
\[
\left\{ Y_{\text{train}}, Y_{\text{test}} \right\}
=
\left\{
\left\{\vect{y}^i_{\text{train}}\right\}_{i=1}^{n_{\text{train}}},
 \left\{\vect{y}^i_{\text{test}}\right\}_{i=1}^{n_{\text{test}}}
\right\} \subset \mathbb{R}^{n_{\text{feat}}},
\]

\noindent corresponding to the training in testing dataset, respectively. We generate a total of \inum{15000} samples. We use fixed subsets (with randomly chosen entries) of these training and testing samples. As we increase the size of training data, we keep the samples we included in the smaller datasets.

If we only consider time series data as features, we have $n_{\text{feat}} = n_t$, where $n_t$ corresponds to the number of time steps used to integrate~\cref{e:fhnfwd} in time; we use an explicit, second-order accurate Runge--Kutta method. The time series are limited to $\tau = 200$\,ms. The number of time steps is set to $n_t = \inum{2000}$, which yields a time step size of $\tau/n_t = 0.1$\,ms. The identifier for these features is \iacr{\ts} (for ``time series''). Consequently, the examples $\vect{y}_{\text{train}}^i$ or $\vect{y}_{\text{test}}^i$ are the time series data associated with the $i$th sample of the model parameters $\vect{\theta}_{\text{dyn}}$ in $\Theta_{\text{train}}$ or $\Theta_{\text{test}}$, respectively. We also consider the associated $n_{\text{feat}} = n_t$ Fourier coefficients as features. We denote this feature setting by \iacr{\fc} (for ``Fourier coefficients''). Lastly, we consider a combination of these features, resulting in $n_{\text{feat}} = 2n_t$ features for each training or testing sample, respectively. The identifier for this setting is \iacr{\tsfc}.

Likewise to traditional methods for parameter estimation, we expect the performance of the proposed approach to deteriorate as a function of increasing noise. Using NNs for parameter estimation, a question that arises is how the presence of noise in the training data affects performance. To explore this, we consider the following settings:
\begin{inparaenum}[\it (i)]
\item training and testing data is noise-free (identifier: \nfnf),
\item training data is noise-free, testing data is noisy (identifier: \nfn),
and \item training data and testing data are noisy (identifier: \nn)\end{inparaenum}.

\begin{remark}
This question was already explored in~\cite{rudi2022:parameter}. Like in~\cite{rudi2022:parameter}, we observed that including noise in the training data yields better predictions for noisy observations. Although our setup is different from~\cite{rudi2022:parameter}, we decided to omit these experiments in the main manuscript and fix the setting to \nn. We report results in the supplementary materials that support this decision (see \Cref{s:effects_of_noise_training_supp}).
\end{remark}

Next, we detail how we generate training data for estimating the covariance matrix $\mat{\Gamma}_{\text{post}} \in \spd{3}$ of $\pi_{\text{post}}$ in~\cref{e:posterior}. This enables us to expose uncertainties as they propagate through our framework. One approach to construct the dataset $\{\{\mat{\Gamma}_{\text{train}}^i\}_{i=1}^{n_{\text{train}}},\{\mat{\Gamma}_{\text{test}}^i\}_{i=1}^{n_{\text{train}}}\}$ for $\mat{\Gamma}_{\text{post}}$ is to draw samples from $\pi_{\text{post}}$ in the vicinity of a trial $\vect{\theta}_{\text{dyn}}$. We tested several MCMC sampling strategies (see \cite{villalobos2023:scientific} for details). In the present case, sampling from $\pi_{\text{prior}}$ is \emph{not} intractable by virtue of the simple forward model~\cref{e:fhnfwd}. However, to obtain an accurate estimate for $\mat{\Gamma}_{\text{post}}$ we still require a large number of samples, and consequently a large number of evaluations of $\vect{f}$. To avoid this, we adopt ideas from~\cite{Stuart:2010a,BuiThanh:2013a,Martin:2012a} and approximate $\mat{\Gamma}_{\text{post}}$ based on the inverse of the Hessian $\mat{H}$ associated with~\cref{e:neglogpost}. We compared this strategy to the results obtained via MCMC sampling. The estimates were in good agreement. Thus, we only report results for the data generated via $\mat{H}^{-1} \approx \mat{\Gamma}_{\text{post}}$.

Under the assumption that the parameter-to-observation map $\vect{f}$ is differentiable, we can linearize the right hand side of~\cref{e:yobs} around a given sample $\vect{\theta}^i_{\text{dyn}}$ to obtain
\[
\vect{y}_{\text{obs}}
\approx \vect{f}(\vect{\theta}^i_{\text{dyn}})
+ \mat{F}^i(\vect{\theta}_{\text{dyn}} - \vect{\theta}_{\text{dyn}}^i) + \vect{\eta},
\]

\noindent where $\mat{F}^i$ is the Jacobian of $\vect{f}$ evaluated at $\vect{\theta}_{\text{dyn}}^i \in \left\{\Theta_{\text{train}}, \Theta_{\text{test}}\right\}$. It follows that we can approximate the posterior covariance matrix for the $i$-th sample of our training dataset as
\[
\mat{\Gamma}_{\text{post}}^i \approx \mat{\Gamma}_{\text{train}}^i = ((\mat{F}^i)^\mathsf{H} \mat{\Gamma}_{\text{noise}}^{-1} \mat{F}^i + \mat{\Gamma}_{\text{prior}}^{-1})^{-1}.
\]

Comparing this linearization to a quadratic approximation of the objective function in~\cref{e:neglogpost} reveals that we can approximate $\mat{\Gamma}_{\text{post}}^i$ associated with a sample $\vect{\theta}_{\text{dyn}}^i$ based on the inverse of the Hessian matrix $\mat{H}^i \in \ns{R}^{3,3}$,
\[
\mat{\Gamma}_{\text{post}}^i \approx (\mat{H}^i)^{-1} = (\mat{H}^i_{\text{dat}} + \mat{H}^i_{\text{reg}})^{-1}.
\]

\noindent If the parameter-to-observation map $\vect{f}$ were linear this approximation would be exact. We present the derivation as well as steps necessary to compute the Hessian matrix $\mat{H}$ for~\cref{e:varopt} in the supplementary materials in \Cref{s:curvature_supp}. The matrices $\mat{H}_{\text{dat}}$ and $\mat{H}_{\text{reg}}$ are discretized versions of the continuous operators $\mathcal{H}_{\text{reg}}$ and $\mathcal{H}_{\text{dat}}$ of the reduced space Hessian operator $\mathcal{H}$ associated with~\cref{e:varopt}. Since $\mat{H}$ is a $3\times 3$ matrix, we can simply invert it using direct methods.

\subsubsection{Mapping \texorpdfstring{$\spd{3}$}{Sym+(3)} into \texorpdfstring{$\mathbb{R}^6$}{R6}}\label{s:logeuclidean}

The proposed approach for learning the posterior covariance matrices builds upon ideas described in~\cite{Arsigny:2005a,Arsigny:2006a,Arsigny:2007a,Pennec:2006a}. Since we use MSE as a loss function, we assume that the parameters $\vect{\theta}$ live in an Euclidean space. This does not hold for the posterior covariance $\mat{\Gamma}_{\text{post}} \in \spd{3}$. The space $ \spd{3}$ when endowed with a Riemannian metric forms a Riemannian manifold that is locally similar to an Euclidean space. Computing distances between two elements of $\spd{3}$ does, in general, require the computation of a geodesic. While the geodesics can, in general, be evaluated in closed form for $\spd{3}$, this strategy would introduce additional complications to our NN framework. To simplify the associated computations, we consider a log-Euclidean metric on $\spd{3}$~\cite{Arsigny:2005a,Arsigny:2006a,Arsigny:2007a}. Using this approach allows us to use classical Euclidean computations in the domain of matrix logarithms~\cite{Arsigny:2007a}. In particular, we are mapping elements of $\spd{3}$ to the space of symmetric matrices $\sym{3}$. These matrices can uniquely be represented as a vector in $\mathbb{R}^6$. This is the data we use for training our NN. There are several alternative strategies; this approach is particularly simple and does not require us to change the loss function for the NN or the NN architecture.

In particular, we use this log-Euclidean framework to map elements in $\spd{3}$ to their tangent space using the Riemannian logarithmic map $\log_{\mat{\Sigma}} : \spd{3} \to T_{\mat{\Sigma}} \spd{3}$,
\[
\log_{\mat{\Sigma}}(\mat{W}) = \mat{\Sigma}^{1/2} \log(\mat{\Sigma}^{-1/2}\mat{W}\mat{\Sigma}^{-1/2})\mat{\Sigma}^{1/2}
\]

\noindent for $\mat{\Sigma} \in \spd{3}$~\cite{Pennec:2006a}; here, $\log$ denotes the usual matrix logarithm. This allows us to perform classical Euclidean computations in the domain of matrix logarithms. The inverse of the logarithmic map $\log_{\mat{\Sigma}}$is given by the Riemannian exponential map $\exp_{\mat{\Sigma}} : T_{\mat{\Sigma}} \spd{3} \to \spd{3}$,
\[
\exp_{\mat{\Sigma}}(\mat{V}) = \mat{\Sigma}^{1/2} \exp(\mat{\Sigma}^{-1/2}\mat{V}\mat{\Sigma}^{-1/2})\mat{\Sigma}^{1/2}
\]

\noindent for $\mat{\Sigma} \in \spd{3}$ that associates to each $\mat{V} \in T_\Sigma \spd{3}$ a point of the manifold $\spd{3}$~\cite{Pennec:2006a}. Here, $\exp(\cdot)$ denotes the usual matrix exponential. Notice that the exponential map and its inverse, the logarithm, are both smooth and invertible maps, i.e., diffeomorphisms~\cite{Arsigny:2007a}.

We evaluate $\log_{\mat{\Sigma}}$ and $\exp_{\mat{\Sigma}}$ at $\mat{\Sigma} = \mat{\operatorname{Id}} = \operatorname{diag}(1,1,1) \in \mathbb{R}^{3,3}$. After applying $\log_{\mat{\operatorname{Id}}}$, we are working with symmetric matrices with $(3(3+1)/2) = 6$ independent coefficients (e.g., the upper triangular part). As such, we can map these matrices from $\sym{3}$ to $\mathbb{R}^6$. Let $\mat{W} = \log_{\mat{\operatorname{Id}}}(\mat{\Gamma}_{\text{post}})$. The associated assignment $\operatorname{vec}_{\mat{\operatorname{Id}}}: \sym{3} \to \mathbb{R}^6$ is given by
\[
\operatorname{vec}_{\mat{\operatorname{Id}}}(\mat{W}) = (w_{00}, \sqrt{2}w_{01}, \sqrt{2}w_{02}, w_{11}, \sqrt{2}w_{12}, w_{22}) \in \mathbb{R}^6.
\]

We apply this mapping to the entire training and testing dataset. Consequently, the entries predicted by the NN architecture are all vectors in $\mathbb{R}^6$, and by that elements of $\sym{3} = T_{\mat{\operatorname{Id}}} \spd{3}$ after appropriate rescaling and re-arrangement. To obtain the prediction of $\mat{\Gamma}_{\text{post}}$ we apply the inverse of the projection $\operatorname{vec}_{\mat{\operatorname{Id}}}$ and subsequently the Riemannian exponential $\exp_{\mat{\operatorname{Id}}}$ to map from $T_{\mat{\operatorname{Id}}} \spd{3}$ to $\spd{3}$.

\subsubsection{Noise Models}\label{s:noisemodel}

Our first noise model is first-order autoregressive in time, following~\cite{rudi2022:parameter}. These types of models are widely used for representing time series. The autoregressive process is parametrized by a correlation parameter $\rho$, which determines the dependence of the process on its previous value~\cite{Mills:1990a}. The $j$th entry of $\vect{y}_{\text{obs}}\in\ns{R}^{n_t}$ for the additive noise model is given by
\begin{equation}\label{e:addnoise}
y_{\text{obs},j} = u(t^j) + \eta(t^j),
\quad
\eta(t^j) = \rho \eta(t^{j-1}) + \varepsilon(t^j),
\end{equation}

\noindent at time point $t^j = jh_t \in [0,\tau]$, $j = 1,\ldots, n_t$, with $h_t = \tau/n_t$, $\eta(t^1) \sim \mathcal{N}(0,\sigma^2/h_t^2)$ and $|\rho| < 1$. Notice that the variance of the noise model depends on the time step size $h_t = \tau / n_t$. Moreover, the variance of $\eta$ is constant across time since the process is stationary due to $|\rho| < 1$ (i.e., $\eta(t^j) \sim \mathcal{N}(0,\sigma^2/h_t^2)$ for all $j = 1,\ldots, n_t$). Moreover, we select $\varepsilon(t^j) \sim \mathcal{N}(0,(1-\rho^2) \sigma^2/h_t^2)$. The noise parameters are $\vect{\theta}_{\text{noise}} = (\rho,\sigma) \in \ns{R}^2$.

In our second noise model, which we call intrinsic noise, we replace the equation for $u$ in~\cref{e:fhnfwd} by the SDE
\begin{equation}\label{e:sdenoise}
\d Y_t - \mu(Y_t, t)\d t - \beta\d W_t = 0,
\end{equation}

\noindent with initial condition $Y_0 = u_0$. The term $\mu(Y_t, t)$ can be interpreted as the expectation and is given by $\mu(Y_t, t) = \theta_2 (Y_t - (Y_t^3/3) + v + z)$. The parameter $\beta > 0$ can be interpreted as the variance of the model, where $\{W_t\}_{t \geq 0}$ is a Wiener process (standard Brownian motion). The $j$th entry of $\vect{y}_{\text{obs}}\in\ns{R}^{n_t}$ for the intrinsic noise model is given by $y_{\text{obs},j} = Y_{t^j}$ for all $t^j = jh_t$, $j = 1,\ldots,n_t$. We use a Euler--Maruyama scheme (e.g., see \cite{Higham:2000}) for the integration of the SDE in~\cref{e:sdenoise}. We have $\vect{\theta}_{\text{noise}} = \beta \in \ns{R}$.

\subsubsection{Parameter Selection}\label{s:parameters}

We discuss the selection of the ODE model parameters $\vect{\theta}_{\text{dyn}}$ first. We use
\begin{equation}\label{e:priorpara}
\vect{\sigma}_{\text{prior}} = (0.3, 0.4, 0.4)
\quad\text{and}\quad
\bar{\vect{\theta}}_{\text{prior}} = (0.4, 0.4, 3.4),
\end{equation}

\noindent respectively, such that $[\bar{\theta}_{\text{prior},i} - 2\sigma_{\text{prior},i}, \bar{\theta}_{\text{prior},i} + 2\sigma_{\text{prior},i}]$ for $i = 0,1,2$. We additionally limit the samples from the prior distribution for each component of $\vect{\theta}_{\text{dyn}}$ to construct $\{\Theta_{\text{train}},\Theta_{\text{test}}\}$ to the ranges
\begin{equation}\label{e:parabounds}
\theta_0 \in [-0.2,1.0], \quad \theta_1 \in [-0.4,1.2], \quad \theta_2 \in [2.0, 5.0].
\end{equation}

\noindent We discard samples for $\vect{\theta}_{\text{dyn}}$ that do not fall into these ranges.

The values of the parameter $\rho$ and $\sigma$ controlling the additive noise perturbation in~\cref{e:addnoise} are varied randomly across different samples $\vect{y}_{\text{obs}}^i$. We generate 100 independent samples, with $\rho \sim\mathcal{N}(0.8,0.05^2)$ and $\sigma \sim\mathcal{N}(0.07,0.01^2)$. Selecting $\rho$ this way allows us to achieve a correlation of 0.65 to 0.95. Likewise, selecting $\sigma$ this way yields a noise level of 4\% to 10\%. For each pair $(\rho^j, \sigma^j)$, independent replicates of the process $\eta$ are generated and added to training and/or testing data. The parameter $\beta$ that controls the intrinsic noise in~\cref{e:sdenoise} is chosen according to $\beta \sim \mathcal{N}(0.15, 0.05^2)$. In addition, we restrict the samples for $\beta$ to the range $[0.01, 0.27]$. This range has been determined empirically. We summarize these (and other) parameters in \Cref{t:parameter-values}.

\begin{table}
\caption{Parameters considered for generating the training and testing data.\label{t:parameter-values}}
\adjtab\centering
\begin{tabular}{llll}\toprule
\bf variable & \bf  meaning & \bf value/range & \bf equation \\
\midrule
$\tau$                                & final time                   & 200\,ms                         & see \cref{e:fhnfwd}   \\
$n_t$                                 & number of time steps         & \inum{2000}                     & ---                   \\
$z$                                   & membrane potential           & $-0.4$                          & see \cref{e:fhnfwd}   \\
$\theta_0$                            & model parameter (excitability threshold)             & $[-0.2, 1.0]$                   & see \cref{e:fhnfwd}   \\
$\theta_1$                            & model parameter (recovery scaling)            & $[-0.4, 1.2]$                   & see \cref{e:fhnfwd}   \\
$\theta_2$                            & model parameter (time scale separation) & $[2.0, 5.0]$                    &see \cref{e:fhnfwd}   \\
$\vect{\sigma}_{\text{prior}}$        & prior standard deviation     & $(0.3,0.4, 0.4)$                & see \cref{e:prior}    \\
$\vect{\bar{\theta}}_{\text{prior}}$  & prior mean                   & $(0.3,0.4, 3.4)$                & see \cref{e:prior}    \\
$\rho$                                & noise correlation            & $\sim\mathcal{N}(0.8,0.05^2)$   & see \cref{e:addnoise} \\
$\sigma^2$                            & noise variance               & $\sim\mathcal{N}(0.07,0.01^2)$  & see \cref{e:addnoise} \\
$\beta^2$                             & noise variance               & $\sim\mathcal{N}(0.15, 0.05^2)$ & see \cref{e:sdenoise} \\
\midrule
$n_{\text{train}}$ & number of training samples   & \inum{1000}, \inum{2000}, \inum{4000}, \inum{8000} & --- \\
$n_{\text{test}}$  & number of testing samples    & \inum{4000}                                        & --- \\
$n_{\text{val}}$   & number of validation samples & \inum{2000}                                        & --- \\
\bottomrule
\end{tabular}
\end{table}

\subsection{Performance Measures}\label{s:performance}

We assess the performance of the proposed approach qualitatively and quantitatively. We note that, since we use simulations to generate our data we have access to the true parameters and hence can use this data to evaluate the performance of the proposed methodology. We split the generated data into training and testing data to provide an unbiased estimate of the model's predictive capabilities. We report results in \Cref{s:results}. Here, we summarize quantitative performance measures. We consider various measures of model performance to get a better understanding of the sources of error in our predictive model. We explain the benefits and differences for each measure below. We decompose the MSE into its squared bias (\iacr{SQB}) and centered MSE (\iacr{CMSE}) components to assess the contributions of the mean and of the variability of the prediction mismatch~\cite{Taylor:2001a}. In particular, a large SQB (high bias) indicates that the model is too simple (i.e., is indicative of underfitting) whereas a large CMSE (high variance) suggests that the model is too sensitive to the training data (i.e., is indicative of overfitting). These measures are given by
\begin{equation}\label{e:sqb}
\textstyle\varepsilon_{\text{SQB}}(\{\vect{\theta}_{\text{pred}}^i\}_{i=1}^{n_{\text{test}}}) =
\|\vect{\bar{\theta}}_{\text{test}} - \vect{\bar{\theta}}_{\text{pred}}\|^2_2
\end{equation}

\noindent with means $\vect{\bar{\theta}}_{\text{test}} = \sum_{i=1}^{n_{\text{test}}} \vect{\theta}_{\text{test}}^i/n_{\text{test}}$ and $\vect{\bar{\theta}}_{\text{pred}} = \sum_{i=1}^{n_{\text{test}}} \vect{\theta}_{\text{pred}}^i/n_{\text{test}}$, and
\begin{equation}\label{e:cmse}
\varepsilon_{\text{CMSE}}(\{\vect{\theta}_{\text{pred}}^i\}_{i=1}^{n_{\text{test}}}) =
\frac{1}{n_s} \sum_{i=1}^{n_{\text{test}}}
\left\|
(\vect{\theta}^i_{\text{test}} - \vect{\bar{\theta}}_{\text{test}}) - (\vect{\theta}^i_{\text{pred}} - \vect{\bar{\theta}}_{\text{pred}})
\right\|_2^2.
\end{equation}

\noindent We also report values for the median of absolute percentage error (\iacr{MdAPE}) given by
\begin{equation}\label{e:mdape}
\varepsilon_{\text{MdAPE}}(\{\vect{\theta}_{\text{pred}}^i\}_{i=1}^{n_{\text{test}}}) =
\operatorname{median}_{i=1,\ldots,n_{\text{test}}}
\left\{
\|\vect{\theta}^i_{\text{test}} - \vect{\theta}^i_{\text{pred}}\|_2
/\|\vect{\theta}^i_{\text{test}}\|_2
\right\}
\end{equation}

This yields a scale invariant performance measure that is less sensitive to outliers (we use the median) compared to the mean absolute percentage error (\iacr{MAPE}) and MSE, i.e., it is more robust. Since MdAPE is a relative error it is interpretable. We also report values for the coefficient of determination (\iacr{CDET}), also called $R^2$, given by
\begin{equation}\label{e:cdet}
\varepsilon_{\text{CDET}}(\{\vect{\theta}_{\text{pred}}^i\}_{i=1}^{n_{\text{test}}}) =
1 -
\sum_{i=1}^{n_{\text{test}}}
\|\vect{\theta}^i_{\text{test}} - \vect{\theta}^i_{\text{pred}}\|_2^2/
\sum_{i=1}^{n_{\text{test}}}
\|\vect{\theta}^i_{\text{test}} - \vect{\bar{\theta}}_{\text{test}}\|_2^2.
\end{equation}

\noindent This measure is used to assess the percentage of variability explained by the prediction. It takes values in $(-\infty,1]$ with 1 being optimal. The closer $\varepsilon_{\text{CDET}}$ is to a value of 1, the more variability is explained by the prediction. This measure is more sensitive to outliers than MdAPE. However, compared to the  other measures it accounts for class imbalance. Throughout our experiments, the number of training samples $n_{\text{train}}$ is varied ranging from $\inum{1000}$ to $\inum{8000}$. We summarize these choices in \Cref{t:parameter-values}.

\section{Numerical Experiments}\label{s:results}

We report results under various conditions. If not otherwise stated, the parameters for the NN architectures are as in \Cref{t:nnconfig}; the parameters for data generation are as in \Cref{t:parameter-values}. We consider a \nn\ noise configuration for all experiments (training and testing samples contain noise). Due to page limitations we have moved additional experiments with varying configurations to the supplementary materials (see \Cref{s:supplementary_results}). The software environment is described in \Cref{s:swhw}.

\subsection{ODE Model and Noise Parameter Estimation}\label{s:model_noise_results}

We report results for point estimates for the ODE parameters $\vect{\theta}_{\text{dyn}} = (\theta_0,\theta_1,\theta_2) \in \ns{R}^3$ along with the parameters $\vect{\theta}_{\text{noise}}$ controlling the noise. We consider additive and/or intrinsic noise. The experiment in \Cref{s:nn_exploration} explores how different NN architectures affect the estimation of $\vect{\theta}_{\text{dyn}}$. Subsequently, we study the performance when additionally estimating additive noise parameters (see \Cref{s:perf_additive_noise}) and a combination of additive and intrinsic noise parameters (see \Cref{s:perf_comb_noise}). The results for solely estimating intrinsic noise parameters (i.e., without additive noise) can be found in the supplementary materials (see \Cref{s:noise_model_para_intrinsic_supp}).

\subsubsection{Exploration of NN Architectures}\label{s:nn_exploration}

This experiment assesses the predictive capabilities of a variety of NN architectures when we restrict ourselves to the estimation of model parameters $\vect{\theta}_{\text{dyn}}$; hence, each NN's output dimension is $p=3$. We intend to empirically identify the architecture with optimal performance to be used for the remaining numerical study. We refer to \cite{villalobos2023:scientific} for a more detailed analysis of the NN architecture and hyperparameter turning.

\ipoint{Setup:} As reported in~\cite{rudi2022:parameter}, a learning rate of 0.002, a batch size of 32, and 64 epochs (i.e., \inum{8000} optimization steps) are adequate for training. We train with \inum{4000} samples, which gives representative results for selecting NN architectures, and validate with \inum{2000} samples. The input data (features) are of type \ts. For DNN architectures we vary
\begin{inparaenum}[\it (i)]
\item the number of layers from 4 to 20 and
\item the number of units per layer $n_u$ from 4 to 256\end{inparaenum}.
For CNN architectures we vary
\begin{inparaenum}[\it (i)]
\item the number of convolutional layers from 2 to 6 and
\item the filter setting $n_f$ from 2 to 128\end{inparaenum}.
Recall that our CNN architecture has two additional dense layers after the convolutions (see details in \Cref{s:cnn}).

\ipoint{Results:} We report the MdAPE and CDET evaluation measures for the DNN architectures in \Cref{t:dnn_exploration_mdape_cdet} and for the CNN in \Cref{t:cnn_exploration_mdape_cdet}.

\begin{table}
\caption{Exploration of DNN architectures. Each cell shows the prediction accuracy on validation data for estimating $\vect{\theta}_{\text{dyn}}$ using the MdAPE and CDET (in brackets) measures. We vary the number of nodes/units $n_u$ per layer (rows) and the number of layers (columns).}
\label{t:dnn_exploration_mdape_cdet}
\centering\adjtab
\begin{tabular}{rccccc}
\toprule
$n_u$ & \textbf{4 layers} & \textbf{8 layers}  & \textbf{12 layers} & \textbf{16 layers} & \textbf{20 layers} \\
\midrule
  4 & \decc{0.130296} (\decc{0.651815}) & \decc{0.127854} (\decc{0.635456}) & \decc{0.165107} (\decc{0.525304}) & \decc{0.179563} (\decc{0.449644}) & \decc{0.189365} (\decc{0.412267})\\
  8 & \decc{0.115974} (\decc{0.727150}) & \decc{0.103461} (\decc{0.767597}) & \decc{0.101814} (\decc{0.773822}) & \decc{0.134969} (\decc{0.660060}) & \decc{0.183756} (\decc{0.454346})\\
 16 & \decc{0.081502} (\decc{0.851448}) & \decc{0.079739} (\decc{0.858554}) & \decc{0.075696} (\decc{0.845268}) & \decc{0.077593} (\decc{0.858490}) & \decc{0.110943} (\decc{0.711280})\\
 32 & \decc{0.070495} (\decc{0.869497}) & \decc{0.066619} (\decc{0.878074}) & \decc{0.074638} (\decc{0.867082}) & \decc{0.084130} (\decc{0.849438}) & \decc{0.109177} (\decc{0.818431})\\
 64 & \decc{0.063871} (\decc{0.893590}) & \decc{0.072298} (\decc{0.872880}) & \decc{0.066560} (\decc{0.886807}) & \decc{0.066673} (\decc{0.892906}) & \decc{0.080249} (\decc{0.848662})\\
128 & \decc{0.072605} (\decc{0.878481}) & \decc{0.068796} (\decc{0.888196}) & \decc{0.064691} (\decc{0.892048}) & \decc{0.068380} (\decc{0.884710}) & \decc{0.208190} (\decc{0.432143})\\
256 & \decc{0.065729} (\decc{0.892368}) & \decc{0.061797} (\decc{0.892323}) & \decc{0.108884} (\decc{0.713904}) & \decc{0.106209} (\decc{0.716812}) & \decc{0.173088} (\decc{0.511089})\\
\bottomrule
\end{tabular}
\end{table}

\begin{table}
\caption{Exploration of CNN architectures. Each cell shows the prediction accuracy on validation data for estimating $\vect{\theta}_{\text{dyn}}$ using the MdAPE and CDET (in brackets) measures. We vary the filter setting $n_f$ per layer (rows) and the number of convolutional layers (columns). Recall that a full CNN has an additional two dense layers after the convolutions.}
\label{t:cnn_exploration_mdape_cdet}
\centering\adjtab
\begin{tabular}{rccccc}
\toprule
$n_f$ & \textbf{2 layers} & \textbf{3 layers}  & \textbf{4 layers} & \textbf{5 layers} & \textbf{6 layers}  \\
\midrule
  2 & \decc{0.059008} (\decc{0.900448}) & \decc{0.053433} (\decc{0.923805}) & \decc{0.053167} (\decc{0.930770}) & \decc{0.053218} (\decc{0.925716}) & \decc{0.052064} (\ms\decc{0.927544})  \\
  4 & \decc{0.061381} (\decc{0.894812}) & \decc{0.053538} (\decc{0.922970}) & \decc{0.050357} (\decc{0.934373}) & \decc{0.051069} (\decc{0.935552}) & \decc{0.051517} (\ms\decc{0.926882}) \\
  8 & \decc{0.063428} (\decc{0.893834}) & \decc{0.055269} (\decc{0.912723}) & \decc{0.050875} (\decc{0.930282}) & \decc{0.051623} (\decc{0.935068}) & \decc{0.053464} (\ms\decc{0.931126}) \\
 16 & \decc{0.069896} (\decc{0.885347}) & \decc{0.060578} (\decc{0.903873}) & \decc{0.058659} (\decc{0.914716}) & \decc{0.055171} (\decc{0.922843}) & \decc{0.056675} (\ms\decc{0.919264}) \\
 32 & \decc{0.068277} (\decc{0.884943}) & \decc{0.063872} (\decc{0.898209}) & \decc{0.058041} (\decc{0.917503}) & \decc{0.054536} (\decc{0.915978}) & \decc{0.057621} (\ms\decc{0.917683}) \\
 64 & \decc{0.073525} (\decc{0.881121}) & \decc{0.063968} (\decc{0.894534}) & \decc{0.067809} (\decc{0.902709}) & \decc{0.061082} (\decc{0.918326}) & \decc{0.061196} (\ms\decc{0.906260}) \\
128 & \decc{0.071374} (\decc{0.872677}) & \decc{0.067779} (\decc{0.889505}) & \decc{0.060340} (\decc{0.907794}) & \decc{0.063307} (\decc{0.911815}) & \decc{0.953771} (\decc{-7.268274}) \\
\bottomrule
\end{tabular}
\end{table}

\ipoint{Observations:} DNN architectures yield reasonable prediction accuracy with CDET $>0.80$ for a wide range of settings, namely 4 to 16 layers and for $n_u$ between 16 and 128. The peak accuracy over all DNNs is $\sim\!\!0.89$. However, across CNN architectures prediction accuracy is significantly better with CDET $>0.90$ for a wide range of networks: 3 to 6 convolutional layers and for $n_f$ between 2 and 128. The peak accuracy over all CNNs reaches $\sim\!\!0.93$, which is $0.04$ points better than for the peak for DNNs; since CDET is close to its maximum value of one, an improvement of $0.04$ is significant. Complementary results in the supplementary materials show that the CMSE metric behaves analogously to CDET and reveal the same trends for DNN and CNN architectures as we observe here; hence, both measures are adequate to assess prediction accuracy. The comparison between DNN and CNN architectures suggests that the use of convolution operations is advantageous for the considered time series data. We attribute this to the fact that convolutions take into account the locality of neighboring time steps in the data, which is not true for dense layers. The CNNs have a lower number of trainable weights, because convolution kernels are $3\times3$ matrices. For instance, the CNN with 5 layers and $n_f=8$ has \inum{38195} trainable parameters, compared with the (slightly worse performing) DNN with 12 layers and $n_u=128$, which amounts to \inum{695179} trainable parameters. On the other hand, DNNs have an advantage in terms of computational efficiency since their arithmetic operations can utilize hardware accelerated kernels for matrix-matrix multiplications. Based on these observations we restrict the remaining experiments to a CNN architecture with 5 layers and $n_f = 8$.

\subsubsection{Estimating ODE and Noise Parameters -- Additive Noise}\label{s:perf_additive_noise}

We assess the performance of NN predictions for the simultaneous estimation of ODE model parameters $\vect{\theta}_{\text{dyn}} = (\theta_0,\theta_1,\theta_2)\in\mathbb{R}^3$ and noise parameters $\vect{\theta}_{\text{noise}} = (\rho,\sigma) \in \mathbb{R}^2$ of the additive noise model (see \Cref{s:noisemodel}); hence, each NN will have output dimension $p=5$. This setup allows to quantify properties of the noise of a data example at nearly the same computational cost as estimating only $\vect{\theta}_{\text{dyn}} \in \mathbb{R}^3$. As will be demonstrated, the NNs are able to estimate parameters of noise---at effectively no additional cost---along with parameters of the ODE, which is a unique capability compared with other methods for parameter estimation or inference.

\ipoint{Setup:} We consider the three different types of NN inputs: \ts, \fc, and \tsfc, where the latter is obtained through stacking of \ts\ and \fc. We will see that choice of features has an impact on the accuracy of estimating $\vect{\theta}_{\text{noise}}$. In addition, we vary the number of training samples between \inum{500} and \inum{8000} to study the effects on estimating $\vect{\theta}_{\text{dyn}}$ and $\vect{\theta}_{\text{noise}}$ as the amount of training data increases. The trained CNNs are assessed with \inum{4000} testing samples that are distinct from the validation data set used in \Cref{s:nn_exploration}. The selection of the NN architecture is based on the experiments reported in \Cref{s:nn_exploration}; we choose a CNN with 5 layers and $n_f=8$.

\begin{figure}
\centering
\includegraphics[width=0.8\textwidth]{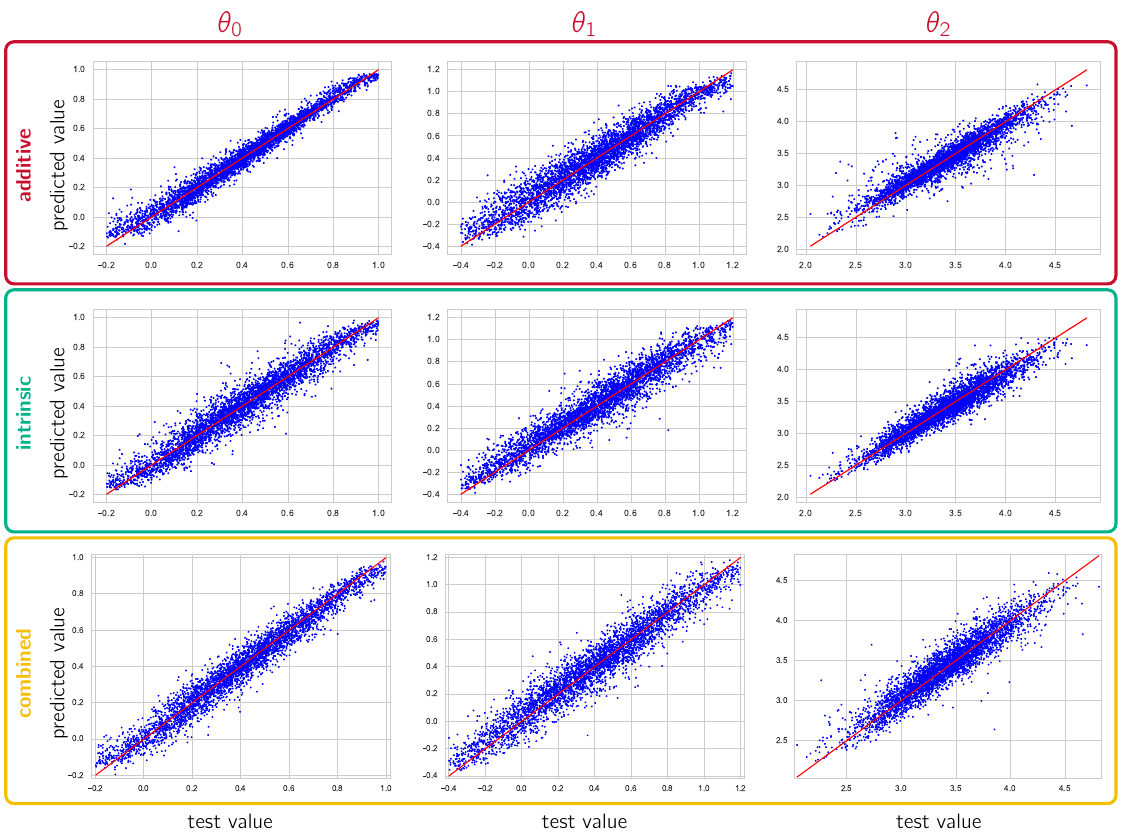}
\caption{Scatter plots for CNN predictions using testing data when estimating $\vect{\theta}_{\text{dyn}} = (\theta_0,\theta_1,\theta_2) \in \mathbb{R}^3$ for the data type \tsfc\ and a training size of $n_{\text{train}} = \inum{8000}$. The CNN also infers $\vect{\theta}_{\text{noise}}$, but this is not shown. The three vertical blocks correspond to additive noise (top block), intrinsic noise (middle block), and combined noise (bottom block). We plot the value of the test samples ($x$-axis) versus the predicted values ($y$-axis) for $\theta_0$ (left), $\theta_1$ (center), and $\theta_2$ (right). The diagonal red line indicates optimal predictions. Overall, the CNNs deliver similar predictions across different noise configurations.}
\label{f:para_recovery_various_noises_scatter}
\end{figure}

\ipoint{Results:} We report the MdAPE and CDET evaluation measures in \Cref{t:cnn_nzadd_noise_recovery}. We show scatter plots for the prediction of the model parameters $\vect{\theta}_{\text{dyn}}$ in \Cref{f:para_recovery_various_noises_scatter} for different noise models, data type \tsfc, and a training size of $n_{\text{train}} = \inum{8000}$; top block ``additive.'' The horizontal and vertical axes correspond to the true and predicted values, respectively. Optimal predictions are indicated by a diagonal red line.

\begin{table}
\centering
\caption{Estimation accuracy of CNN on testing data when inferring $\vect{\theta}_{\text{dyn}} = (\theta_0,\theta_1,\theta_2) \in \mathbb{R}^3$ and $\vect{\theta}_{\text{noise}} = (\rho,\sigma) \in \mathbb{R}^2$. Each cell shows MdAPE and CDET (in parenthesis). Blocks of three rows have a varying number of training samples $n_{\text{train}}$, while within each block the input data to the NN is different (\ts: time series data; \fc: Fourier coefficients; \tsfc: combination of \ts\ and \fc).}
\label{t:cnn_nzadd_noise_recovery}
\centering\adjtab
\begin{tabular}{rlccccc}
\toprule
$n_{\text{train}}$ & \bf input & $\theta_0$ & $\theta_1$ & $\theta_2$  & $\rho$ & $\sigma$ \\
\midrule
            & \ts   & \decc{0.101439} (\decc{0.933649}) & \decc{0.266450} (\decc{0.806941}) & \decc{0.029839} (\decc{0.730307}) & \decc{0.057937} (\decc{-0.287224})   & \decc{0.090218} (\decc{-0.193864})   \\
 \inum{500} & \fc   & \decc{0.262849} (\decc{0.416492}) & \decc{0.429213} (\decc{0.360127}) & \decc{0.047466} (\decc{0.548418}) & \decc{0.027010} (\ms\decc{0.669485}) & \decc{0.060297} (\ms\decc{0.554192}) \\
            & \tsfc & \decc{0.113781} (\decc{0.927969}) & \decc{0.350972} (\decc{0.697450}) & \decc{0.034559} (\decc{0.647204}) & \decc{0.028099} (\ms\decc{0.622094}) & \decc{0.060980} (\ms\decc{0.532075}) \\
\hline
            & \ts   & \decc{0.101208} (\decc{0.938824}) & \decc{0.256867} (\decc{0.815738}) & \decc{0.029027} (\decc{0.765650}) & \decc{0.054221} (\decc{-0.147207})   & \decc{0.085869} (\decc{-0.125742})   \\
\inum{1000} & \fc   & \decc{0.262265} (\decc{0.505308}) & \decc{0.396063} (\decc{0.553181}) & \decc{0.038625} (\decc{0.663139}) & \decc{0.024207} (\ms\decc{0.703232}) & \decc{0.053644} (\ms\decc{0.659086}) \\
            & \tsfc & \decc{0.092851} (\decc{0.949637}) & \decc{0.241944} (\decc{0.825879}) & \decc{0.026018} (\decc{0.777421}) & \decc{0.022113} (\ms\decc{0.758500}) & \decc{0.047686} (\ms\decc{0.717913}) \\
\hline
            & \ts   & \decc{0.098635} (\decc{0.937374}) & \decc{0.206227} (\decc{0.872403}) & \decc{0.028488} (\decc{0.793807}) & \decc{0.049890} (\decc{-0.146354})   & \decc{0.088641} (\ms\decc{0.007138}) \\
\inum{2000} & \fc   & \decc{0.217402} (\decc{0.581193}) & \decc{0.357994} (\decc{0.627635}) & \decc{0.034308} (\decc{0.730496}) & \decc{0.021610} (\ms\decc{0.759084}) & \decc{0.048124} (\ms\decc{0.702317}) \\
            & \tsfc & \decc{0.093610} (\decc{0.954207}) & \decc{0.192823} (\decc{0.886672}) & \decc{0.023397} (\decc{0.823516}) & \decc{0.020576} (\ms\decc{0.788201}) & \decc{0.043652} (\ms\decc{0.754122}) \\
\hline
            & \ts   & \decc{0.093524} (\decc{0.945307}) & \decc{0.186488} (\decc{0.889913}) & \decc{0.025013} (\decc{0.825750}) & \decc{0.046883} (\decc{-0.005214})   & \decc{0.064941} (\ms\decc{0.480397}) \\
\inum{4000} & \fc   & \decc{0.199514} (\decc{0.532929}) & \decc{0.320374} (\decc{0.664997}) & \decc{0.032567} (\decc{0.721381}) & \decc{0.019821} (\ms\decc{0.789644}) & \decc{0.045283} (\ms\decc{0.745332}) \\
            & \tsfc & \decc{0.070591} (\decc{0.965736}) & \decc{0.166104} (\decc{0.900889}) & \decc{0.021403} (\decc{0.858754}) & \decc{0.018455} (\ms\decc{0.833946}) & \decc{0.039849} (\ms\decc{0.804174}) \\
\hline
            & \ts   & \decc{0.070273} (\decc{0.965882}) & \decc{0.159959} (\decc{0.916073}) & \decc{0.021728} (\decc{0.870171}) & \decc{0.026125} (\ms\decc{0.684876}) & \decc{0.044998} (\ms\decc{0.761624}) \\
\inum{8000} & \fc   & \decc{0.186824} (\decc{0.571823}) & \decc{0.290906} (\decc{0.700348}) & \decc{0.029644} (\decc{0.769661}) & \decc{0.019116} (\ms\decc{0.815396}) & \decc{0.040955} (\ms\decc{0.787872}) \\
            & \tsfc & \decc{0.070485} (\decc{0.969701}) & \decc{0.157937} (\decc{0.919323}) & \decc{0.020204} (\decc{0.876198}) & \decc{0.016319} (\ms\decc{0.861656}) & \decc{0.034751} (\ms\decc{0.834789}) \\
\bottomrule
\end{tabular}
\end{table}

\ipoint{Observations:} The Fourier coefficients play a critical role in the recovery of the noise parameters $\vect{\theta}_{\text{noise}}$. For instance, in \Cref{t:cnn_nzadd_noise_recovery} with $n_{\text{train}}=\inum{2000}$, we obtain a CDET for $\rho$ and $\sigma$ with \ts\ inputs that is at or below zero---a very poor accuracy---but \fc\ inputs yield CDET values $>0.70$. However, the Fourier coefficients by themselves (\fc\ rows in the table) do not yield as good of an accuracy for $\vect{\theta}_{\text{dyn}}$ as with the \ts\ data. The best results for both $\vect{\theta}_{\text{dyn}}$ and $\vect{\theta}_{\text{noise}}$ are consequently obtained with \tsfc, a concatenation of the \ts\ and \fc\ features.

\subsubsection{Estimating ODE and Noise Parameters -- Combined Noise}\label{s:perf_comb_noise}

We extend the previous experiments in \Cref{s:perf_additive_noise} to include intrinsic noise in addition to additive noise. Therefore, we jointly estimate ODE model parameters $\vect{\theta}_{\text{dyn}} = (\theta_0,\theta_1,\theta_2) \in \mathbb{R}^3$ and noise parameters $\vect{\theta}_{\text{noise}} = (\rho,\sigma,\beta) \in \mathbb{R}^3$. Consequently, each NN will have output dimension $p=6$.

\ipoint{Setup:} This experiment uses the same configuration for data and CNN as the experiment in \Cref{s:perf_additive_noise}. Since we combine two noise models, the intensity in each noise model is reduced by 50\% to obtain a similar signal-to-noise ratio as in the former experiment.

\ipoint{Results:} We report the MdAPE and CDET evaluation measures in \Cref{t:comb_half_noise_cnn2}. We complement the table with scatter plots for the prediction of the model parameters $\vect{\theta}_{\text{dyn}}$ in \Cref{f:para_recovery_various_noises_scatter}, bottom block ``combined.''

\ipoint{Observations:} The overall behavior is in line with previous experiments:
\begin{inparaenum}[\it (i)]
\item relatively low CDET for the recovery of noise parameters with \ts\ data;
\item \fc\ data significantly improves the accuracy for the noise parameters;
\item \tsfc\ data maintains the same accuracy for ODE model parameters as for pure \ts\ data, and achieves similar accuracy for noise parameters as \fc\ data\end{inparaenum}.
Comparing the accuracy obtained here with the experiment in \Cref{s:perf_additive_noise}, we observe slightly worse values for MdAPE (and CDET) for the noise parameters.
For example, the best accuracy in \Cref{s:perf_additive_noise} reached CDET values $(0.862, 0.835)$ for $(\rho, \sigma)$ versus $(0.795, 0.729)$ now. This indicates that the presence of both autocorrelated and intrinsic noise makes the parameter estimation more challenging. However, the accuracy for the ODE model parameters $\vect{\theta}_{\text{dyn}}$ remains mostly similar to the experiment in \Cref{s:perf_additive_noise}, which is also confirmed by the scatter plots in the top and bottom rows of \Cref{f:para_recovery_various_noises_scatter}.

\begin{table}
\caption{Estimation accuracy of CNN on testing data when inferring $\vect{\theta}_{\text{dyn}} = (\theta_0,\theta_1,\theta_2)$ and $\vect{\theta}_{\text{noise}} = (\rho,\sigma,\beta)$. Each cell shows MdAPE (and CDET in parenthesis). Blocks of three rows have varying number of training samples $n_{\text{train}}$, while within each block the input data to NN is different (\ts: time series data; \fc: Fourier coefficients; \tsfc: combination of \ts\ and \fc).}
\label{t:comb_half_noise_cnn2}
\centering\adjtab
%\resizebox{\textwidth}{!}{
\begin{tabular}{rlcccccc}
\toprule
$n_{\text{train}}$ & \bf input  & $\theta_0$ & $\theta_1$  & $\theta_2$  &  $\beta$ & $\rho$  &  $\sigma$ \\
\midrule
& \ts   & \decc{0.087856} (\decc{0.931032}) & \decc{0.151545} (\decc{0.842867}) & \decc{0.036869} (\decc{0.689685}) & \decc{0.263299} (\decc{-0.102905})   & \decc{0.048271} (\decc{-0.112096}) & \decc{0.110772} (\decc{-0.696250}) \\
\inum{500}
& \fc   & \decc{0.226797} (\decc{0.392128}) & \decc{0.290173} (\decc{0.439774}) & \decc{0.053434} (\decc{0.469541}) & \decc{0.150072} (\ms\decc{0.646217}) & \decc{0.027256} (\ms\decc{0.653993}) & \decc{0.057015} (\ms\decc{0.616751}) \\
& \tsfc & \decc{0.101416} (\decc{0.943279}) & \decc{0.196804} (\decc{0.882914}) & \decc{0.033404} (\decc{0.773820}) & \decc{0.132117} (\ms\decc{0.696717}) & \decc{0.023524} (\ms\decc{0.705094}) & \decc{0.050396} (\ms\decc{0.664925}) \\
\hline
& \ts   & \decc{0.093305} (\decc{0.928214}) & \decc{0.140426} (\decc{0.880715}) & \decc{0.033675} (\decc{0.753156}) & \decc{0.213159} (\ms\decc{0.235619}) & \decc{0.038944} (\ms\decc{0.235849}) & \decc{0.114066} (\decc{-0.655869}) \\
\inum{1000}
& \fc   & \decc{0.213786} (\decc{0.348985}) & \decc{0.201330} (\decc{0.605284}) & \decc{0.041078} (\decc{0.668570}) & \decc{0.132193} (\ms\decc{0.697395}) & \decc{0.024129} (\ms\decc{0.702071}) & \decc{0.052254} (\ms\decc{0.638311}) \\
& \tsfc & \decc{0.098424} (\decc{0.943663}) & \decc{0.195479} (\decc{0.868887}) & \decc{0.029092} (\decc{0.833750}) & \decc{0.122914} (\ms\decc{0.752381}) & \decc{0.021716} (\ms\decc{0.754401}) & \decc{0.045948} (\ms\decc{0.706680}) \\
\hline
&  \ts  & \decc{0.075069} (\decc{0.945826}) & \decc{0.107314} (\decc{0.913046}) & \decc{0.028432} (\decc{0.812520}) & \decc{0.182375} (\ms\decc{0.433945}) & \decc{0.034151} (\ms\decc{0.433614}) & \decc{0.113941} (\decc{-0.670531}) \\
\inum{2000}
& \fc   & \decc{0.203285} (\decc{0.362672}) & \decc{0.200241} (\decc{0.634591}) & \decc{0.036058} (\decc{0.740631}) & \decc{0.127002} (\ms\decc{0.723347}) & \decc{0.022695} (\ms\decc{0.727000}) & \decc{0.050126} (\ms\decc{0.667707}) \\
& \tsfc & \decc{0.097986} (\decc{0.941880}) & \decc{0.176958} (\decc{0.907081}) & \decc{0.026568} (\decc{0.855300}) & \decc{0.117531} (\ms\decc{0.767103}) & \decc{0.020655} (\ms\decc{0.768247}) & \decc{0.045025} (\ms\decc{0.711821}) \\
\hline
&  \ts  & \decc{0.075004} (\decc{0.951502}) & \decc{0.098813} (\decc{0.928155}) & \decc{0.025743} (\decc{0.859021}) & \decc{0.163770} (\ms\decc{0.527264}) & \decc{0.029284} (\ms\decc{0.528056}) & \decc{0.111371} (\decc{-0.623537}) \\
\inum{4000}
& \fc   & \decc{0.179181} (\decc{0.412661}) & \decc{0.183971} (\decc{0.689532}) & \decc{0.033842} (\decc{0.777731}) & \decc{0.123315} (\ms\decc{0.752676}) & \decc{0.021676} (\ms\decc{0.753574}) & \decc{0.049573} (\ms\decc{0.679104}) \\
& \tsfc & \decc{0.088314} (\decc{0.959217}) & \decc{0.163748} (\decc{0.920807}) & \decc{0.025307} (\decc{0.868371}) & \decc{0.115815} (\ms\decc{0.777542}) & \decc{0.020368} (\ms\decc{0.776767}) & \decc{0.045301} (\ms\decc{0.717113}) \\
\hline
& \ts   & \decc{0.072967} (\decc{0.958661}) & \decc{0.088472} (\decc{0.943355}) & \decc{0.023922} (\decc{0.871125}) & \decc{0.148733} (\ms\decc{0.651308}) & \decc{0.026946} (\ms\decc{0.651350}) & \decc{0.104678} (\decc{-0.378932}) \\
\inum{8000}
& \fc   & \decc{0.170018} (\decc{0.421513}) & \decc{0.164342} (\decc{0.723233}) & \decc{0.031986} (\decc{0.807645}) & \decc{0.123314} (\ms\decc{0.753598}) & \decc{0.021318} (\ms\decc{0.753598}) & \decc{0.047201} (\ms\decc{0.695808}) \\
& \tsfc & \decc{0.088424} (\decc{0.957964}) & \decc{0.148462} (\decc{0.928605}) & \decc{0.024706} (\decc{0.877221}) & \decc{0.108834} (\ms\decc{0.794461}) & \decc{0.018912} (\ms\decc{0.794536}) & \decc{0.043019} (\ms\decc{0.728900}) \\
\bottomrule
\end{tabular}%}
\end{table}

\subsection{Uncertainty Quantification -- Posterior Covariance Estimation}\label{s:covariance_results}

We jointly estimate ODE parameters $\vect{\theta}_{\text{dyn}} = (\theta_0, \theta_1, \theta_2)\in\mathbb{R}^3$, noise parameters $\vect{\theta}_{\text{noise}} = (\beta, \rho, \sigma) \in \mathbb{R}^3$ of different noise conditions (additive noise, intrinsic noise, and a combination thereof), and entries of the covariance matrix $\mat{\Gamma}_{\text{post}} \in\spd{3}$. For the latter we use a Laplace approximation of the posterior distribution of the model parameters $\vect{\theta}_{\text{dyn}}$ conditioned on $\vect{y}_{\text{obs}}$. Estimating $\mat{\Gamma}_{\text{post}}$ allows us to assess uncertainties locally around $\vect{\theta}_{\text{dyn}}$. By combining the three output vectors/matrices, the NN will have a total output dimension $p=12$, where the covariance matrix $\mat{\Gamma}_{\text{post}}$ contributes six values $(\Gamma_{00}, \Gamma_{01}, \Gamma_{02}, \Gamma_{11}, \Gamma_{21}, \Gamma_{22})$; we note that we use a log-Euclidean framework to project the data from the Riemannian manifold $\spd{3}$ to $\mathbb{R}^6$. This allows us to perform classical Euclidean computations in the domain of matrix logarithms (see \Cref{s:logeuclidean}). With the task to predict $\mat{\Gamma}_{\text{post}}$, we are asking the question if the CNN is capable to determine a property that is not directly determined from its input $\vect{y}_{\text{obs}}$ but predict additional information that we impose, namely $\mat{\Gamma}_{\text{post}}$.

\ipoint{Setup:}
We utilize \tsfc\ as NN inputs. For each sample of $\vect{\theta}_{\text{dyn}}$ we compute $\mat{\Gamma}_{\text{post}}$ by inverting the Hessian matrix $\mat{H}$ of the optimization problem \Cref{e:varoptred} at the (known) solution $\vect{\theta}_{\text{dyn}}$ (see \Cref{s:data} for more details). While the matrix $\mat{H}$ associated with each sample is expected to be positive definite, $\mat{H} \in \spd{3}$, we observe that singular matrices can occur in practice. This is likely caused by the weak prior term (see \Cref{s:parameters}) combined with a highly varying data misfit term (see \Cref{f:lscp}) and numerical inaccuracies. We train only on the subset of samples satisfying $\mat{H}  \in \spd{3}$; we discard samples that violate that condition, because $\mat{H}$ does not satisfy conditions to represent an inverse covariance matrix. We decided to not regenerate discarded samples, because our goal was to remain consistent with previous experiments. As a consequence, the number of the training and testing samples varies slightly compared with the previous experiments: Overall, our dataset consists of \inum{15000} samples. \inum{110} samples yielded negative definite Hessians. An additional \inum{1893} Hessians were ill-conditioned. This leaves a total of \inum{12997} samples for which we construct approximations for the associated covariances $\mat{\Gamma}_{\text{post}}$. We provide additional details in the supplementary material accompanying this manuscript.

\ipoint{Results:} We report MdAPE and CDET evaluation measures in \Cref{t:cnn_param_and_cov_recovery_logeuclidean}. Each block of three rows corresponds to different training sizes $n_{\text{train}}$. Each row in each block corresponds to a different noise setting (additive, intrinsic, and combined).

\begin{table}
\caption{Estimation accuracy of CNNs on testing data when inferring $\vect{\theta}_{\text{dyn}}$, $\vect{\theta}_{\text{noise}}$, and $\mat{\Gamma}_{\text{post}}$. Each cell shows MdAPE (and CDET in parenthesis).
We report values for the model and noise parameters (A: additive noise; I: intrinsic noise; C: combined noise) in the top table and for the covariance entries in the bottom table. The training sizes differ from previous experiments; we refer to the text for more details.}
\label{t:cnn_param_and_cov_recovery_logeuclidean}
\centering\adjtab
%\resizebox{\textwidth}{!}{
\begin{tabular}{rccccccc}
\toprule
$n_{\text{train}}$ & \bf noise  & $\theta_0$ & $\theta_1$  & $\theta_2$  & $\beta$  & $\rho$  &  $\sigma$ \\
\midrule
            & A & \decc{0.120473} (\decc{0.919174}) & \decc{0.280166} (\decc{0.777022}) & \decc{0.029328} (\decc{0.689130}) & --- & \decc{0.031634} (\decc{0.503844}) & \decc{0.071393} (\decc{0.335957}) \\
\inum{432}  & I & \decc{0.151745} (\decc{0.876634}) & \decc{0.231978} (\decc{0.835559}) & \decc{0.042778} (\decc{0.657371}) & \decc{0.089256} (\decc{0.870905}) & --- & --- \\
            & C & \decc{0.107784} (\decc{0.931468}) & \decc{0.210020} (\decc{0.873567}) & \decc{0.036586} (\decc{0.679760}) & \decc{0.144466} (\decc{0.654369}) & \decc{0.026074} (\decc{0.672242}) & \decc{0.055604} (\decc{0.614005})\\
\hline
            & A & \decc{0.094256} (\decc{0.947663}) & \decc{0.228580} (\decc{0.846818}) & \decc{0.025747} (\decc{0.794706}) & --- & \decc{0.025903} (\decc{0.669372}) & \decc{0.055090} (\decc{0.628457}) \\
\inum{872}  & I & \decc{0.137619} (\decc{0.885941}) & \decc{0.201075} (\decc{0.861447}) & \decc{0.036429} (\decc{0.747383}) & \decc{0.082278} (\decc{0.888932}) & --- & --- \\
            & C & \decc{0.097260} (\decc{0.944305}) & \decc{0.185138} (\decc{0.893183}) & \decc{0.030531} (\decc{0.794677}) & \decc{0.127377} (\decc{0.740333}) & \decc{0.023155} (\decc{0.738124}) & \decc{0.049238} (\decc{0.659749})\\
\hline
            & A & \decc{0.092025} (\decc{0.953657}) & \decc{0.194301} (\decc{0.879172}) & \decc{0.023776} (\decc{0.829927}) & --- & \decc{0.025600} (\decc{0.714490}) & \decc{0.055182} (\decc{0.625053})\\
\inum{1742} & I & \decc{0.126289} (\decc{0.905758}) & \decc{0.191210} (\decc{0.874523}) & \decc{0.031732} (\decc{0.810370}) & \decc{0.075339} (\decc{0.904059}) & --- & --- \\
            & C & \decc{0.090455} (\decc{0.954155}) & \decc{0.162899} (\decc{0.915183}) & \decc{0.027618} (\decc{0.828420}) & \decc{0.116912} (\decc{0.786520}) & \decc{0.020540} (\decc{0.784984}) & \decc{0.049570} (\decc{0.672528})\\
\hline
            & A & \decc{0.092672} (\decc{0.956495}) & \decc{0.165577} (\decc{0.913093}) & \decc{0.021614} (\decc{0.863625}) & --- & \decc{0.019918} (\decc{0.812990}) & \decc{0.042264} (\decc{0.783791}) \\
\inum{3480} & I & \decc{0.149439} (\decc{0.886746}) & \decc{0.214246} (\decc{0.874506}) & \decc{0.029285} (\decc{0.835745}) & \decc{0.073105} (\decc{0.919152}) & --- & --- \\
            & C & \decc{0.090406} (\decc{0.955899}) & \decc{0.158824} (\decc{0.923942}) & \decc{0.023564} (\decc{0.878426}) & \decc{0.107307} (\decc{0.806381}) & \decc{0.019002} (\decc{0.805294}) & \decc{0.043148} (\decc{0.748188}) \\
\hline
            & A & \decc{0.062077} (\decc{0.973380}) & \decc{0.135659} (\decc{0.939087}) & \decc{0.016559} (\decc{0.889840}) & --- & \decc{0.016326} (\decc{0.862704}) & \decc{0.039805} (\decc{0.812831}) \\
\inum{6928} & I & \decc{0.125508} (\decc{0.906996}) & \decc{0.165011} (\decc{0.905912}) & \decc{0.026624} (\decc{0.862155}) & \decc{0.065145} (\decc{0.936454}) & --- & ---  \\
            & C & \decc{0.088671} (\decc{0.951811}) & \decc{0.163238} (\decc{0.902820}) & \decc{0.024979} (\decc{0.879952}) & \decc{0.102294} (\decc{0.824335}) & \decc{0.018318} (\decc{0.824203}) & \decc{0.043822} (\decc{0.745046})\\
\bottomrule
&&&&&&& \\
\toprule
$n_{\text{train}}$ & \bf noise  & $\Gamma_{00}$ & $\Gamma_{01}$  & $\Gamma_{02}$  &  $\Gamma_{11}$ & $\Gamma_{12}$  &  $\Gamma_{22}$ \\
\midrule
            & A & \decc{0.056769} (\decc{0.923542}) & \decc{0.450306} (\decc{0.783533}) & \decc{0.126280} (\decc{0.929961}) & \decc{0.132439} (\decc{0.780382}) & \decc{0.348502} (\decc{0.793183}) & \decc{0.057696} (\decc{0.874591}) \\
 \inum{432} & I & \decc{0.086740} (\decc{0.737715}) & \decc{0.415801} (\decc{0.740252}) & \decc{0.148644} (\decc{0.837406}) & \decc{0.131336} (\decc{0.749414}) & \decc{0.346102} (\decc{0.795970}) & \decc{0.082577} (\decc{0.799428}) \\
            & C & \decc{0.065122} (\decc{0.864175}) & \decc{0.341092} (\decc{0.832132}) & \decc{0.128307} (\decc{0.917654}) & \decc{0.131678} (\decc{0.798481}) & \decc{0.294300} (\decc{0.846668}) & \decc{0.066782} (\decc{0.871909}) \\
\hline
            & A & \decc{0.048756} (\decc{0.945958}) & \decc{0.338227} (\decc{0.863958}) & \decc{0.105936} (\decc{0.956767}) & \decc{0.111201} (\decc{0.845151}) & \decc{0.275609} (\decc{0.861147}) & \decc{0.051923} (\decc{0.909932}) \\
 \inum{872} & I & \decc{0.087960} (\decc{0.749160}) & \decc{0.382357} (\decc{0.766137}) & \decc{0.148374} (\decc{0.850931}) & \decc{0.113455} (\decc{0.787180}) & \decc{0.304206} (\decc{0.846524}) & \decc{0.073766} (\decc{0.816626}) \\
            & C & \decc{0.061229} (\decc{0.895155}) & \decc{0.338276} (\decc{0.855845}) & \decc{0.115543} (\decc{0.934016}) & \decc{0.097304} (\decc{0.851844}) & \decc{0.265633} (\decc{0.880667}) & \decc{0.059348} (\decc{0.899473}) \\
\hline
            & A & \decc{0.040406} (\decc{0.955123}) & \decc{0.302682} (\decc{0.891784}) & \decc{0.094179} (\decc{0.963127}) & \decc{0.097138} (\decc{0.885048}) & \decc{0.230511} (\decc{0.889465}) & \decc{0.048509} (\decc{0.925574}) \\
\inum{1742} & I & \decc{0.079066} (\decc{0.782115}) & \decc{0.386955} (\decc{0.792731}) & \decc{0.130067} (\decc{0.868195}) & \decc{0.107702} (\decc{0.815438}) & \decc{0.269986} (\decc{0.865249}) & \decc{0.072174} (\decc{0.830937}) \\
            & C & \decc{0.053256} (\decc{0.911435}) & \decc{0.300216} (\decc{0.880277}) & \decc{0.102554} (\decc{0.942916}) & \decc{0.089048} (\decc{0.880161}) & \decc{0.235407} (\decc{0.896117}) & \decc{0.055133} (\decc{0.912030}) \\
\hline
            & A & \decc{0.045133} (\decc{0.952062}) & \decc{0.285468} (\decc{0.909411}) & \decc{0.089084} (\decc{0.964451}) & \decc{0.085733} (\decc{0.907260}) & \decc{0.195608} (\decc{0.917350}) & \decc{0.043429} (\decc{0.940409}) \\
\inum{3480} & I & \decc{0.089134} (\decc{0.752620}) & \decc{0.442385} (\decc{0.760935}) & \decc{0.137006} (\decc{0.855165}) & \decc{0.111555} (\decc{0.819692}) & \decc{0.304154} (\decc{0.853566}) & \decc{0.075353} (\decc{0.831304}) \\
            & C & \decc{0.052160} (\decc{0.902665}) & \decc{0.285895} (\decc{0.897218}) & \decc{0.096858} (\decc{0.942729}) & \decc{0.086445} (\decc{0.896912}) & \decc{0.226925} (\decc{0.913924}) & \decc{0.053945} (\decc{0.916414}) \\
\hline
            & A & \decc{0.033070} (\decc{0.969617}) & \decc{0.223061} (\decc{0.934648}) & \decc{0.077670} (\decc{0.976941}) & \decc{0.077082} (\decc{0.926089}) & \decc{0.160791} (\decc{0.936642}) & \decc{0.039036} (\decc{0.959481}) \\
\inum{6928} & I & \decc{0.078800} (\decc{0.778226}) & \decc{0.351692} (\decc{0.819258}) & \decc{0.123398} (\decc{0.865146}) & \decc{0.092191} (\decc{0.842558}) & \decc{0.245837} (\decc{0.884544}) & \decc{0.076118} (\decc{0.827159}) \\
            & C & \decc{0.050679} (\decc{0.900514}) & \decc{0.316940} (\decc{0.875516}) & \decc{0.102854} (\decc{0.944785}) & \decc{0.094560} (\decc{0.878225}) & \decc{0.251510} (\decc{0.888930}) & \decc{0.057664} (\decc{0.912421}) \\
\bottomrule
\end{tabular}%}
\end{table}

\ipoint{Observations:} The most important observation is that the NN-based approach for inference is able to scale to larger output dimensions without increasing training size or NN architecture.

 The accuracy of predicting model parameters (see columns $(\theta_0,\theta_1,\theta_2)$) remains nearly identical to the experiment in \Cref{s:perf_comb_noise}, where only $\vect{\theta}_{\text{dyn}}$ and $\vect{\theta}_{\text{noise}}$ were NN outputs (i.e., without $\mat{\Gamma}_{\text{post}}$; see \Cref{t:comb_half_noise_cnn2}). The accuracy of predicting noise parameters is generally consistent with the previous experiment. We observe in many cases a slightly lower accuracy (see columns $(\beta,\rho,\sigma)$ and rows of combined noise type) compared with \Cref{t:comb_half_noise_cnn2}. This reduction in accuracy only shows up in CDET values at lower training sizes $n_{\text{train}}<\inum{2000}$. For larger training sizes ($n_{\text{train}} \in \{\inum{3480}, \inum{6928}\}$) recovery results across all parameters ($\vect{\theta}_{\text{dyn}}$ and $\vect{\theta}_{\text{noise}}$) are comparable to those without trying to predict covariances. This observation holds for all noise conditions. Overall, nearly the same performance for different NN outputs---when compared with the experiment in \Cref{s:perf_comb_noise}---shows a remarkable property of NN to scale to larger data dimensions without sacrificing accuracy.

When estimating the covariance entries, we observe that NN predictions achieve an accuracy that is comparable to the one for the ODE and/or noise parameters. An increase in training size also increases accuracy, which is a trend seen in previous experiments. At the largest training size, $n_{\text{train}} = \inum{6928}$, all but one of the covariance entry reach CDET values above $0.8$.

\subsection{Sensitivity to Initialization}\label{s:nn_initialization}

In this experiment, we assess the sensitivity of the NN architectures to the initialization of the network.

\ipoint{Setup:} We report results for different random seeds (i.e., weight initializations), focusing on recovery of model parameters, noise parameters, and covariance entries under our combined noise model. We consider 10 different random initialization. We use a DNN with $n_u = 128$ nodes per layer and 12 hidden layers, and a CNN with 5 hidden layers and $n_f = 8$. For these experiments, training was conducted for 100 epochs using our largest training set size ($n_{\text{train}} = 6928$ samples) and tested on $n_{\text{test}} = 3456$ samples.

\ipoint{Results:} \Cref{t:rand_seed_inits} summarize the recovery performance of the CNN and DNN architectures across seeds to evaluate the effects of random weight initialization. We report detailed results in the supplementary material in \Cref{s:nn_initialization_supp} and \Cref{s:nn_initialization_supp_cc}, respectively. We also include results for a 6-fold cross validation in the supplementary material.

\ipoint{Observations:}  Across all seeds tested, model performance was consistent, with no significant impact from random weight initialization. Both CNN and DNN architectures exhibited stable metrics (e.g., CDET and MdAPE) with low variability. No anomalies or outliers were observed. The CNN generally achieved slightly higher CDET values and lower errors. Additionally, training times were also consistent across seeds, with the CNN averaging 73.43 seconds total, and the DNN averaging 48.44 seconds total.

\begin{table}
\caption{Sensitivity to initialization. We report the mean and standard deviation for the NN-predictions of the model parameters $\vect{\theta}_{\text{dyn}}$ and covariance matrix entries across different seed initializations. The top block corresponds to the DNN architecture. The bottom block corresponds to the CNN architecture.\label{t:rand_seed_inits}}
\label{t:cnn_summary_stats}
\centering\adjtab
\begin{tabular}{llrrrr}\toprule
\textbf{architecture} & \textbf{metric} & \multicolumn{2}{c}{\textbf{parameter recovery}} & \multicolumn{2}{c}{\textbf{covariance recovery}} \\
\midrule
 & & \textbf{mean} & \textbf{std dev} & \textbf{mean} & \textbf{std dev} \\
\midrule
CNN
& CDET  & \decc{0.85564800} & \sci{0.00631200}  & \decc{0.79114600} & \sci{0.00967900} \\
& MdAPE & \decc{0.07175100} & \sci{0.00198300}  & \decc{0.17410400} & \sci{0.00313700} \\
& CMSE  & \decc{0.03049000} & \sci{0.00179000}  & \sci{0.00127500}  & \sci{0.00006500} \\
& SQB   & \sci{0.00041800}  & \sci{0.00034100}  & \sci{0.00000700}  & \sci{0.00000300} \\
\midrule
DNN
& CDET  & \decc{0.78248000} & \sci{0.00620700} & \decc{0.64649800}  & \sci{0.00421400} \\
& MdAPE & \decc{0.10373000} & \sci{0.00162400} & \decc{0.22701900}  & \sci{0.00255700} \\
& CMSE  & \decc{0.06868400} & \sci{0.00110900}  & \sci{0.00209700}  & \sci{0.00003600} \\
& SQB   & \sci{0.00092300}  & \sci{0.00136200}  & \sci{0.00002300}  & \sci{0.00000900} \\
\bottomrule
\end{tabular}
\end{table}

\subsection{Sensitivity of NN Predictions to Noise}\label{s:noise_sensitivity}

We explore the performance of the proposed framework as a function of increasing noise perturbations. We expect the prediction accuracy to deteriorate as the noise level increases. We note that these results are baked into those reported in prior experiments. Here, we attempt to disentangle the errors and report prediction performance as a function of increasing noise perturbations. Note that we draw 100 noise parameter samples; and each noise parameter is used multiple times to generate multiple different realizations of noise. Consequently, we report statistics averaged across different predictions for each sample of noise parameter.

\ipoint{Setup:} These results correspond to experiments in \Cref{s:covariance_results}, in particular, those reported in \Cref{t:cnn_param_and_cov_recovery_logeuclidean}, where inference is performed for model parameter, noise parameters, and covariance entries.

\begin{figure}
\centering
\includegraphics[width=0.7\textwidth]{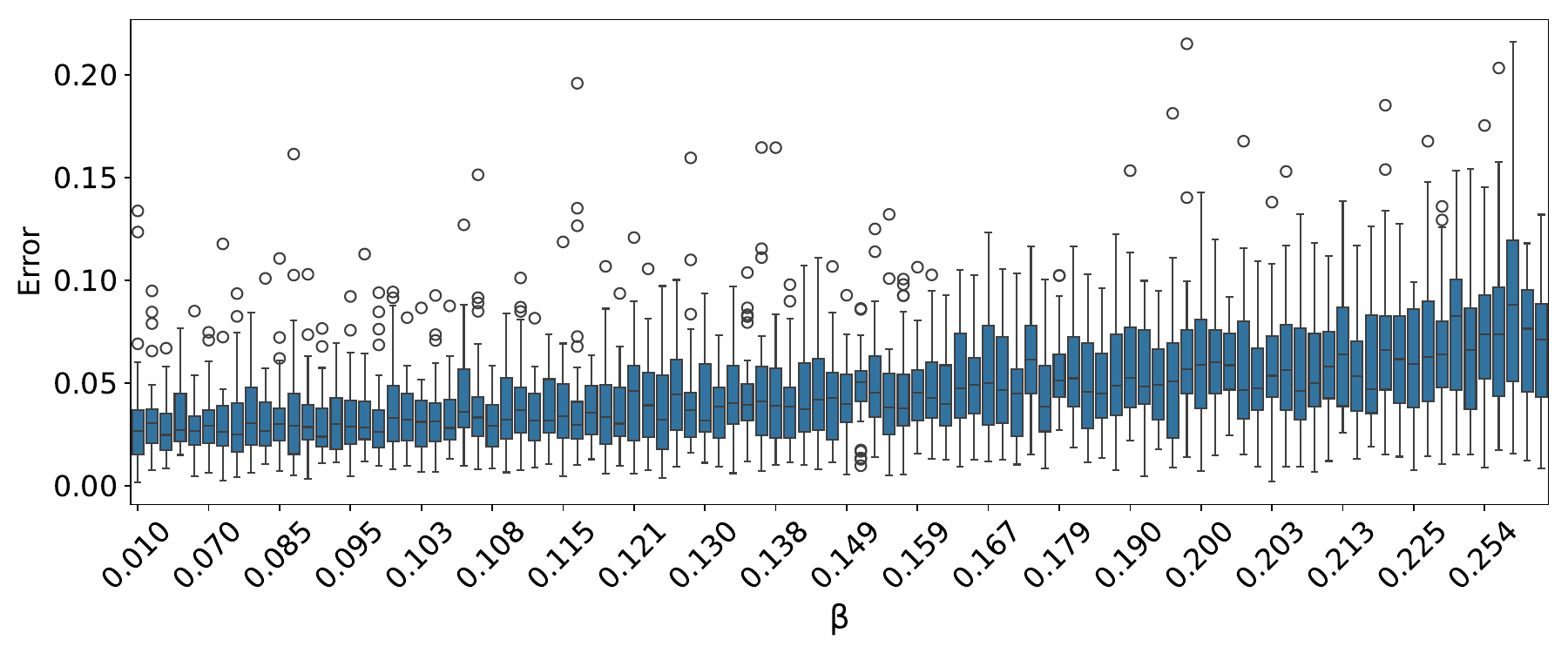}
\includegraphics[width=0.7\textwidth]{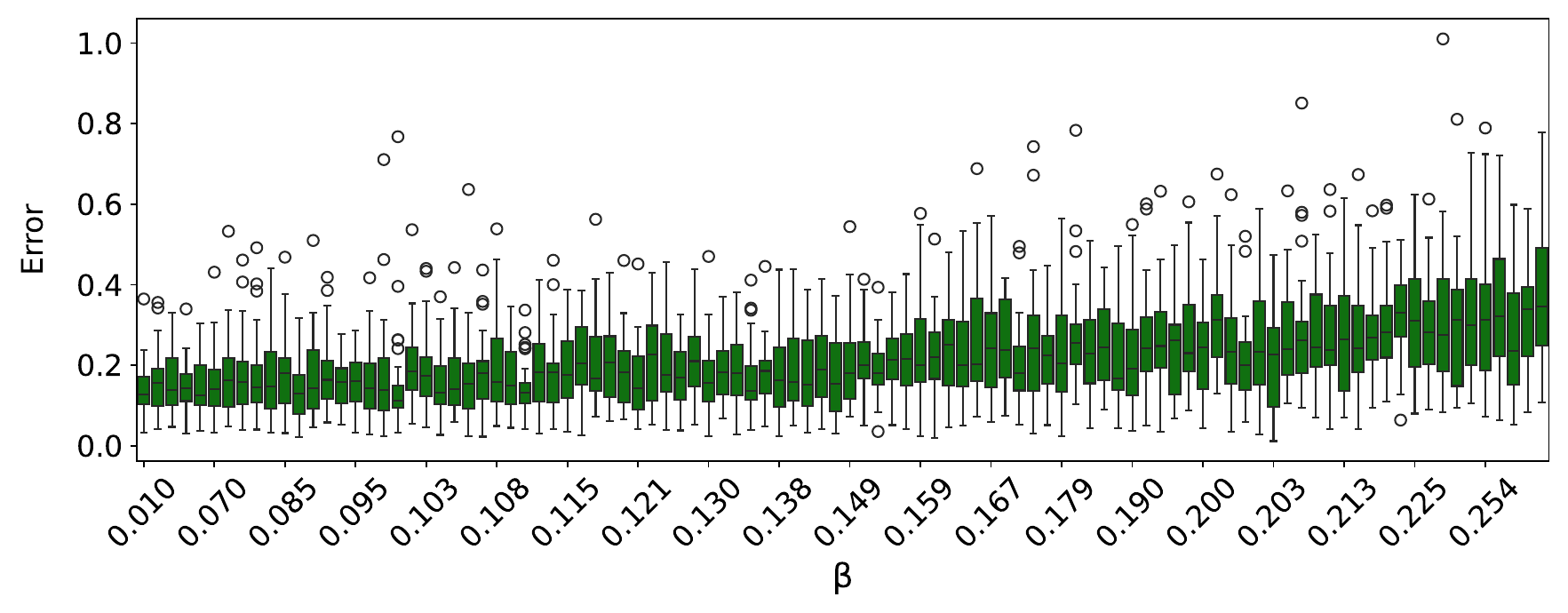}
\caption{Statistics of the relative $\ell^2$-error for the intrinsic noise model as a function of the noise parameter $\beta$. We report the $\ell^2$-norm for the predicted entries of the model parameter vector $\vect{\theta}_{\text{dyn}}$ (top; blue) and the entries of the covariance matrix (bottom; green). The boundaries of each of these boxes correspond to the 25th and 75th percentile (lower and upper quartile), respectively. The line inside the box represents the 50th percentile. The whiskers correspond to the data points that are within 1.5 times the interquartile range from the lower and upper quartiles. The circles above or below the whiskers denote points outside this range, i.e. outliers.}\label{fig:box_whisker_intrinsic_noise}
\end{figure}

\ipoint{Results:} Since the intrinsic noise is controlled by a single parameter $\beta$, reporting the error as a function of increasing $\beta$ is straightforward. We provide box-whisker plots for the error in the ODE parameters and the covariance matrix entries in \Cref{fig:box_whisker_intrinsic_noise}. We report the relative $\ell^2$-error, which corresponds to MdAPE in \cref{e:mdape}. For the additive noise model the noise is controlled by two parameters $(\rho,\sigma) \in \mathbb{R}^2$. We provide box-whisker plots for the estimation errors as each of these parameters increases in \Cref{fig:box_whisker_additive_noise}. We provide 2D plots in the supplementary material (see \Cref{s:noise_sensitivity_supp}).

\begin{figure}
\centering
\includegraphics[width=0.7\textwidth]{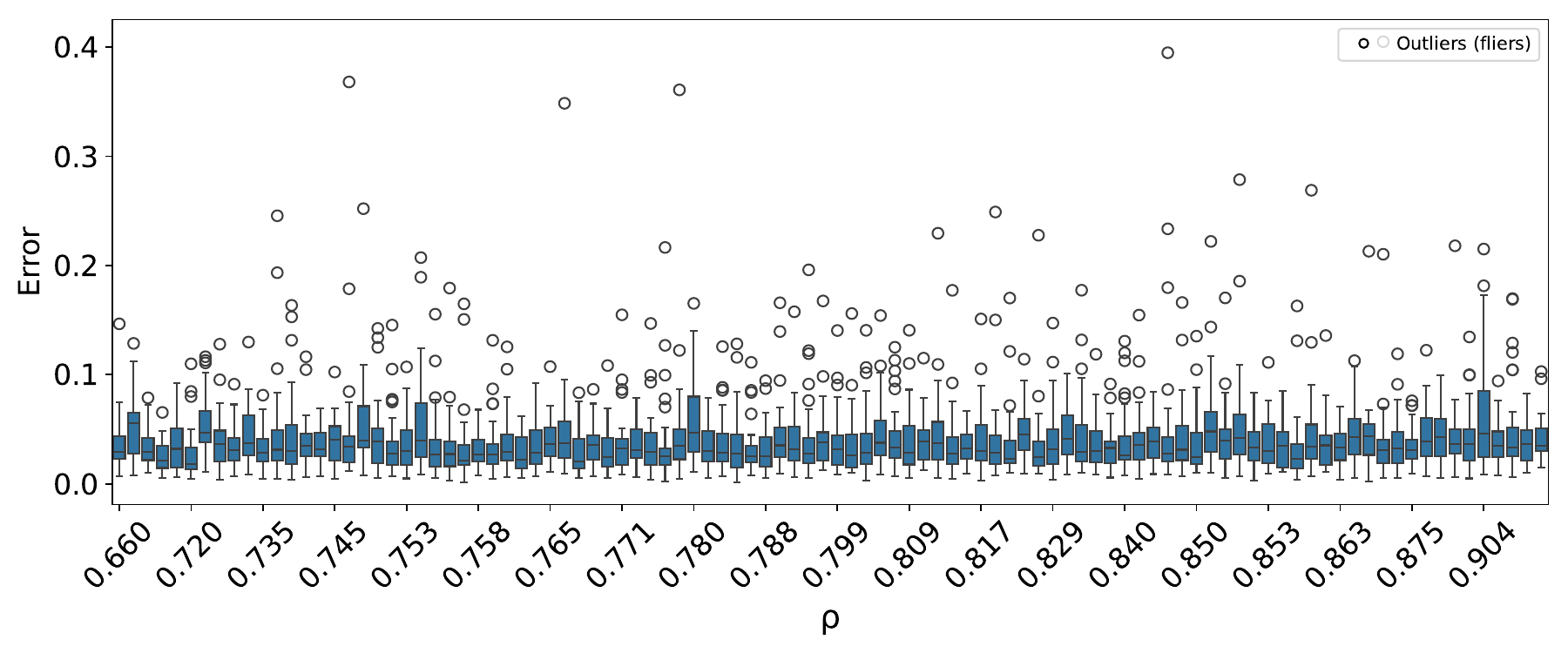}
\includegraphics[width=0.7\textwidth]{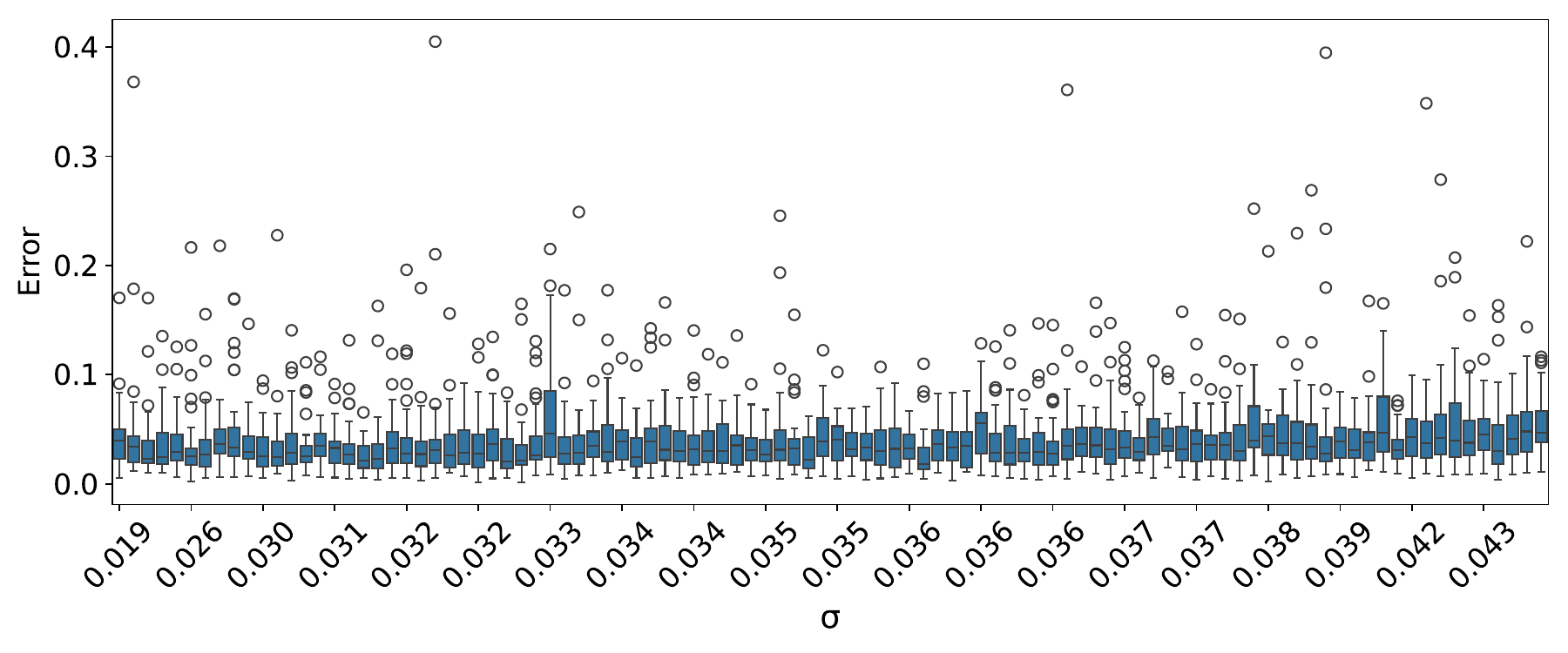}
\includegraphics[width=0.7\textwidth]{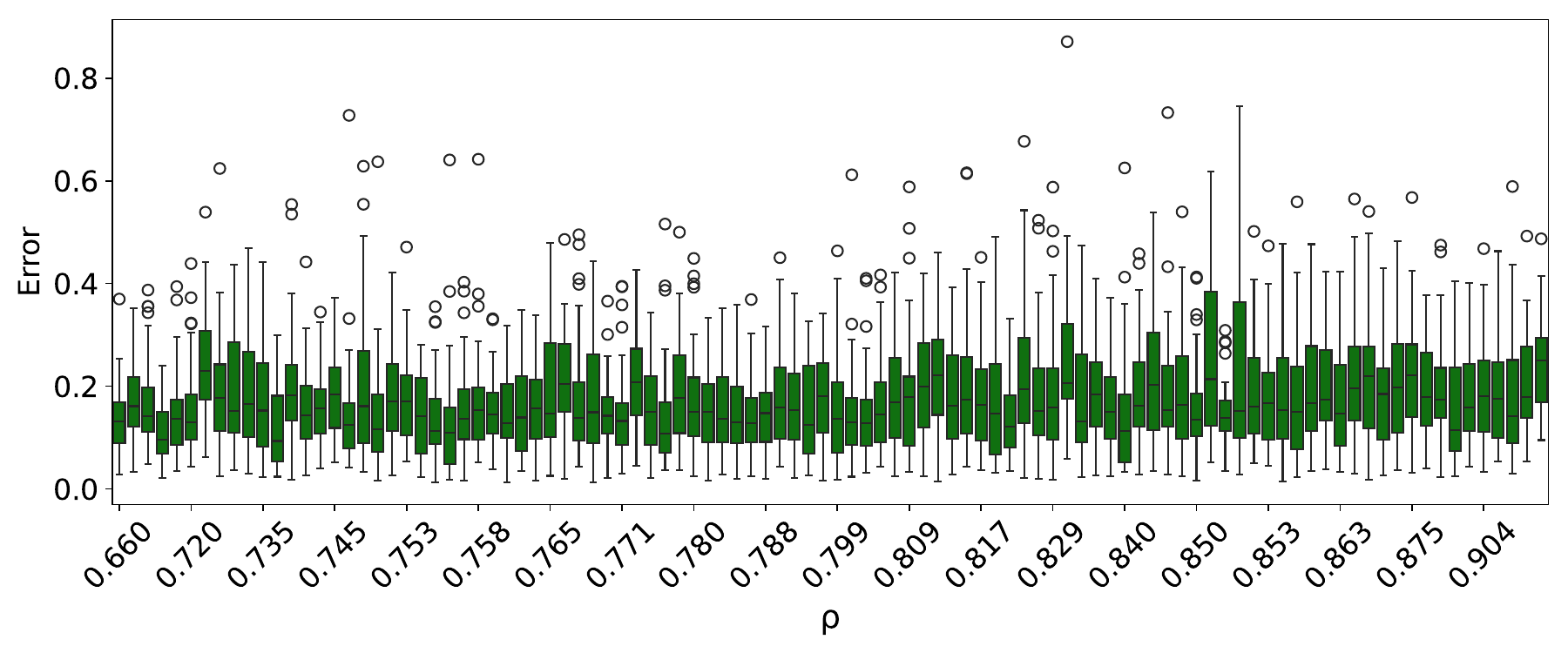}
\includegraphics[width=0.7\textwidth]{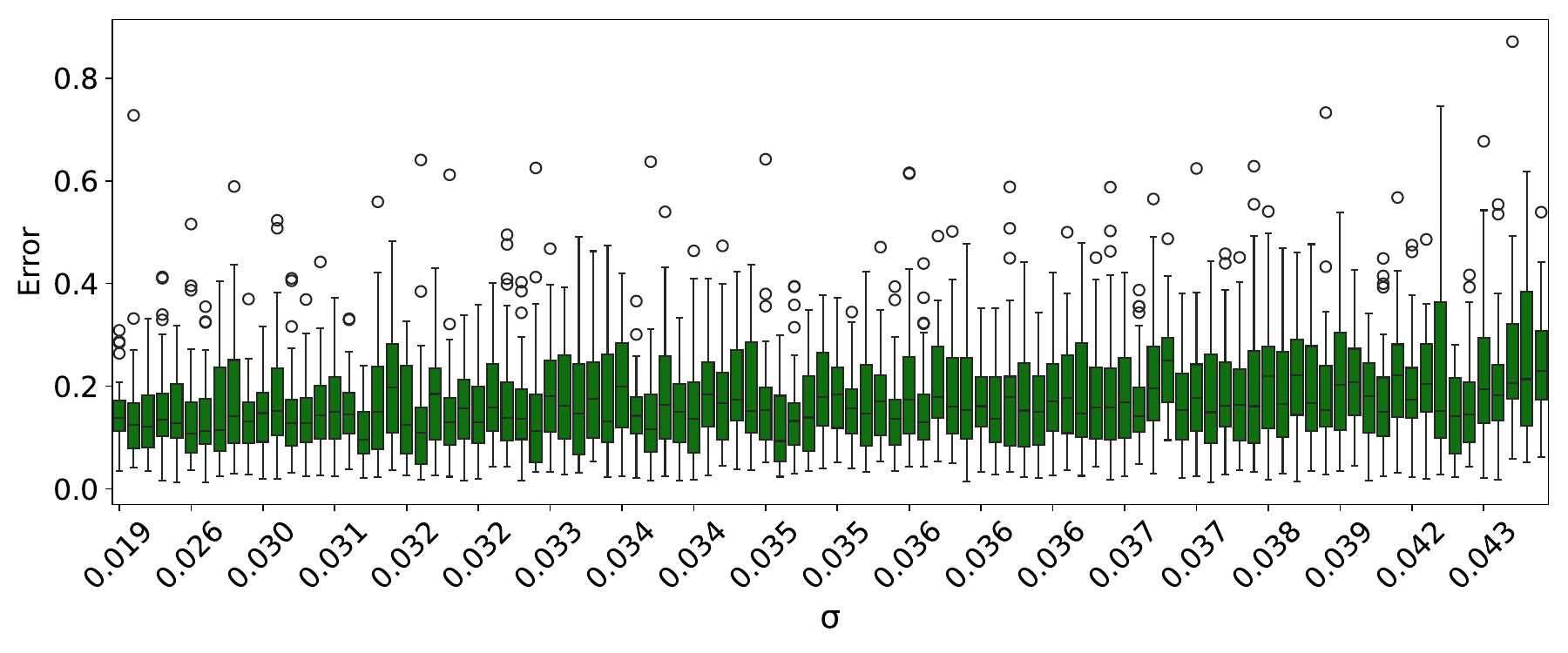}
\caption{Statistics of the relative $\ell^2$-error for the additive noise model as a function of the noise parameter $\rho$ and $\sigma$. We report the $\ell^2$-norm for the predicted entries of the model parameter vector $\vect{\theta}_{\text{dyn}}$ (top two blocks in blue) and the entries of the covariance matrix (bottom two blocks in green). The boundaries of each box correspond to the 25th and 75th percentile (lower and upper quartile), respectively. The line inside the box represents the 50th percentile. The whiskers correspond to the data points that are within 1.5 times the interquartile range from the lower and upper quartiles. The circles above or below the whiskers denote points outside this range, i.e. outliers.}\label{fig:box_whisker_additive_noise}
\end{figure}

\ipoint{Observations:} For the results for the intrinsic noise model, we can observe that increasing $\beta$ results in a slight increase (on average) in the prediction error (see \Cref{fig:box_whisker_intrinsic_noise}). We note that the errors overall remain relative stable for the wide range of values for $\beta$ considered in this study.

For the results reported for the additive noise model in \Cref{fig:box_whisker_additive_noise} we cannot observe a clear trend as each individual noise parameter increases. This can be attribute to the fact that these plots mix results for two parameters. The 2D plots included in the supplementary material seem to indicate that as both parameters increase, the errors seem to increase. But even for these plots the trend is not striking.

Lastly, these plots reaffirm that our framework struggles more when it comes to predicting the entries for the covariance matrices.

\subsection{Performance on High-Performance Computing Hardware}\label{s:hpc_performance}

The goal of these experiments is to scale up the sizes the NN architectures in terms of their trainable parameters and measure improvements in accuracy of NN predictions. Furthermore, we document the performance of NN training on high-performance computing (\iacr{HPC}) hardware. The main computational cost for NN-based parameter inference, after generating a data set for training, is the optimization of the NN parameters during the training phase.
Here, we investigate this cost on two different types of HPC hardware. We compare the contemporary graphics processing unit (GPU), NVIDIA A100, with the dedicated NN hardware, Cerebras CS-2. The Cerebras CS-2 is a system that is optimized for the workloads of training deep neural networks. More information about the two platforms can be found in \Cref{s:swhw}. The TensorFlow code that we run is identical on both platforms.

\ipoint{Setup:} The inputs and outputs of the NN are the same as in the previous experiment; see \Cref{s:covariance_results} with additive noise. Namely, we utilize \tsfc\ as NN inputs and estimate ODE model parameters $\vect{\theta}_{\text{dyn}} = (\theta_0, \theta_1, \theta_2)$, noise parameters $\vect{\theta}_{\text{noise}} = (\rho, \sigma)$ of additive noise noise, and a covariance matrix $\mat{\Gamma}_{\text{post}} \in\ns{R}^{3 \times 3}$. Hence, the NN will have the total output dimension $p=11$. The training size is $n_\text{train}=\inum{6928}$, which is the largest quantity in the previous experiment of \Cref{s:covariance_results}.

The goal of this section is to scale NN sizes to measure the effects on accuracy and the runtime performance on HPC systems. Therefore, we employ CNNs and DNNs architectures that are slightly modified compared with the ones previously used in \Cref{s:nn_exploration}. The new configurations of NNs and the settings for training are summarized in \Cref{t:hpc_nnconfig}, which lists only the differences to the baseline configurations shown in \Cref{t:nnconfig} from \Cref{s:architectures}. Notably, the dense units and the filter settings increase to up to $n_u=4096$ and $n_f=128$, respectively, amouting to up to $n_f\times2^5=2048$ filters in the CNN.  The increased complexity of the DNNs requires additional considerations for regularization and stability of the optimization.  To this end, we employ dropout layers \cite{Srivastava:2014} with probability $0.2$ between hidden layers of DNNs. Note that the dropout probability is zero for CNNs.  Furthermore, we decrease the learning rate along with an increase in optimization steps and batch size (see \Cref{t:hpc_nnconfig}).

\begin{table}
\caption{DNN and CNN configurations for scaling and performance. For brevity, this table lists only the differences relative to the baseline settings in \Cref{t:hpc_nnconfig}.}
\label{t:hpc_nnconfig}
\adjtab\centering
\begin{tabular}{lll}\toprule
\bf architecture & \bf option            & \bf value(s)                   \\
\midrule
DNN \& CNN       & $n_\text{train}$      & $\inum{6928}$                  \\
                 & batch size            & 250                            \\
                 & optimization steps    & 12800 ($\approx 462$ epochs)   \\
                 & learning rate         & $10^{-4}$                      \\
\midrule
DNN              & layers                & 12                             \\
                 & nodes/units $n_u$     & \{128,256,512,1024,2048,4096\} \\
                 & dropout               & 0.2                            \\
\midrule
CNN              & convolution layers    & 5                              \\
                 & filters setting $n_f$ & \{2,4,8,16,32,64,128\}         \\
                 & kernel size           & 3                              \\
                 & kernel stride         & 2                              \\
                 & pooling type          & \emph{none}                    \\
                 & post-flattening layers& 8 dense layers with $n_u=1024$ \\
                 & dropout               & 0                              \\
\bottomrule
\end{tabular}
\end{table}

We would ideally employ CNN architectures on the Cerebras CS-2, but the CS-2 software stack of the particular system that we have access to does not support the 1-dimensional convolutions required by our models.\footnote{Our additional efforts to find workarounds by using 2-dimensional convolutions, while not being ideal for the type of our data, have not been successful.} Therefore, we carry out performance experiments on the CS-2 with only DNN architectures, while performance with CNN and DNN architectures are shown for the GPU.

\ipoint{Results:} The results are presented in three tables. \Cref{tab:gpu_performance_cnn} for CNNs and \Cref{tab:gpu_performance_dnn} for DNNs show the number of trainable parameters and the associated accuracy on the testing data using the metrics MSE (i.e., loss function) and CDET. The best results are in bold and they are reached for $n_f=64$ for CNNs and $n_u=2048$ for DNNs.  In addition, those two tables show the runtime for training on the NVIDIA A100 GPU. We list the total training time for all optimizations steps (lower is better) and the average number of samples processed per second (higher is better). The latter is a relative metric that is independent of the number of optimizaiton steps.

\Cref{tab:cerebras_performance} shows the performance on the Cerebras CS-2 hardware for the same DNNs as in \Cref{tab:gpu_performance_dnn} making the performance on the GPU comparable to the CS-2.  The results for CS-2 encompass additional quantities. The utilization of the CS-2 chip is the percentage of programmable fabric and compute cores that are utilized to execute a program. Furthermore once before training, the CS-2 requires compilation of TensorFlow code followed by a setup, when the CS-2 fabric is configured to execute the program. Note that comilation and setup phases do not exist when running on the GPU.

\begin{table}
  \caption{Performance and scaling of CNNs on NVIDIA A100 GPU. Rows correspond to increasing NN sizes via $n_f$, while keeping fixed: 5 convolutional layers, 8 dense layers, $n_u=1024$ units per dense layer. The performance is expressed in total training time (seconds) and in a nondimensional speed: samples/second. The accuracy on the testing data (MSE, CDET) is shown when inferring $(\vect{\theta}_{\text{dyn}}, \vect{\theta}_{\text{noise}}, \mat{\Gamma}_{\text{post}})$. (For comparison, the NN used in \Cref{s:covariance_results} achieves $\text{CDET}=0.838$).}
\label{tab:gpu_performance_cnn}
\centering\adjtab
\begin{tabular}{rrrrrr}
  \toprule
  \multicolumn{2}{c}{\thead{CNN size}} & \multicolumn{2}{c}{\thead{Testing accuracy}} & \multicolumn{2}{c}{\thead{NVIDIA A100}}\\
  \cmidrule(lr){1-2}
  \cmidrule(lr){3-4}
  \cmidrule(lr){5-6}
  $n_f$ & parameters & MSE & CDET & training [sec] & samples/sec\\
  \midrule
  \inum{  2} & \inum{ 11455499} &    0.007579 &    0.798 & \decs{168.93} & \decs{18942.64}\\
  \inum{  4} & \inum{ 15551499} &    0.006571 &    0.825 & \decs{165.25} & \decs{19365.16}\\
  \inum{  8} & \inum{ 23743499} &    0.006688 &    0.822 & \decs{230.67} & \decs{13872.80}\\
  \inum{ 16} & \inum{ 40127499} &    0.006050 &    0.838 & \decs{267.10} & \decs{11980.58}\\
  \inum{ 32} & \inum{ 72895499} &    0.004958 &    0.868 & \decs{274.57} & \decs{11654.67}\\
  \inum{ 64} & \inum{138431499} &\bf 0.004578 &\bf 0.878 & \decs{418.49} & \decs{ 7646.61}\\
  \inum{128} & \inum{269503499} &    0.004968 &    0.868 & \decs{772.15} & \decs{ 4144.29}\\
  \bottomrule
\end{tabular}
\vskip 4ex
\caption{Performance and scaling of DNNs on NVIDIA A100 GPU. Rows correspond to increasing NN sizes via $n_u$, while fixing 12 dense layers. The performance is expressed in total training time (seconds) and in a nondimensional speed: samples/second. The accuracy on the testing data (MSE, CDET) is shown when inferring $(\vect{\theta}_{\text{dyn}}, \vect{\theta}_{\text{noise}}, \mat{\Gamma}_{\text{post}})$.}
\label{tab:gpu_performance_dnn}
\centering\adjtab
\begin{tabular}{rrrrrr}
  \toprule
  \multicolumn{2}{c}{\thead{DNN size}} & \multicolumn{2}{c}{\thead{Testing accuracy}} & \multicolumn{2}{c}{\thead{NVIDIA A100}}\\
  \cmidrule(lr){1-2}
  \cmidrule(lr){3-4}
  \cmidrule(lr){5-6}
  $n_u$ & parameters & MSE & CDET & training [sec] & samples/sec\\\midrule
  \inum{ 128} & \inum{   695179} &    0.009745 &    0.740 & \decs{128.74} & \decs{24906.5}\\
  \inum{ 256} & \inum{  1750795} &    0.006169 &    0.836 & \decs{130.67} & \decs{24532.8}\\
  \inum{ 512} & \inum{  4943371} &    0.005819 &    0.845 & \decs{130.58} & \decs{24550.7}\\
  \inum{1024} & \inum{ 15653899} &    0.005687 &    0.848 & \decs{132.47} & \decs{24178.3}\\
  \inum{2048} & \inum{ 54376459} &\bf 0.005456 &\bf 0.855 & \decs{151.06} & \decs{21205.2}\\
  \inum{4096} & \inum{201027595} &    0.005631 &    0.850 & \decs{284.88} & \decs{11238.8}\\
  \bottomrule
\end{tabular}
\vskip 4ex
\caption{Performance and scaling of DNNs on Cerebras CS-2. We increase utilization of the CS-2 chip from 1\% to 94\% while maintaining a high throughput that only weakly depends on the NN size (see samples/sec), because it decreases by only $\sim$25\%. Rows correspond to increasing NN sizes via $n_u$, while fixing 12 dense layers. Explanations for utilization and compile, setup, and training times are given in the text. The CS-2 processes 8--14.5 times more samples per second compared with the A100 GPU (see \Cref{tab:gpu_performance_dnn}).}
\label{tab:cerebras_performance}
\centering\adjtab
\begin{tabular}{rrrrrrr}
  \toprule
  \multicolumn{2}{c}{\thead{DNN size}} & \multicolumn{5}{c}{\thead{Cerebras CS-2}}\\
  \cmidrule(lr){1-2}
  \cmidrule(lr){3-7}
  $n_u$ & parameters & utilization\,[\%] & compile [sec] & setup [sec] & training [sec] & samples/sec\\
  \midrule
  \inum{ 128} & \inum{   695179} & \decs{ 1.08} & \decs{ 204.67} & \decs{ 53.66} & \decs{31.07} & \decs{215375.0}\\
  \inum{ 256} & \inum{  1750795} & \decs{ 2.10} & \decs{ 234.82} & \decs{ 55.06} & \decs{31.27} & \decs{214106.3}\\
  \inum{ 512} & \inum{  4943371} & \decs{ 4.89} & \decs{ 228.18} & \decs{ 59.64} & \decs{32.16} & \decs{215527.1}\\
  \inum{1024} & \inum{ 15653899} & \decs{12.69} & \decs{ 279.07} & \decs{ 73.21} & \decs{30.52} & \decs{201337.7}\\
  \inum{2048} & \inum{ 54376459} & \decs{39.64} & \decs{1132.18} & \decs{125.92} & \decs{34.04} & \decs{197587.8}\\
  \inum{4096} & \inum{201027595} & \decs{94.10} & \decs{1445.67} & \decs{250.82} & \decs{51.82} & \decs{162383.4}\\
  \bottomrule
\end{tabular}
\end{table}

\ipoint{Observations:} We first note the effect of scaling NN sizes on the accuracy of predictions. We grow the number of trainable parameters exponentially and are able to reach a pivotal point in prediction accuracy for both CNNs and DNNs. The accuracy is better with CNNs by around 5\% compared with DNNs, which supports our observations in \Cref{s:nn_exploration}.  On the flipside, is is significantly faster to train DNNs: for the best models of each class of architecture, the DNN ($n_u=2048$ in \Cref{tab:gpu_performance_dnn}) takes around 1/3 of the training time as the CNN ($n_f=64$ in \Cref{tab:gpu_performance_cnn}).

The training time on the GPU can depend on the NN size. For the ``smaller'' NN sizes with $<10^5$ parameters, the time is nearly constant.  It only increases for sufficiently large networks (e.g., CNNs with $n_f \ge 8$ and DNNs with $n_u \ge 2048$). The effect of steady training time despite larger NNs is even more impressive on the CS-2, because the training time decreases only slightly as the DNN size grows exponentially to more than 50 million parameters. It is only near the limit of CS-2 utilization, when the training slows by around 40\%  ($n_u=2048$ vs.\ $n_u=4096$ in \Cref{tab:cerebras_performance}). The training performance when considering samples per second on Cerebras CS-2 is better than NVIDIA A100 by a factor of $\sim$8 when DNNs have $n_u=128,\ldots,512$ units. As the DNN sizes increase further, the gap between the two systems widens to a factor of $\sim$14.5 at $n_u=\inum{4096}$. Currently, the speed of computations of the CS-2 comes with an overhead cost for program compilation and setup of the programmable hardware. Additional limitations are caused by the less mature software stack compared with established GPU hardware.  If these limitations can be lifted by the vendor in the future, we see greater potential for using such specialized hardware for scientific applications.

\section{Conclusions}\label{s:conclusions}

We have presented a study of the performance of NNs for inference of model parameters, noise parameters, and approximate parameter uncertainties in inverse parameter estimation problems governed by nonlinear ODEs. We chose this problem as a testing bed since it poses significant challenges for traditional variational approaches or sampling methods for Bayesian inference; the associated optimization landscape features sharp gradients, a narrow, flat optimality zone, and is strongly nonconvex. In addition, this nonlinear dynamical system is simple enough to gain greater insight into the performance of our approach than more complex systems would allow. Moreover, the forward and adjoint operators are relatively cheap to evaluate and manipulate. This provides an excellent benchmark for the proposed NN framework for statistical inference.

We reported results for a joint estimation of model parameters and noise parameters. We considered different noise models; an autocorrelated additive noise model, an intrinsic noise model (based on an SDE), as well as a combination thereof. In addition, we developed a framework that allows us to estimate the covariance of the posterior distribution of the model parameters conditioned on the data. To provide training data for the NN, we have \begin{inparaenum}[\it (i)] \item  used the prior of the parameters and \item exploited a Laplace approximation that allows us to efficiently generate approximations of the covariance of the posterior distribution\end{inparaenum}. We have not only considered different noise settings but also explored the performance of the considered NNs for different features (time series, spectral coefficients, and a combination thereof). We considered a log-Euclidean framework to project the covariance data from a Riemannian manifold into an Euclidean space in an attempt to simplify our computations and guarantee that the estimated covariances are symmetric positive definite matrices.

The \textit{most important observations} of our experiments are as follows:
\begin{itemize}[leftmargin=*]
\item
We obtain an excellent agreement (i.e., relative errors are at single-digit percentages) between NN predictions and the true, underlying model parameters, noise parameters, and covariance entries for an extremely challenging problem.
\item
The performance of our methodology is robust with respect to different noise settings in observational data of the inverse problem.
\item
To reliably estimate model and noise parameters we require a combination of spectral coefficients and time series data.
\item
A minimum training size of \inum{2000} training samples is required to achieve what we consider acceptable results.
\item
The NN framework scales with increasing output dimension, where accuracy remains high without significant increase in training data.
\item
We could reduce the training time by a factor of 8--14 when utilizing dedicated hardware platforms for training NNs.
\end{itemize}

In this work, we have limited the our explorations to DNN and CNN architectures. These architectures performed well for our purposes. In future work, we intend to push the proposed methodology to large-scale inverse problems, for which the forward and adjoint operators are more challenging to evaluate. In this context, we plan to explore modern ML frameworks that have gained popularity in the inverse problems community. Potential approaches include generative adversarial networks~\cite{adler2024:deep, Patel:2022a, Ray:2023}, normalizing flows~\cite{Rezende:2015, papamakarios2021:normalizing, dasgupta2024:dimension}, or diffusion models~\cite{dasgupta2025:conditional,daras2024:survey, chung2022:improving}.

\begin{appendix}

\section{Software and Hardware Environment}\label{s:swhw}

The numerical experiment in \Cref{s:model_noise_results} and \Cref{s:covariance_results} are implemented in \texttt{Python} \texttt{3.7.10}. We utilize the following software environment:
\begin{inparaenum}[\it (i)]
  \item \texttt{Numpy} \texttt{1.18.5},
  \item \texttt{Scipy} \texttt{1.4.1},
  \item \texttt{Sklearn} \texttt{0.24.2},
  \item \texttt{Tensorflow} \texttt{2.3.0}.
\end{inparaenum} We execute these runs on UH's Sabine system. Each node is equipped with a 28 core Intel Xeon E5-2680v4 CPU with 128 or 256\, GB memory.

The GPU system used in \Cref{s:hpc_performance} is hosted by the Department of Mathematics at Virginia Tech. The GPU is the NVIDIA A100 connected via PCIe; it contains 108 Multiprocessors and 64 CUDA Cores per Multiprocessor (6912 CUDA Cores in total), 40GB of device memory and a memory bandwidth of 1.56TB/s. The processor has 54.2 billion transistors over a silicon area of \inum{826}mm\textsuperscript{2}.  The software environment is
\begin{inparaenum}[\it (i)]
  \item \texttt{CUDA 11},
  \item \texttt{cuDNN 8.9.6.50},
  \item \texttt{cuBLAS 11.11.3.6},
  \item \texttt{TensorFlow-GPU 2.9.3}.
\end{inparaenum}

The Cerebras CS-2 system used in \Cref{s:hpc_performance} is hosted at the Pittsburgh Supercomputing Center within their Neocortex system, which houses two Cerebras CS-2 systems. The CS-2 main component is the processor: Cerebras Wafer Scale Engine 2 (WSE2). The WSE2 consists of \inum{850000} compute cores in total that are interconnected via a so-called programmable fabric, where each core contains a programmable router to communicate with its four neighboring cores. The WSE2 has 40GB of on-chip memory with a memory bandwidth of 20PB/s. The processor has 2.6 trillion transistors over a silicon area of \inum{46225}mm\textsuperscript{2}.  The software environment is
\begin{inparaenum}[\it (i)]
  \item \texttt{Cerebras Model Zoo 1.6.0},
  \item \texttt{TensorFlow 2.2.0}.
\end{inparaenum}

\section{Supplementary Material}

The supplementary materials provide additional background information on the computation of the reduced-space Hessian and complementary results to those reported in the main manuscript. In \Cref{s:curvature_supp}, we discuss the reduced-space Hessian. First, we outline the steps required to evaluate the Hessian matvec in \Cref{s:comphess_supp}. We formally derive the reduced-space Hessian in \Cref{s:variations_supp}. We address issues with ill-conditioning, negative definite, and indefinite Hessians in the training data in \Cref{s:cov_instabilities_supp}. We report additional numerical results in \Cref{s:supplementary_results}.
In particular, we provide additional results for exploring the NN architectures in \Cref{s:exploration_nns_addres_supp}
We study the effects of including or omitting noise in the training data in \Cref{s:effects_of_noise_training_supp}.
We report results for estimating noise parameters for the intrinsic noise in \Cref{s:noise_model_para_intrinsic_supp}.
In \Cref{s:predicted_state_supp} we explore if the end-to-end predictions of the parameters of the dynamical system yield state predictions that are in on good agreement with the data and the true state of the system.
We discuss issues with an off-the-shelf prediction of covariance matrices in \Cref{s:pred_npsd_mat_fix_supp}.
We provide additional results for the evaluation of the sensitivity of the NN predictions with respect to the initialization of the NN in \Cref{s:nn_initialization_supp}.
We include $k$-fold cross-validation results in \Cref{s:nn_initialization_supp_cc} to shed additional light on how our framework generalizes to unseen data.
We provide additional results for the sensitivity of our NN framework with respect to the amount of noise perturbations in \Cref{s:noise_sensitivity_supp}.

\subsection{Curvature Information}\label{s:curvature_supp}

As outlined in \Cref{s:data}, we use the Hessian $\mat{H} \in \spd{3}$ associated with the variational problem formulation in~\cref{e:varopt} to approximate the covariance matrix $\mat{\Gamma}_{\text{post}} \in \spd{3}$ of the posterior distribution. In the following, we discuss the computations for the construction of the Hessian.

We consider an \emph{optimize-then-discretize} approach to derive the second-order derivatives~\cite{Gunzburger:2003a}. That is, we compute first- and second-order variations in the continuum and subsequently discretize the resulting mathematical expressions. To handle the ODE constraints in~\cref{e:varopt}, we introduce the Lagrange multipliers $\lambda : [0,\tau] \to \ns{R}$ for~\cref{e:stateu} and $\nu : [0,\tau]\to\ns{R}$ for~\cref{e:statev}. Let $\vect{\Xi}$ denote the collection of the parameter vector, and the primal and dual variables,  i.e., $\vect{\Xi} \defeq (\vect{\theta}_{\text{dyn}},u,v,\lambda,\nu)$. The Lagrangian functional $\ell$ is given by
\begin{equation}\label{e:lagrangian-fhn}
\begin{aligned}
\ell(\vect{\Xi})
&= \frac{\gamma}{2}\int_0^\tau (u(t) - y_{\text{obs}}(t))^2 \,\d t
+ \frac{1}{2} \|\vect{\theta}_{\text{dyn}} - \vect{\theta}_{\text{ref}}\|_{\mat{L}}^2 \\
& \quad + \int_0^1 \lambda (\d_t u - \theta_2(u - (u^3/3) + v + z)) \,\d t
+ \lambda(t=0) (u(t=0) - u_0) \\
& \quad + \int_0^1(\d_t v + (u - \theta_0 - \theta_1 v)/\theta_2)\nu \, \d t
+ \nu(t=0)(v(t=0) - v_0),
\end{aligned}
\end{equation}

\noindent where $\gamma > 0$, $\mat{L} \in \spd{3}$, $\vect{\theta}_{\text{ref}} \in \ns{R}^3$, $z > 0$, $u_0 \in \ns{R}$ and $v_0 \in \ns{R}$ are fixed variables defined in \Cref{s:methods}.

Below, we denote the perturbations of the parameter vector, and the primal and dual variables associated with the first variations by $\vect{\hat{\Xi}} \defeq (\vect{\hat{\theta}},\hat{u},\hat{v},\hat{\lambda},\hat{\nu})$. For the second variations, we use $\vect{\tilde{\Xi}} \defeq (\vect{\tilde{\theta}},\tilde{u},\tilde{v},\tilde{\lambda},\tilde{\nu})$. We refer to $\vect{\tilde{\theta}}$ as the incremental parameter vector and to $\tilde{u}$, $\tilde{v}$, $\tilde{\lambda}$, and $\tilde{\nu}$, as the incremental state and adjoint variables. In the context of optimization, these variables correspond to the search directions associated with the Newton step.

In \Cref{s:comphess_supp} we overview the steps associated with constructing the Hessian matrix $\mat{H}$. In \Cref{s:variations_supp} we formally derive the first- and second-order variations.

\subsubsection{Computation of Hessian}\label{s:comphess_supp}

We consider a reduced space method to compute the Hessian~\cite{Biros:2005a,Biros:2005b}. The Hessian matvec, i.e., the action of the Hessian on a vector $\vect{\tilde{\theta}}$, is given by
\begin{equation}\label{e:hessianmatvec}
\mathcal{H}\vect{\tilde{\theta}}
=
(\mathcal{H}_{\text{reg}} + \mathcal{H}_{\text{dat}})\vect{\tilde{\theta}}
\defeq
\mat{L}\vect{\tilde{\theta}}
-\int_0^\tau
\mat{Q}
\begin{pmatrix}
\tilde{\lambda} \\
\tilde{\nu}
\end{pmatrix}
+
\mat{\tilde{Q}}
\begin{pmatrix}
\lambda \\
\nu
\end{pmatrix}
\d t,
\end{equation}

\noindent  where $\mathcal{H}_{\text{reg}} = \mat{L} \in \spd{3}$ is the contribution of the regularization model and
\[
\mathcal{H}_{\text{dat}} = -\int_0^\tau
\mat{Q}
\begin{pmatrix}
\tilde{\lambda} \\
\tilde{\nu}
\end{pmatrix}
+
\mat{\tilde{Q}}
\begin{pmatrix}
\lambda \\
\nu
\end{pmatrix}
\d t
\]

\noindent represents the data misfit term. The expressions for $\mat{Q}$ and $\mat{\tilde{Q}}$ are given by
\[
\mat{Q} =
\begin{pmatrix}
0                 & 1/\theta_2 \\
0                 & -v/\theta_2 \\
u-(u^3/3) + v + z & (u-\theta_0 +\theta_1 v)/\theta_2^2
\end{pmatrix}
\]

\noindent and
\[
\mat{\tilde{Q}} =
\begin{pmatrix}
0 & -\tilde{\theta}_2/ \theta_2^2 \\
0 & -(\tilde{v} / \theta_2) + (v\tilde{\theta}_2/\theta_2^2) \\
(1-u^2)\tilde{u}+\tilde{v} &
((
\tilde{u} + \theta_1 \tilde{v}
+v\tilde{\theta}_1 - \tilde{\theta}_0
)
/ \theta_2^2)
- (2(u-\theta_0 +\theta_1 v)\tilde{\theta}_2/\theta_2^3)
\end{pmatrix},
\]

\noindent respectively. Since the discrete representation $\mat{H}$ of $\mathcal{H}$ is a $3\times 3$ matrix, we can explicitly form and store $\mat{H}$. To do so, we sample the columns of $\mat{H}$ by applying the discretized version of $\mathcal{H}$ in~\cref{e:hessianmatvec} to the standard basis vectors $\vect{e}_k \in \mathbb{R}^3$, $k = 1,2,3$. That is, we invoke~\cref{e:hessianmatvec} by replacing $\vect{\tilde{\theta}}$ with $\vect{e}_k$. To evaluate~\cref{e:hessianmatvec} we require all variables of $\vect{\Xi}$ and $\vect{\tilde{\Xi}}$, respectively. The candidate vector $\vect{\theta}_{\text{dyn}} \in \ns{R}^3$ is given. We find the state variables $u:[0,\tau] \to \ns{R}$ and $v:[0,\tau] \to \ns{R}$ for the candidate $\vect{\theta}_{\text{dyn}}$ by solving \cref{e:fhnfwd} forward in time. Subsequently, we compute $\lambda:[0,\tau] \to \ns{R}$ and $\nu:[0,\tau] \to \ns{R}$ by solving the adjoint equations
\begin{subequations}
\label{e:adj}
\begin{align}
-\d_t\lambda - \theta_2 (1-u^2) \lambda + (\nu/\theta_2) + u - y_{\text{obs}} &= 0
&& \text{in}\,\,[0,\tau),
\\
-\d_t\nu + (\theta_1 \nu/\theta_2) - \theta_2 \lambda & = 0
&& \text{in}\,\,[0,\tau),
\end{align}
\end{subequations}

\noindent with final conditions $\lambda(t=\tau) = 0$ and $\nu(t=\tau) = 0$, respectively, backward in time. The incremental state variables $\tilde{u}:[0,\tau] \to \ns{R}$ and $\tilde{v}:[0,\tau] \to \ns{R}$ are found by solving
\[
\begin{aligned}
\d_t \tilde{u}
- \theta_2((1-u^2)\tilde{u} + \tilde{v}) - \tilde{\theta}_2(u-(u^3/3) + v + z) &= 0
&& \text{in}\,\,(0,\tau],
\\
 \d_t \tilde{v}
+ [\theta_1\tilde{v}
+ (\tilde{u}
- \tilde{\theta}_0
+ \tilde{\theta}_1v)
- (\tilde{\theta}_2(u - \theta_0 + \theta_1 v)/\theta_2)
]
/\theta_2
&=0
&& \text{in}\,\,(0,\tau],
\end{aligned}
\]

\noindent with initial conditions $\tilde{u}(t=0) = 0$ and $\tilde{v}(t=0) = 0$, respectively, forward in time. The incremental adjoint variables $\tilde{\lambda}:[0,\tau] \to \ns{R}$ and $\tilde{\nu}:[0,\tau] \to \ns{R}$ are found by solving
\[
\begin{aligned}
-\d_t \tilde{\lambda}
- \theta_2(1-u^2)\tilde{\lambda}
- \tilde{\theta}_2(1-u^2)\lambda
+ (\tilde{\nu}/\theta_2)
-(\nu\tilde{\theta}_2/\theta_2^2)
+ (1+2\theta_2 u \lambda) \tilde{u}
&=0
&& \text{in}\,\,[0,\tau),
\\
- \d_t\tilde{\nu}
+
(
[
\theta_1\tilde{\nu}
+ \tilde{\theta}_1 \nu
- (\theta_1 \tilde{\theta}_2\nu/\theta_2)
]
/\theta_2
)
- \theta_2 \tilde{\lambda}
- \tilde{\theta}_2\lambda
&= 0
&& \text{in}\,\,[0,\tau),
\end{aligned}
\]

\noindent with final conditions $\tilde{\lambda}(t=\tau) = 0$ and $\tilde{\nu}(t=\tau) = 0$, respectively, backward in time.

\subsubsection{Variations}\label{s:variations_supp}

Below, we formally derive the first- and second-order variations. We assume that all variables that appear in the derivations below satisfy suitable regularity requirements to carry out the necessary computations.

The first variations with respect to the perturbations $\hat{\lambda}$, $\hat{\nu}$, of the Lagrange multipliers $\lambda$, $\nu$, are given by
\begin{align*}
\ell_{\lambda}[\vect{\Xi}](\hat{\lambda}) &=
\int_0^\tau (\d_t u - \theta_2(u - (u^3/3) + v + z))  \hat{\lambda}\,\d t
+ (u(t=0) - u_0)\hat{\lambda}(t=0) ,
\\
\ell_{\nu}[\vect{\Xi}](\hat{\nu}) &=
\int_0^\tau (\d_t v + (u - \theta_0 - \theta_1 v)/\theta_2) \hat{\nu}\, \d t
+ (v(t=0) - v_0)\hat{\nu}(t=0).
\end{align*}

The variations with respect to the perturbations $\hat{u}$, $\hat{v}$, of the state variables $u$, $v$, are given by
\begin{align*}
\ell_u[\vect{\Xi}](\hat{u}) &=
\int_0^\tau (-\d_{t}\lambda - \theta_2 (1-u^2) \lambda + (\nu/\theta_2) + u - y_{\text{obs}})\hat{u}\,\d t
+ \lambda(t=\tau)\hat{u}(t=\tau),
\\
\ell_v[\vect{\Xi}](\hat{v}) &=
\int_0^\tau(-\d_t\nu + (\theta_1 \nu/\theta_2) - \theta_2 \lambda)\hat{v}\,\d t
+ \nu(\tau=1)\hat{v}(\tau=1).
\end{align*}

The variations with respect to the perturbations $\vect{\hat{\theta}} = (\hat{\theta}_0, \hat{\theta}_1, \hat{\theta}_2) \in \ns{R}^3$ of the parameter vector $\vect{\theta}_{\text{dyn}} = (\theta_0, \theta_1, \theta_2) \in \ns{R}^3$ are given by
\[
\begin{aligned}
\ell_{\vect{\theta}}[\vect{\Xi}](\hat{\vect{\theta}}) &=
\langle\mat{L}\vect{\theta}_{\text{dyn}},\vect{\hat{\theta}}\rangle
\!
-
\!\!\int_0^\tau
\!\!
\Big\langle
\!\!
\begin{pmatrix}
0 &  1/\theta_2 \\
0 & -v/\theta_2 \\
u-(u^3/3) + v + z & (u-\theta_0 +\theta_1 v)/\theta_2^2
\end{pmatrix}\!\!
\begin{pmatrix}
\lambda \\
\nu
\end{pmatrix}\!,
\vect{\hat{\theta}}
\Big\rangle\, \d t.
\end{aligned}
\]

At optimality, i.e., for any primal-dual optimal $\vect{\Xi}^\star = (u^\star,v^\star,\nu^\star,\lambda^\star,\vect{\theta}_{\text{dyn}}^\star)$, these variations vanish. Setting $\ell_{\lambda}[\vect{\Xi}](\hat{\lambda})$, $\ell_{\nu}[\vect{\Xi}](\hat{\nu})$, $\ell_u[\vect{\Xi}](\hat{u})$ and $\ell_v[\vect{\Xi}](\hat{v})$ to zero yields the strong form of the state and adjoint equations for the variational problem~\cref{e:varopt} (see \cref{e:fhnfwd} and \cref{e:adj}). Their numerical time integration allows us to compute the state variables $u$ and $v$ as well as the adjoint variables $\lambda$ and $\nu$ for a candidate $\vect{\theta}_{\text{dyn}}$.

Since we are interested in deriving curvature information for training the NNs for the prediction of the posterior covariance matrices we have to compute second variations. These follow from the first variations in a straight-forward way.

The second variations of the Lagrangian in with respect to perturbations $\tilde{\lambda}$ in $\lambda$ are given by
\[
\begin{aligned}
\ell_{\lambda\lambda}[\vect{\Xi},\vect{\tilde{\Xi}}](\hat{\lambda})
& = 0,
\\
\ell_{\nu\lambda}[\vect{\Xi},\vect{\tilde{\Xi}}](\hat{\nu})
& = 0,
\\
\ell_{u \lambda}[\vect{\Xi},\vect{\tilde{\Xi}}](\hat{u})
& = \int_0^\tau
(- \d_t\tilde{\lambda} - \theta_2(1-u^2)\tilde{\lambda})  \hat{u}\, \d t
+ \tilde{\lambda}(t=\tau)\hat{u}(t=\tau),
\\
\ell_{v \lambda}[\vect{\Xi},\vect{\tilde{\Xi}}](\hat{v})
& = \int_0^\tau (- \theta_2 \tilde{\lambda})  \hat{v}\, \d t,
\\
\ell_{\vect{\theta}\lambda}[\vect{\Xi},\vect{\tilde{\Xi}}](\vect{\hat{\theta}})
& = \int_0^\tau - (u - (u^3/3)+v+z) \tilde{\lambda}\hat{\theta}_2\, \d t.
\end{aligned}
\]

The second variations of the Lagrangian in with respect to perturbations $\tilde{\nu}$ in $\nu$ are given by
\[
\begin{aligned}
\ell_{\lambda\nu}[\vect{\Xi},\vect{\tilde{\Xi}}](\hat{\lambda})
& = 0,
\\
\ell_{\nu\nu}[\vect{\Xi},\vect{\tilde{\Xi}}](\hat{\nu})
& = 0,
\\
\ell_{u\nu}[\vect{\Xi},\vect{\tilde{\Xi}}](\hat{u})
&= \int_0^\tau (\tilde{\nu}/\theta_2)\hat{u}\, \d t,
\\
\ell_{v\nu}[\vect{\Xi},\vect{\tilde{\Xi}}](\hat{v})
&= \int_0^\tau (- \d_t\tilde{\nu} + (\theta_1\tilde{\nu}/\theta_2)) \hat{v}\, \d t
+ \tilde{\nu}(t=\tau)\hat{v}(t=\tau),
\\
\ell_{\vect{\theta}\nu}[\vect{\Xi},\vect{\tilde{\Xi}}](\vect{\hat{\theta}})
&= \int_0^\tau
(-\tilde{\nu}/\theta_2)\hat{\theta}_0
+(v\tilde{\nu}/\theta_2)\hat{\theta}_1
-((u-\theta_0 +\theta_1 v)\tilde{\nu}/\theta_2^2)\hat{\theta}_2
\, \d t.
\end{aligned}
\]

The second variations of the Lagrangian in with respect to perturbations $\tilde{u}$ in $u$ are given by
\[
\begin{aligned}
\ell_{\lambda u}[\vect{\Xi},\vect{\tilde{\Xi}}](\hat{\lambda})
&= \int_0^\tau (\d_t \tilde{u} - \theta_2 (1-u^2)\tilde{u})\hat{\lambda}\, \d t
+  \tilde{u}(t=0)\hat{\lambda}(t=0),
\\
\ell_{\nu u}[\vect{\Xi},\vect{\tilde{\Xi}}](\hat{\nu})
&= \int_0^\tau (\tilde{u}/\theta_2)\hat{\nu} \, \d t,
\\
\ell_{u u}[\vect{\Xi},\vect{\tilde{\Xi}}](\hat{u})
&=\int_0^\tau (1+2\theta_2 u \lambda) \tilde{u} \hat{u}\,\d t,
\\
\ell_{v u}[\vect{\Xi},\vect{\tilde{\Xi}}](\hat{v})
&= 0,
\\
\ell_{\vect{\theta} u}[\vect{\Xi},\vect{\tilde{\Xi}}](\vect{\hat{\theta}})
&= \int_0^\tau -((1-u^2)\lambda\tilde{u} + (\nu\tilde{u}/\theta_2^2)) \hat{\theta}_2\, \d t.
\end{aligned}
\]

The second variations of the Lagrangian in with respect to perturbations $\tilde{v}$ in $v$ are given by
\[
\begin{aligned}
\ell_{\lambda v}[\vect{\Xi},\vect{\tilde{\Xi}}](\hat{\lambda})
&=
\int_0^\tau - \theta_2 \tilde{v}\hat{\lambda}\, \d t,
\\
\ell_{\nu v}[\vect{\Xi},\vect{\tilde{\Xi}}](\hat{\nu})
&=
\int_0^\tau (\d_t \tilde{v} + (\theta_1\tilde{v}/\theta_2))\hat{\nu}\, \d t
+ \tilde{v}(t=0)\hat{\nu}(t=0),
\\
\ell_{u v}[\vect{\Xi},\vect{\tilde{\Xi}}](\hat{u})
&= 0,
\\
\ell_{v v}[\vect{\Xi},\vect{\tilde{\Xi}}](\hat{v})
&= 0,
\\
\ell_{\vect{\theta} v}[\vect{\Xi},\vect{\tilde{\Xi}}](\vect{\hat{\theta}})
&=
\int_0^\tau
(\nu\tilde{v}/\theta_2) \hat{\theta}_1
 - (\lambda\tilde{v} + (\theta_1\nu\tilde{v}/\theta_2^2))\hat{\theta}_2\, \d t.
\end{aligned}
\]

The second variations of the Lagrangian in with respect to perturbations $\vect{\tilde{\theta}} = (\tilde{\theta}_0, \tilde{\theta}_1, \tilde{\theta}_2) \in \ns{R}^3$ in $\vect{\theta}_{\text{dyn}} = (\theta_0, \theta_1, \theta_2) \in \ns{R}^3$ are given by
\[
\begin{aligned}
\ell_{\lambda \vect{\theta}}[\vect{\Xi},\vect{\tilde{\Xi}}](\hat{\lambda})
&=
\int_0^\tau -(u-(u^3/3) +v + z)\tilde{\theta}_2 \hat{\lambda} \, \d t,
\\
\ell_{\nu \vect{\theta}}[\vect{\Xi},\vect{\tilde{\Xi}}](\hat{\nu})
&=
\int_0^\tau
[
(-(u - \theta_0+ \theta_1 v)\tilde{\theta}_2/\theta_2^2)
- (\tilde{\theta}_0 / \theta_2)
+ (v\tilde{\theta}_1/\theta_2)
]\hat{\nu} \, \d t,
\\
\ell_{u \vect{\theta}}[\vect{\Xi},\vect{\tilde{\Xi}}](\hat{u})
&=
\int_0^\tau
[
-(1-u^2)\lambda\tilde{\theta}_2
-(\nu\tilde{\theta}_2/\theta_2^2)
]\hat{u} \, \d t,
\\
\ell_{v \vect{\theta}}[\vect{\Xi},\vect{\tilde{\Xi}}](\hat{v})
&=
\int_0^\tau
[
-\tilde{\theta}_2 \lambda
- (\theta_1 \tilde{\theta}_2\nu/\theta_2^2)
+ (\tilde{\theta}_1 \nu/\theta_2)
]\hat{v} \, \d t,
\\
\ell_{\vect{\theta}\vect{\theta}}[\vect{\Xi},\vect{\tilde{\Xi}}](\vect{\hat{\theta}})
&=
\langle\mat{L} \vect{\tilde{\theta}},\vect{\hat{\theta}}\rangle
+
\int_0^\tau
\Big\langle
\frac{\nu}{\theta_2^2}
\begin{pmatrix}
\tilde{\theta}_2\\
-v\tilde{\theta}_2\\
(2(u-\theta_0 + \theta_1 v)\tilde{\theta}_2/\theta_2)
+\tilde{\theta}_0
-v\tilde{\theta}_1
\end{pmatrix}
,\vect{\hat{\theta}}
\Big\rangle \,\d t.
\end{aligned}
\]

\subsection{Instabilities in Covariance Dataset}\label{s:cov_instabilities_supp}

We generate $\inum{15000}$ samples for the model parameter $\vect{\theta}_{\text{dyn}}\in\ns{R}^3$. Since we evaluate the Hessian $\mat{H}$ at the solution of our problem (i.e., for a given $\vect{\theta}_{\text{dyn}}$ that was used to generate the observation $\vect{y}_{\text{obs}}$), we expect that $\mat{H} \in \spd{3}$. However, we do not necessarily observe this in computation. Due to numerical errors, the Hessians $\mat{H}$ we compute is only symmetric up to a certain level. That is, $\|\mat{H} - \mat{H}^{\mathsf{T}}\|$ is not exactly zero. As a remedy, we set $\mat{H}$ to $(\mat{H} + \mat{H}^{\mathsf{T}})/2$. In addition, we observed that for certain parameters $\vect{\theta}_{\text{dyn}}$ the associated Hessian is no longer positive definite. Because the number of negative definite matrices is relatively small (110 samples), we decided to exclude these examples for the experiments reported in \Cref{s:covariance_results} of the main manuscript.

\begin{figure}
\centering
\includegraphics[width=0.8\textwidth]{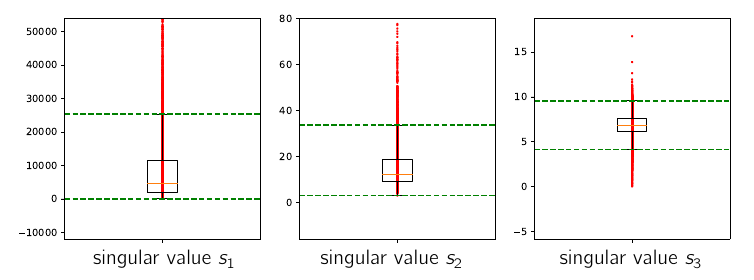}
\caption{Singular values $s_i$, $i=1,2,3$, for the positive-definite Hessian matrices $\mat{H}$. From left to right we show the box-whisker plots for the singular values $s_1$, $s_2$, and $s_3$, respectively. The original dataset consists of \inum{15000} samples. 110 of these matrices were negative definite in computation. Here, we show the singular values for the remaining \inum{14890} matrices. Dashed horizontal lines correspond to the upper and lower whiskers of the plots.}
\label{f:sing_vals_hess3D}
\end{figure}

We also observed that some of the Hessian matrices are severely ill-conditioned. This is a significant challenge, since we use the inverse of the Hessian to approximate the covariance $\mat{\Gamma}_{\text{post}}$ of the negative log posterior. Let $\mat{H} = \mat{U}\mat{S}\mat{V}^{\mathsf{T}}$ with $\mat{S} = \diag(s_1, s_2, s_3)$, $s_1 \geq s_2 \geq s_3 > 0$. We observed that the first singular values $s_1$ associated with a subset of the training samples for $\vect{\theta}_{\text{dyn}}$ were extremely large. We show box-whisker plots for the singular values $s_i$, $i = 1,2,3$, in \Cref{f:sing_vals_hess3D}, with the green horizontal lines further highlighting the whiskers of the plots. The upper whiskers are computed by adding 1.5 times the interquartile range (IQR) to the 75th percentile, while the lower whiskers are determined by subtracting 1.5 times the IQR from the 25th percentile. Notably, significant variations are observed in the ranges among singular values. A high ratio between the largest and smallest singular values can lead to numerical instabilities, indicating an ill-conditioned matrix. To mitigate potential issues, we exclude Hessian matrices with a first singular value $s_1$ exceeding the upper whisker, as well as those with a third singular value $s_3$ below the lower whisker. Overall, these thresholds result in an exclusion of a total of \inum{1893} samples. This leaves a total of \inum{12997} sample pairs $\Theta_{\text{train}} = \{\vect{\theta}_{\text{dyn}}^i, \vect{\theta}_{\text{noise}}^i, \mat{\Gamma}^i\}$, with $\mat{\Gamma}^i_{\text{post}} \approx \mat{\Gamma}^i = (\mat{H}^i)^{-1}$. To remain consistent with other experiments reported in the main manuscript (where we estimate noise and model parameters), we kept the same samples, resulting in different training sizes. For example, we keep \inum{3456} of the matrices of the dataset that compromises \inum{4000} samples. We summarize the remaining training data sizes in \Cref{t:sample_update}.

\begin{table}
\caption{Training samples after removal of Hessian matrices that were negative definite (110 samples) and severely ill-conditioned (\inum{1893} samples) from our original dataset of \inum{15000} samples. We kept the same training sets that we considered in other experiments. Consequently, the size of each training dataset was reduced.\label{t:sample_update}}
\adjtab\centering
\begin{tabular}{rrr}
\toprule
$n_{\text{train}}$ (original) &   $n_{\text{train}}$ (after removal)   & samples dropped  \\
\midrule
\inum{500}   &  \inum{432}~(80\%) & \inum{68}   \\
\inum{1000}  &  \inum{872}~(87\%) & \inum{128}  \\
\inum{2000}  & \inum{1742}~(87\%) & \inum{258}  \\
\inum{4000}  & \inum{3480}~(87\%) & \inum{520}  \\
\inum{8000}  & \inum{6928}~(86\%) & \inum{1072} \\
\bottomrule
\end{tabular}
\end{table}

\begin{figure}
\centering
\includegraphics[width=0.8\textwidth]{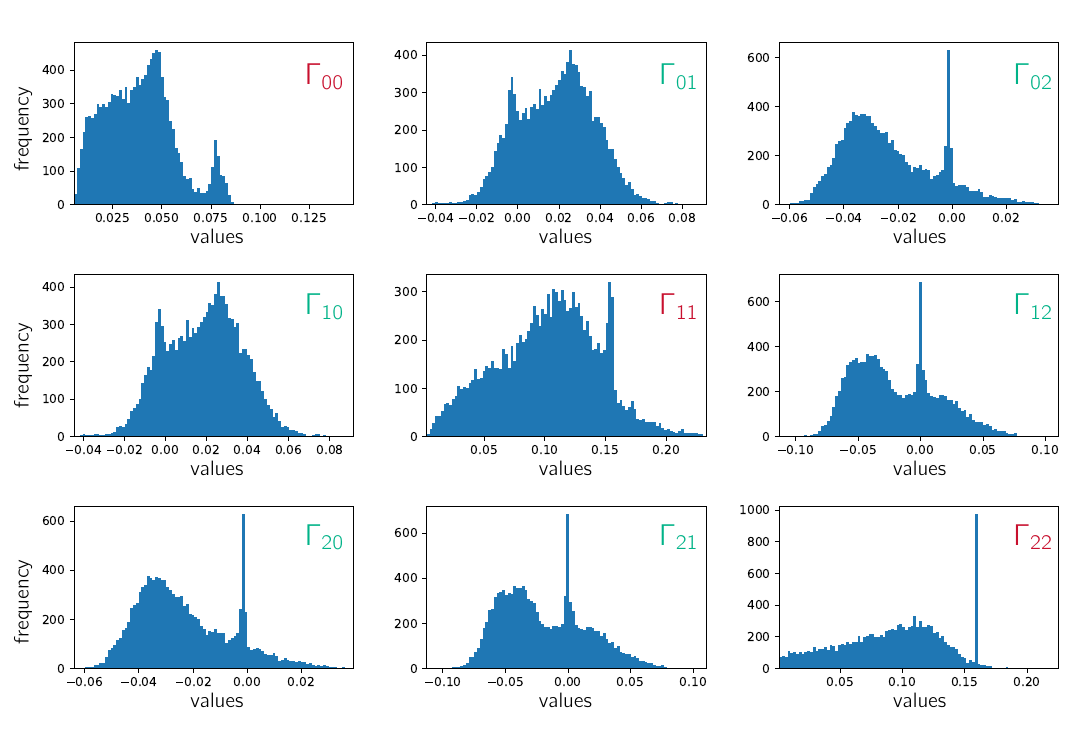}
\caption{Histograms for the approximations $\mat{\Gamma} = (\mat{H})^{-1}$ of the posterior covariance matrices $\mat{\Gamma}_{\text{post}}$ based on the Hessian matrices $\mat{H}$ included in our dataset. We show the frequency for each individual entry $\Gamma_{ij}$, $i,j = 0,1,2$.}
\label{f:refined_covariance3D}
\end{figure}

The entries of the resulting matrices $\mat{\Gamma} \in \spd{3}$ approximating the covariance matrices $\mat{\Gamma}_{\text{post}}$ are shown in \Cref{f:refined_covariance3D}.

To determine if a matrix is symmetric positive definite, we first assess symmetry by comparing the matrix to its transpose. In particular, we quantify the relative error calculated as the norm of the difference divided by the norm of the matrix itself. We use a tolerance of \sci{1e-3}. To assess if a matrix is positive definite we compute an eigenvalue decomposition.

\subsection{Additional Numerical Results}\label{s:supplementary_results}

Below, we provide more insight into the performance of our method.

\subsubsection{Exploration of NN Architectures}\label{s:exploration_nns_addres_supp}

Here, we report additional results for the exploration of the DNN and CNN architectures. In the main manuscript we have reported values for the MdAPE and CDET score. Here, we report values for the SQB and CMSE. The values for the DNN are reported in \Cref{t:dnn_exploration_sqb_cmse}. The values for the CNN are reported in \Cref{t:cnn_exploration_sqb_cmse}. In the main manuscript we selected a CNN architecture with 5 layers and $n_f = 8$ based on these experiments.

\begin{table}
\caption{Exploration of DNN architectures. We report values for the SQB and the CMSE (in brackets). We consider time series data as features (\ts\ setting). We report results for a varying number of nodes/units $n_u$ per layer, where $n_u$ ranges from 4 to 256. We vary the number of layers from 2 to 20. \label{t:dnn_exploration_sqb_cmse}}
\centering\adjtab
\begin{tabular}{rcccccc}
\toprule
$n_u$ & \textbf{4 layers} & \textbf{8 layers}  & \textbf{12 layers} & \textbf{16 layers} & \textbf{20 layers} \\\midrule
  4 & \sci{0.000174} (\decc{0.139788}) & \sci{0.001143} (\decc{0.148505}) & \sci{0.000613} (\decc{0.182531}) & \sci{0.001308} (\decc{0.213447}) & \sci{0.000069} (\decc{0.226219})\\
  8 & \sci{0.001089} (\decc{0.107527}) & \sci{0.002054} (\decc{0.089362}) & \sci{0.000583} (\decc{0.088555}) & \sci{0.002093} (\decc{0.134545}) & \sci{0.000316} (\decc{0.207014})\\
 16 & \sci{0.000101} (\decc{0.059693}) & \sci{0.000052} (\decc{0.056713}) & \sci{0.000124} (\decc{0.062659}) & \sci{0.000551} (\decc{0.055941}) & \sci{0.000029} (\decc{0.115627})\\
 32 & \sci{0.000299} (\decc{0.052456}) & \sci{0.000309} (\decc{0.048951}) & \sci{0.002413} (\decc{0.051344}) & \sci{0.003185} (\decc{0.056672}) & \sci{0.003360} (\decc{0.067949})\\
 64 & \sci{0.000041} (\decc{0.042973}) & \sci{0.001130} (\decc{0.050413}) & \sci{0.000269} (\decc{0.045581}) & \sci{0.000054} (\decc{0.042783}) & \sci{0.001207} (\decc{0.059791})\\
128 & \sci{0.000390} (\decc{0.048594}) & \sci{0.004545} (\decc{0.040882}) & \sci{0.000348} (\decc{0.043334}) & \sci{0.002353} (\decc{0.044391}) & \sci{0.002278} (\decc{0.212177})\\
256 & \sci{0.000524} (\decc{0.043088}) & \sci{0.000160} (\decc{0.043556}) & \sci{0.002280} (\decc{0.113831}) & \sci{0.005122} (\decc{0.108363}) & \sci{0.000960} (\decc{0.184936})\\
\bottomrule
\end{tabular}
\end{table}

\begin{table}
\caption{Exploration of CNN architectures. We report values for the SQB and the CMSE (in brackets). We consider time series data as features (\ts\ setting). We report results for a varying number of filters $n_f$ per layer, where $n_f$ ranges from 2 to 128. We vary the number of layers from 2 to 6.\label{t:cnn_exploration_sqb_cmse}}
\centering\adjtab
\begin{tabular}{rccccc}
\toprule
$n_f$ & \textbf{2 layers} & \textbf{3 layers}  & \textbf{4 layers} & \textbf{5 layers} & \textbf{6 layers}  \\
\midrule
  2 & \sci{0.000073} (\decc{0.040272}) & \sci{0.000259} (\decc{0.030574}) & \sci{0.000335} (\decc{0.027666}) & \sci{0.000070} (\decc{0.029799}) & \sci{0.000310} (\decc{0.029026}) \\
  4 & \sci{0.000070} (\decc{0.042606}) & \sci{0.000347} (\decc{0.030797}) & \sci{0.000206} (\decc{0.026240}) & \sci{0.000612} (\decc{0.025537}) & \sci{0.000384} (\decc{0.028846}) \\
  8 & \sci{0.000601} (\decc{0.042380}) & \sci{0.000400} (\decc{0.035104}) & \sci{0.000118} (\decc{0.027971}) & \sci{0.000449} (\decc{0.025696}) & \sci{0.000563} (\decc{0.027174}) \\
 16 & \sci{0.000305} (\decc{0.045562}) & \sci{0.000305} (\decc{0.038399}) & \sci{0.000104} (\decc{0.034216}) & \sci{0.000212} (\decc{0.030594}) & \sci{0.000215} (\decc{0.032068}) \\
 32 & \sci{0.000224} (\decc{0.045756}) & \sci{0.000166} (\decc{0.040999}) & \sci{0.000121} (\decc{0.032905}) & \sci{0.000661} (\decc{0.033128}) & \sci{0.000369} (\decc{0.032778}) \\
 64 & \sci{0.000469} (\decc{0.046862}) & \sci{0.000892} (\decc{0.041518}) & \sci{0.003220} (\decc{0.035706}) & \sci{0.000484} (\decc{0.032018}) & \sci{0.000602} (\decc{0.037002}) \\
128 & \sci{0.000309} (\decc{0.050395}) & \sci{0.000202} (\decc{0.044156}) & \sci{0.000160} (\decc{0.036821}) & \sci{0.000638} (\decc{0.034460}) & \sci{2.747946} (\decc{0.344337}) \\
\bottomrule
\end{tabular}
\end{table}

\subsubsection{Effects of Noise in Training Data}\label{s:effects_of_noise_training_supp}

In the main manuscript, we consistently included noise in both training and testing data. Here, we provide some additional experiments to identify if this was indeed necessary. While the setup is slightly different (we invert for three and not two model parameters), we note that similar experiments have been reported in~\cite{rudi2022:parameter}. We repeat them to confirm that the results reported in~\cite{rudi2022:parameter} also hold true for our setup.

\ipoint{Purpose:} To show how the presence (or lack thereof) of noise in the training data affects the prediction quality.

\ipoint{Setup:} These experiments are performed on a validation set of size \inum{4000}, comparing the following three different noise configurations:
\begin{itemize}
\item \nfnf: Both training and testing data do not contain noise.
\item \nfn: Only the testing data contain noise.
\item \nn: Both training and testing data contain noise.
\end{itemize}

\noindent Aside from these different setups, we also vary the training size. We set the learning rate to 0.002 and train for \inum{8000} steps (64 epochs). We consider a DNN with 4 layers and $n_u = 64$ and a CNN with 5 layers and $n_f=8$, which consist of a total \inum{140739} and \inum{17223} trainable parameters, respectively. We only estimate the model parameters $\vect{\theta}_{\text{dyn}}$ and explore how the different settings above affect the performance of our methodology. We only include time series data as features with $n_{\text{feat}} = \inum{2000}$. We vary the number of training data.

\ipoint{Results:} We report results in \Cref{t:dnn_supp1} and \Cref{t:cnn_supp1}. Here, we show results for two results for the selected DNN ($n_u = 64$, with 4 layers) and CNN ($n_f=8$, with 5 layers).

\begin{table}
\caption{MdAPE and CDET values (in brackets) of model parameter predictions with noisy observational data. We consider the following setups: no noise in training and testing (label: \nfnf); no noise in training, noise in testing (label: \nfn); noise in training and testing (label: \nn). We consider a DNN with 4 layers with $n_u = 64$.\label{t:dnn_supp1}}
\adjtab\centering
\begin{tabular}{rccc}
\toprule
$n_{\text{train}}$ &   \nfnf   &   \nfn  & \nn\\
\midrule
 \inum{500} & \decc{0.0793} (\decc{0.7801}) & \decc{0.1000} (\decc{0.7426}) & \decc{0.1267} (\decc{0.6950}) \\
\inum{1000} & \decc{0.0465} (\decc{0.8881}) & \decc{0.0797} (\decc{0.8341}) & \decc{0.0958} (\decc{0.7921}) \\
\inum{2000} & \decc{0.0326} (\decc{0.9251}) & \decc{0.0801} (\decc{0.8549}) & \decc{0.0794} (\decc{0.8525}) \\
\inum{4000} & \decc{0.0462} (\decc{0.9470}) & \decc{0.0976} (\decc{0.8380}) & \decc{0.0628} (\decc{0.8962}) \\
\inum{8000} & \decc{0.0239} (\decc{0.9762}) & \decc{0.0927} (\decc{0.8184}) & \decc{0.0597} (\decc{0.9109}) \\
\bottomrule
\end{tabular}
\end{table}

\begin{table}
\caption{MdAPE and CDET values (in brackets) of model parameter predictions with noisy observational data. We consider the following setups: no noise in training and testing (label: \nfnf); no noise in training, noise in testing (label: \nfn); noise in training and testing (label: \nn). We consider a CNN with 5 layers and $n_f = 32$.\label{t:cnn_supp1}}
\adjtab\centering
\begin{tabular}{rccc}
\toprule
$n_{\text{train}}$ & \nfnf   &   \nfn  & \nn\\
\midrule
 \inum{500} & \decc{0.0221} (\decc{0.9570}) & \decc{0.0917} (\decc{0.8331}) & \decc{0.0681} (\decc{0.8880}) \\
\inum{1000} & \decc{0.0219} (\decc{0.9636}) & \decc{0.0859} (\decc{0.8295}) & \decc{0.0605} (\decc{0.8988}) \\
\inum{2000} & \decc{0.0203} (\decc{0.9687}) & \decc{0.0868} (\decc{0.8390}) & \decc{0.0640} (\decc{0.9183}) \\
\inum{4000} & \decc{0.0119} (\decc{0.9907}) & \decc{0.1037} (\decc{0.7796}) & \decc{0.0490} (\decc{0.9412}) \\
\inum{8000} & \decc{0.0078} (\decc{0.9956}) & \decc{0.1260} (\decc{0.6295}) & \decc{0.0465} (\decc{0.9473}) \\
\bottomrule
\end{tabular}
\end{table}

\ipoint{Observations:} The observations are in line with those reported in \cite{rudi2022:parameter}. The best performance is obtained for noise free data (both, in testing and training; setting \nfnf). Moreover, including noise in the training data is critical if the observations are noise (setting \nn). The performance deteriorate significantly if we do not include noise in the training data, trying to infer model parameters from noise observations (setting \nfn).

\subsubsection{Recovering Noise and Model Parameters: Intrinsic Noise}\label{s:noise_model_para_intrinsic_supp}

In this experiment, we estimate the model parameters $\vect{\theta}_{\text{dyn}} \in \ns{R}^3$ along with the noise parameter $\vect{\theta}_{\text{noise}} = \beta \in \ns{R}$ of the intrinsic noise model (see \Cref{s:noisemodel}).

\ipoint{Purpose:} To explore the performance of the proposed framework for estimating noise and model parameters in the presence of intrinsic noise (modeled as the solution of an SDE; see \Cref{s:noisemodel}).

\ipoint{Setup:} We consider a testing data set of size $n_{\text{test}} = \inum{4000}$. No validation data is used. Training and testing data sets both include noise. Like in the former experiment, we consider three different feature configurations, namely, \ts, \fc, and a concatenation of the two, \tsfc. The \ts\ data has length $n_t=\inum{2000}$. We recover $\vect{\theta} = (\theta_0,\theta_1,\theta_2,\beta) \in \mathbb{R}^4$. The learning rate is set to 0.002 and the batch size is 32. We consider different training sizes. We report results after \inum{8000} training steps (64 epochs).

\ipoint{Results:} We report values for the MdAPE and CDET scores for the recovery of the individual parameters in \Cref{t:cnn_noise_recovery_sde}. We show scatter plots for the prediction of the model parameters in \Cref{f:para_recovery_various_noises_scatter} (middle block) for a training size of $n_{\text{train}} = \inum{8000}$. The horizontal and vertical axes correspond to the true and predicted values, respectively. Optimal predictions are indicated by a diagonal red line.

\begin{table}
\centering
\caption{Assessments of joint recovery of ODE model and noise parameters using observational data perturbed by an intrinsic noise model. We report values for $\varepsilon_{\text{MdAPE}}$ and $\varepsilon_{\text{CDET}}$ (in parenthesis). We use a CNN architecture with 5 layers and $n_f = 8$. We report results for a varying number of training samples, ranging from \inum{500} to \inum{8000}. We report results for different choices of the features we consider during training (\ts: time series data; \fc: Fourier coefficients; \tsfc: combination thereof). The remaining parameters are as identified in the text.} \label{t:cnn_noise_recovery_sde}
\centering\adjtab
\begin{tabular}{rlccccc}
\toprule
$n_{\text{train}}$ & \bf feature  & $\theta_0$ & $\theta_1$  & $\theta_2$  &  $\beta$ \\
\midrule
           & \ts    & \decc{0.136811} (\decc{0.902170}) & \decc{0.206097} (\decc{0.849769}) & \decc{0.042116} (\decc{0.642868}) & \decc{0.187234} (\decc{0.416194}) \\
\inum{500} & \fc    & \decc{0.309822} (\decc{0.363905}) & \decc{0.422941} (\decc{0.527187}) & \decc{0.039582} (\decc{0.700946}) & \decc{0.068423} (\decc{0.921352}) \\
           & \tsfc  & \decc{0.124536} (\decc{0.909712}) & \decc{0.224529} (\decc{0.817417}) & \decc{0.034031} (\decc{0.787709}) & \decc{0.069457} (\decc{0.914216}) \\
\hline
            & \ts   & \decc{0.131153} (\decc{0.901045}) & \decc{0.216924} (\decc{0.848321}) & \decc{0.041409} (\decc{0.654827}) & \decc{0.165043} (\decc{0.510944}) \\
\inum{1000} & \fc   & \decc{0.298006} (\decc{0.315114}) & \decc{0.360054} (\decc{0.575342}) & \decc{0.036929} (\decc{0.743434}) & \decc{0.069492} (\decc{0.924951}) \\
            & \tsfc & \decc{0.126625} (\decc{0.914066}) & \decc{0.180343} (\decc{0.896579}) & \decc{0.027967} (\decc{0.850921}) & \decc{0.064924} (\decc{0.921001}) \\
\hline
            & \ts   & \decc{0.134962} (\decc{0.897484}) & \decc{0.210217} (\decc{0.856630}) & \decc{0.041829} (\decc{0.651745}) & \decc{0.160820} (\decc{0.581400}) \\
\inum{2000} & \fc   & \decc{0.251017} (\decc{0.282973}) & \decc{0.201250} (\decc{0.606228}) & \decc{0.035827} (\decc{0.755409}) & \decc{0.064199} (\decc{0.934919}) \\
            & \tsfc & \decc{0.109777} (\decc{0.926248}) & \decc{0.163526} (\decc{0.906402}) & \decc{0.025674} (\decc{0.867268}) & \decc{0.061696} (\decc{0.940567}) \\
\hline
            & \ts   & \decc{0.133277} (\decc{0.897772}) & \decc{0.200016} (\decc{0.863724}) & \decc{0.040145} (\decc{0.687695}) & \decc{0.146561} (\decc{0.638560}) \\
\inum{4000} & \fc   & \decc{0.215323} (\decc{0.320782}) & \decc{0.189849} (\decc{0.620172}) & \decc{0.033627} (\decc{0.784847}) & \decc{0.066331} (\decc{0.935167}) \\
            & \tsfc & \decc{0.109455} (\decc{0.931169}) & \decc{0.180231} (\decc{0.904188}) & \decc{0.025762} (\decc{0.880200}) & \decc{0.055618} (\decc{0.950992}) \\
\hline
            & \ts   & \decc{0.123414} (\decc{0.909785}) & \decc{0.178593} (\decc{0.892297}) & \decc{0.037505} (\decc{0.706447}) & \decc{0.130753} (\decc{0.731572}) \\
\inum{8000} & \fc   & \decc{0.275825} (\decc{0.311528}) & \decc{0.301668} (\decc{0.670045}) & \decc{0.033184} (\decc{0.790156}) & \decc{0.060447} (\decc{0.943866}) \\
            & \tsfc & \decc{0.108636} (\decc{0.930940}) & \decc{0.154504} (\decc{0.919471}) & \decc{0.022932} (\decc{0.895462}) & \decc{0.054758} (\decc{0.955802}) \\
\bottomrule
\end{tabular}%
\end{table}

\ipoint{Observations:} We notice similar patterns as in the previous experiment with additive noise. Time series features alone (with smaller training sizes) yields a poor recovery of the noise model parameter $\beta$. Overall, the performance is slightly better than for additive noise; we do no longer observe negative values for $\varepsilon_{\text{CDET}}$. Using Fourier coefficients alone yields a poor recovery of the model parameters $\vect{\theta}_{\text{dyn}}$ of the dynamical system. Combining both features, time series and Fourier coefficients (\tsfc\ setting) yields the overall best recovery. Compared with the former experiment, we also require less training data.

\subsubsection{Predicted State}\label{s:predicted_state_supp}

Here, we include additional visualizations for the results reported in the main manuscript. In particular, we include plots for the state variable associated with the predicted parameters.

\ipoint{Purpose:} We want to explore how the accuracy in the prediction of the parameters of the dynamical system translates relates to errors in the predicted state; this is critical since during inference the NN predicts quantities of interest without considering (i.e., minimizing) the mismatch between predicted and observed state. That is, our predictor is ``likelihood free.'' This might be a significant issue when it comes to actual applications. If we are, for example, solely interested in the underlying model parameters, errors in the parameter predictions for the dynamical system are our primary concern. However, if we are actually interested in the associated state of our system, errors in the predicted state are critical. Since we do not control these errors during inference, it is important to gain an understanding of how well the predicted states match the true state and observed data.

\ipoint{Setup:} We jointly estimate model parameters $\vect{\theta}_{\text{model}}$, noise parameters ($\vect{\theta}_{\text{noise}} = (\rho, \sigma)$ for additive noise, $\vect{\theta}_{\text{noise}} = \beta$ for intrinsic noise, and $\vect{\theta}_{\text{noise}} = (\rho,\sigma,\beta)$ for the combination of these two noise models), as well as entries of the covariance matrix $\mat{\Gamma}_{\text{post}}$. The results reported below are obtained for a CNN with 5 layers and $n_f = 8$ and a training size of $n_{\text{train}} = \inum{6928}$. We consider additive and intrinsic noise, as well as a combination thereof.

\begin{figure}
\centering
\includegraphics[width=0.7\textwidth]{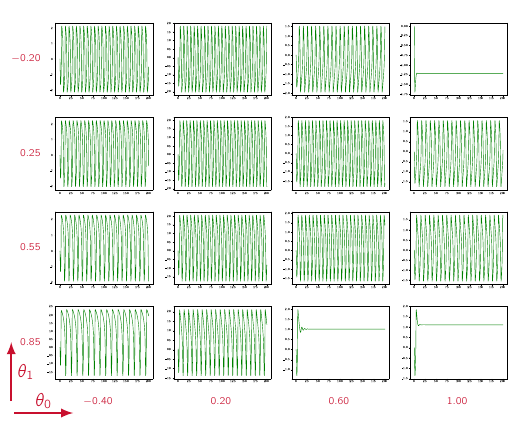}
\caption{Forward simulation results for the membrane potential $u$ for varying parameters $\vect{\theta}_{\text{dyn}}$. We set $\theta_2$ to 0.2 and vary $\theta_0$ and $\theta_1$ as indicated. The simulation time is 200\,ms.\label{f:fwdsim-range}}
\end{figure}

\begin{figure}
\centering
\includegraphics[trim={1.5cm 3cm 1.5cm 3cm}, clip,width=0.9\textwidth]{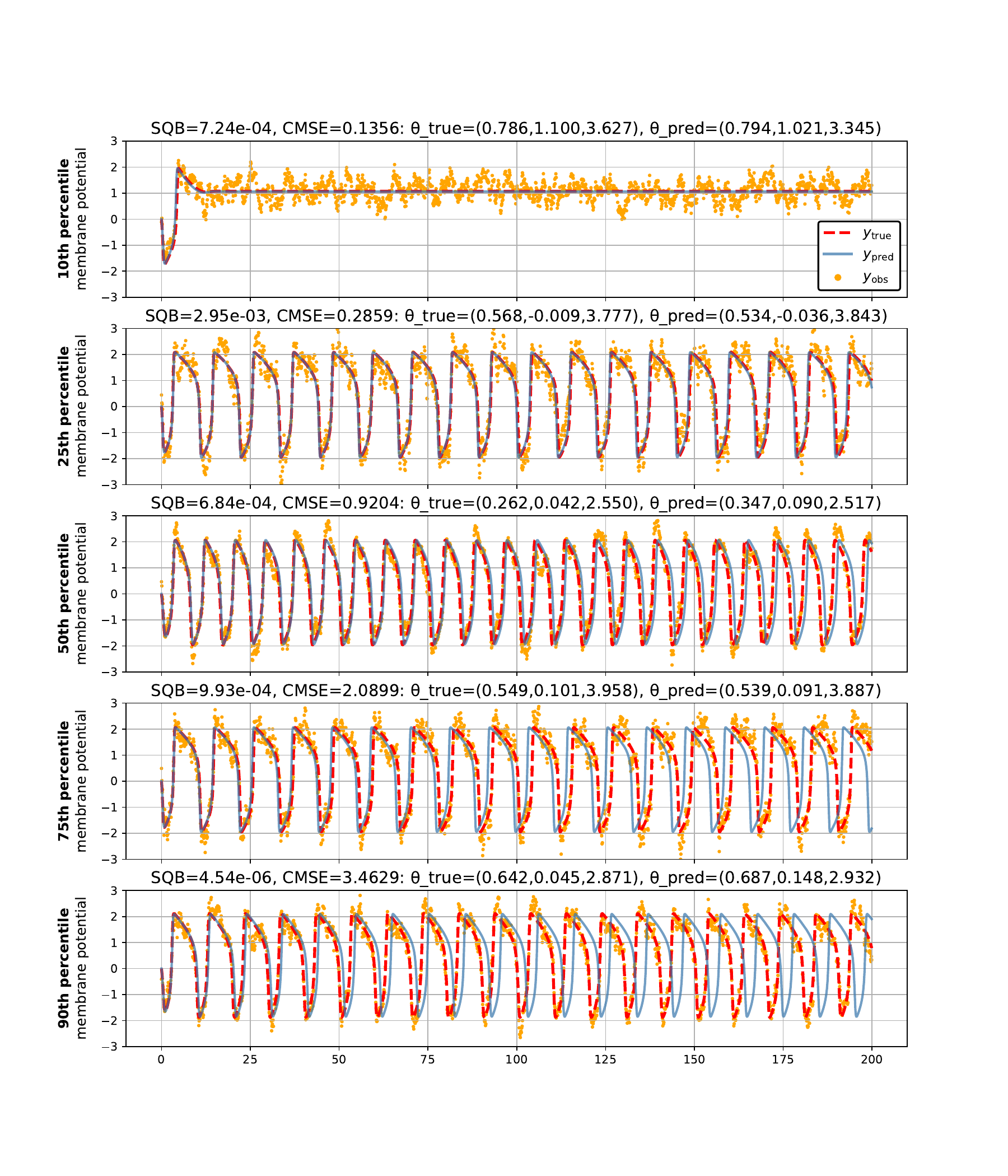}
\caption{Simulation results for the state variable $u$ for representative estimates $\vect{\theta}_{\text{pred}}$ of the model parameters $\vect{\theta}_{\text{dyn}}$. We consider an additive noise model. We show from top to bottom increasing quantiles of the MSE between the observed and estimated time series. The blue line corresponds to the membrane potential $u$ generated using predicted parameters (i.e., solving the forward model). The orange dots represent the observations (perturbed membrane potential using an additive noise model). The red dashed line corresponds to the ground truth membrane potential. We jointly estimate model parameters $\vect{\theta}_{\text{dyn}} = (\theta_0,\theta_1,\theta_2)$, noise parameters $\vect{\theta}_{\text{noise}} = (\rho, \sigma)$ (additive noise model), as well as entries of the covariance matrix $\mat{\Gamma}_{\text{post}}$. The number of training samples is $n_{\text{train}} = \inum{6928}$. The number of testing samples is $n_{\text{test}} = \inum{3456}$. On top of each graph we report the squared bias, the C-MSE, as well as the predicted ($\vect{\theta}_{\text{pred}}$) and true values ($\vect{\theta}_{\text{true}}$) for the model parameters $\vect{\theta}_{\text{dyn}}$.\label{f:cnn-add-tscomp}}
\end{figure}

\begin{figure}
\centering
\includegraphics[trim={1.5cm 3cm 1.5cm 3cm}, clip,width=0.9\textwidth]{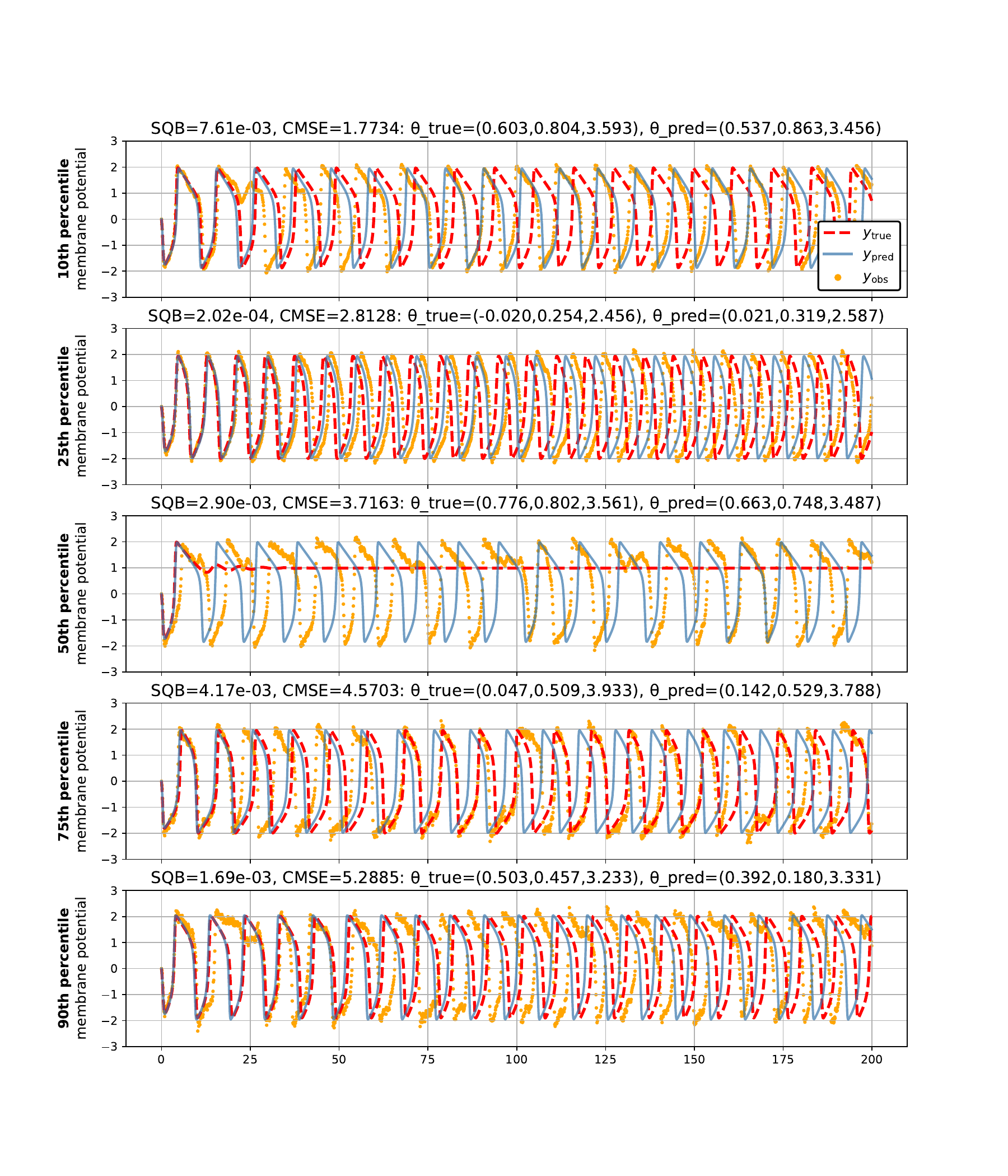}
\caption{Simulation results for the state variable $u$ for representative estimates $\vect{\theta}_{\text{pred}}$ of the model parameters $\vect{\theta}_{\text{dyn}}$. We consider an intrinsic noise model. We show from top to bottom increasing quantiles of the MSE between the observed and estimated time series. The blue line corresponds to the membrane potential $u$ generated using predicted parameters (i.e., solving the forward model). The orange dots represent the observations (perturbed membrane potential using an additive noise model). The red dashed line corresponds to the ground truth membrane potential. We jointly estimate model parameters $\vect{\theta}_{\text{dyn}} = (\theta_0,\theta_1,\theta_2)$, noise parameters $\vect{\theta}_{\text{noise}} = \beta$ (intrinsic noise model), as well as entries of the covariance matrix $\mat{\Gamma}_{\text{post}}$. The number of training samples is $n_{\text{train}} = \inum{6928}$. The number of testing samples is $n_{\text{test}} = \inum{3456}$. On top of each graph we report the squared bias, the C-MSE, as well as the predicted ($\vect{\theta}_{\text{pred}}$) and true values ($\vect{\theta}_{\text{true}}$) for the model parameters $\vect{\theta}_{\text{dyn}}$.\label{f:cnn-int-tscomp}}
\end{figure}

\begin{figure}
\centering
\includegraphics[trim={1.5cm 3cm 1.5cm 3cm}, clip,width=0.9\textwidth]{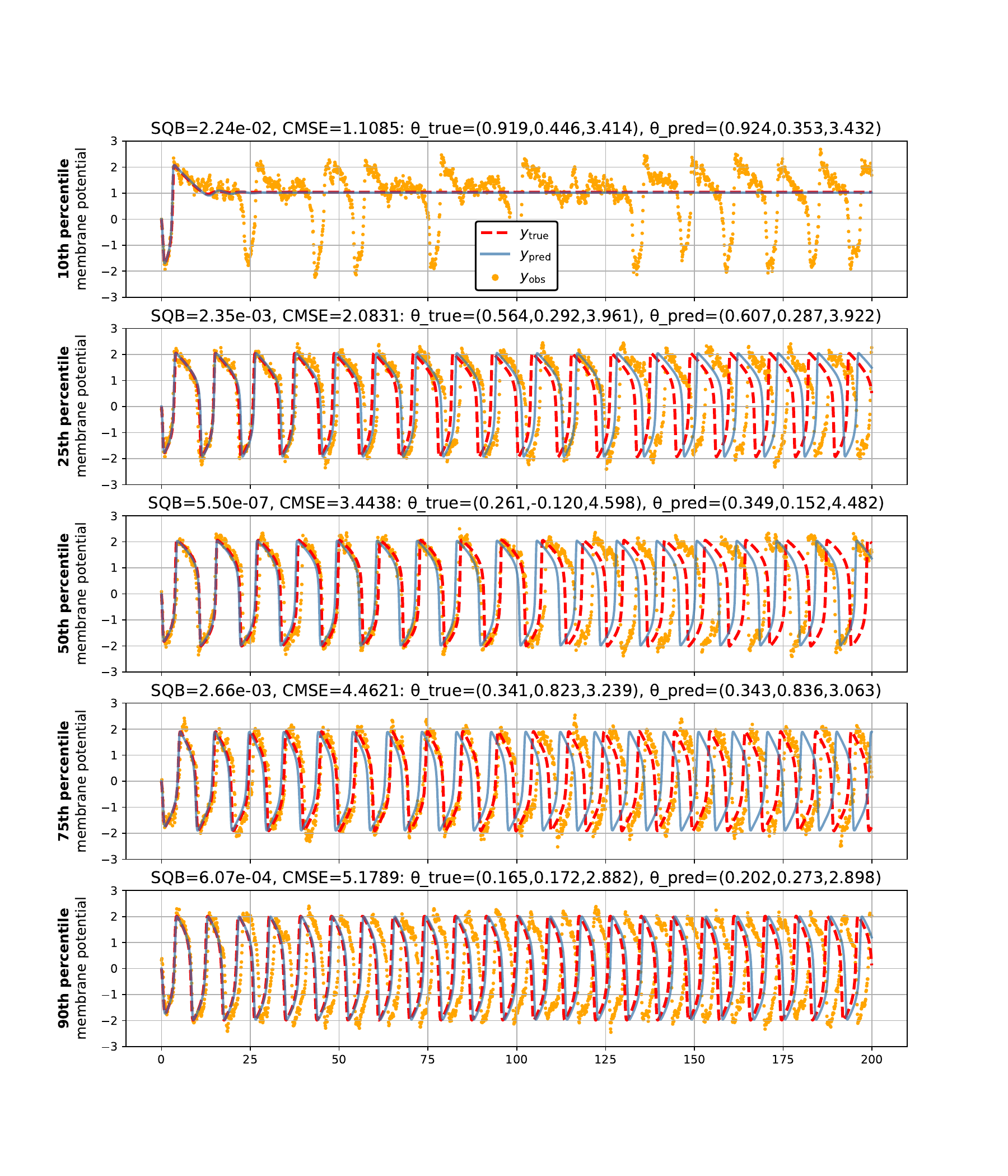}
\caption{Simulation results for the state variable $u$ for representative estimates $\vect{\theta}_{\text{pred}}$ of the model parameters $\vect{\theta}_{\text{dyn}}$. We perturb the simulations with an additive and an intrinsic noise model to obtain the observed membrane potential. We show from top to bottom increasing quantiles of the MSE between the observed and estimated time series. The blue line corresponds to the membrane potential $u$ generated using predicted parameters (i.e., solving the forward model). The orange dots represent the observations (perturbed membrane potential using an additive noise model). The red dashed line corresponds to the ground truth membrane potential. We jointly estimate model parameters $\vect{\theta}_{\text{dyn}} = (\theta_0,\theta_1,\theta_2)$, noise parameters $\vect{\theta}_{\text{noise}} = (\rho,\sigma,\beta)$ (combined noise model), as well as entries of the covariance matrix $\mat{\Gamma}_{\text{post}}$. The number of training samples is $n_{\text{train}} = \inum{6928}$. The number of testing samples is $n_{\text{test}} = \inum{3456}$. On top of each graph we report the squared bias, the C-MSE, as well as the predicted ($\vect{\theta}_{\text{pred}}$) and true values ($\vect{\theta}_{\text{true}}$) for the model parameters $\vect{\theta}_{\text{dyn}}$.\label{f:cnn-combhalf-tscomp}}
\end{figure}

\ipoint{Results:} We show representative forward simulation results for the membrane potential $u$ in \Cref{f:fwdsim-range}. We show a comparison between the predicted membrane potential and the observed membrane potential for additive noise in \Cref{f:cnn-add-tscomp}, for intrinsic noise in \Cref{f:cnn-int-tscomp}, and for a combination of these two noise models in \Cref{f:cnn-combhalf-tscomp}. For each of these plots, we show from top to bottom increasing quantiles of the MSE between true and estimated time series.

\ipoint{Observations:}  The predictions for the cell membrane potential are in good agreement with the true state in most cases, even if the perturbed (noisy) dynamics are quite different. This is particularly true for the additive noise model. For the intrinsic noise we can see larger deviations between the predicted, observed, and true state; small noise perturbation manifest in larger differences in the dynamics.

\subsubsection{Non-Positive Definiteness in Predicted Covariance Data}\label{s:pred_npsd_mat_fix_supp}

If we use off-the-shelf NNs to predict the covariance matrices we cannot expect the predictions to be positive definite. We take advantage of the (in theory) symmetry of the covariance matrices and reduce complexity of the learning task by recovering only the upper triangular entries. Our testing set for the experiments with covariance data consists of of \inum{3456} samples. We remove the data that correspond to parameters that yield close to singular (severely ill-conditioned) and negative definite Hessian matrices from our training and testing datasets. We took this approach to remain consistent with our datasets used in prior experiments that did not involve the estimation of covariance matrices. For the NN predictions we took a different approach.

Exemplary results for NN predictions of the covariance matrix are shown in \Cref{t:posdef_pred_check}. As we can see from these experiments, roughly 30\% of the predictions are negative definite or indefinite. To address this issue one can, e.g., leverage the methodology described in~\cite{Fasiha:2017}; we estimate the nearest positive definite matrix for the predicted covariances. The nearest positive definite matrix is computed by first symmetrizing the predicted matrix, say $\mat{\tilde{\Gamma}}_{\text{pred}}$, followed by an averaging with a symmetric matrix $\mat{\Gamma}_{\text{sym}}$ generated from the singular value decomposition of the matrix $\mat{\tilde{\Gamma}}_{\text{pred}}$. This average is symmetrized and adjusted (iteratively) through the addition of a positive multiple of the identity matrix. Specifically, the identity is scaled by the smallest eigenvalue of the previous iteration until a positive definite matrix is achieved. For details we refer to~\cite{Higham:1988c}. In the main manuscript we have developed a different strategy that maps the data to the tangent space.

\begin{table}
\adjtab\centering
\caption{Predicted covariance matrices. We report the number of covariance matrices that are positive definite and negative definite/indefinite after feeding them through the trained NN, respectively. The testing data consists of \inum{3456} samples.\label{t:posdef_pred_check}}
\adjtab\centering
\begin{tabular}{lrr}
\toprule
\bf noise type      & positive definite &  negative definite/indefinite  \\
\midrule
additive   & 2399~($\sim70\%)$ & 1057~($\sim30\%)$ \\
intrinsic  & 2573~($\sim74\%)$ &  883~($\sim26\%)$ \\
combined   & 2195~($\sim64\%)$ & 1261~($\sim36\%)$ \\
\bottomrule
\end{tabular}
\end{table}

We report results after mapping the predictions to the closest positive definite matrix as just described in \Cref{t:cnn_posdef_cov_recovery}.

\begin{sidewaystable}
\caption{MdAPE and CDET (in brackets) for the estimation of entries of the covariance matrix of the posterior distributions. Here, we map the predicted covariance matrices (from results in \Cref{s:covariance_results}) that are negative definite (indefinite) to the nearest positive definite matrix. We apply the same methodology as described in~\cite{Fasiha:2017}.\label{t:cnn_posdef_cov_recovery}}
\adjtab\centering
\begin{tabular}{rlcccccc}
\toprule
$n_{\text{train}}$ & \bf noise type  & $\Gamma_{00}$ & $\Gamma_{01}$  & $\Gamma_{02}$  &  $\Gamma_{11}$ & $\Gamma_{12}$  &  $\Gamma_{22}$ \\
\midrule
           &  additive & \decc{0.133504292020} (\decc{0.741520786356}) & \decc{0.262364483185} (\decc{0.770455188265}) & \decc{0.181025603309} (\decc{0.849902338237}) & \decc{0.160381276389} (\decc{0.615173469614}) & \decc{0.300546112363} (\decc{0.788002319382}) & \decc{0.131414666763} (\decc{0.816647542402}) \\
\inum{872} & intrinsic & \decc{0.164802266728} (\decc{0.582542285922}) & \decc{0.299485535660} (\decc{0.735995635351}) & \decc{0.203969688873} (\decc{0.767942507163}) & \decc{0.171902923201} (\decc{0.555150675173}) & \decc{0.299764090466} (\decc{0.807148643342}) & \decc{0.144893326077} (\decc{0.750987235337}) \\
           & combined  & \decc{0.141328194906} (\decc{0.671935554419}) & \decc{0.297715641386} (\decc{0.753509631153}) & \decc{0.201709669561} (\decc{0.806772795723}) & \decc{0.174235826219} (\decc{0.581655886039}) & \decc{0.324105107559} (\decc{0.719749800757}) & \decc{0.139306572904} (\decc{0.784312394732}) \\
\hline
       & additive  & \decc{0.128617432517} (\decc{0.771607263388}) & \decc{0.253825456043} (\decc{0.789514489538}) & \decc{0.173329831233} (\decc{0.864081089908}) & \decc{0.155609924137} (\decc{0.651575695515}) & \decc{0.263889439635} (\decc{0.826947759752}) & \decc{0.124539008215} (\decc{0.826817337260}) \\
1742 & intrinsic  & \decc{0.157778424982} (\decc{0.594796776952}) & \decc{0.282876437910} (\decc{0.759117235020}) & \decc{0.190108733795} (\decc{0.800077791426}) & \decc{0.169705500258} (\decc{0.586362944980}) & \decc{0.278555406667} (\decc{0.823720739733}) & \decc{0.146655817379} (\decc{0.749390574843}) \\
       & combined  &    \decc{0.132036015501} (\decc{0.716395467381}) & \decc{0.244236151791} (\decc{0.816805049742}) & \decc{0.169092861915} (\decc{0.858060043673}) & \decc{0.159277841331} (\decc{0.648327918725}) & \decc{0.252851341479} (\decc{0.845084668793}) & \decc{0.128134890439} (\decc{0.822328347436}) \\
\hline
       & additive  & \decc{0.126007833453} (\decc{0.781976274567}) & \decc{0.219276938682} (\decc{0.828266483709}) & \decc{0.169684063887} (\decc{0.861147534916}) & \decc{0.152471757647} (\decc{0.658928438701}) & \decc{0.236518996708} (\decc{0.859458169011}) & \decc{0.117226031030} (\decc{0.840603104082}) \\
3480 &  intrinsic  & \decc{0.150350574700} (\decc{0.614530834044}) & \decc{0.280020335584} (\decc{0.768413143361}) & \decc{0.188817749696} (\decc{0.809048000415}) & \decc{0.156679555397} (\decc{0.604938534659}) & \decc{0.281474847434} (\decc{0.830842127892}) & \decc{0.140723949099} (\decc{0.759130184551})  \\
       &  combined &   \decc{0.125332879874} (\decc{0.708424913141}) & \decc{0.244409447083} (\decc{0.818491489614}) & \decc{0.164572246713} (\decc{0.868018637806}) & \decc{0.149758048797} (\decc{0.651301129150}) & \decc{0.246830287714} (\decc{0.862684194177}) & \decc{0.124442115144} (\decc{0.822849091839}) \\
\hline
       & additive  & \decc{0.116772487973} (\decc{0.797955628218}) & \decc{0.218652215113} (\decc{0.835654061785}) & \decc{0.153688858560} (\decc{0.887929494642}) & \decc{0.141174618222} (\decc{0.672842368635}) & \decc{0.217864640830} (\decc{0.882264549380}) & \decc{0.110675095406} (\decc{0.867844260206}) \\
6928 &  intrinsic & \decc{0.156017794503} (\decc{0.616839352107}) & \decc{0.256923869966} (\decc{0.792879005159}) & \decc{0.182666376362} (\decc{0.807872877033}) & \decc{0.158658007181} (\decc{0.603803724711}) & \decc{0.261039618091} (\decc{0.858943614971}) & \decc{0.133686387823} (\decc{0.773778043631}) \\
       & combined  &  \decc{0.128995287644} (\decc{0.720707286389}) & \decc{0.241161213322} (\decc{0.814559740419}) & \decc{0.173534202020} (\decc{0.857920220607}) & \decc{0.154415431250} (\decc{0.627272715910}) & \decc{0.252977129298} (\decc{0.850567085851}) & \decc{0.128205641101} (\decc{0.817765351452}) \\
\bottomrule
\end{tabular}
\end{sidewaystable}

\subsubsection{Sensitivity to Initialization}\label{s:nn_initialization_supp}

Here, we report more detailed results for the sensitivity of the considered NN architectures to the initialization of the NN weights. We consider 10 different random initialization.

\ipoint{Purpose:} Assess the sensitivity of the NN predictions with respect to the initialization of the NN weights.

\ipoint{Setup:} These results complement the summary results reported in \Cref{s:nn_initialization}. We refer to the main manuscript for details about the setup.

\ipoint{Results:} We report results for CNN architecture in \Cref{t:cnn_seed_inits_supp} for the CNN and \Cref{t:dnn_seed_inits_supp} for the DNN, respectively. The training time for the CNN architecture is roughly 73.43 seconds total; the training time for the DNN averages 48.44 seconds total.

\ipoint{Observations:} Since we report summary results in the main manuscript, the findings in this supplement are in line with the observations we made in the main manuscript: The NNs considered in this work are not sensitive to the initialization in the context we employ them.

\begin{table}
\caption{Performance for the recovery of the model parameters $\vect{\theta}_{\text{dyn}}$ (left block) and covariance entries (right block) using a CNN across different random seeds. We report values for CDET, MdAPE, CMSE, and SQB. The training time is reported in seconds.} \label{t:cnn_seed_inits_supp}
\centering\adjtab
\begin{tabular}{rccccccccc}\toprule
&
\multicolumn{4}{c}{\textbf{parameter recovery}}
&
\multicolumn{4}{c}{\textbf{covariance entries}}
\\
\midrule
\textbf{id}
& \textbf{CDET} & \textbf{MdAPE} & \textbf{CMSE} & \textbf{SQB}
& \textbf{CDET} & \textbf{MdAPE} & \textbf{CMSE} & \textbf{SQB} & \textbf{time}
\\
\midrule
1  & \decc{0.85811373} & \decc{0.07070069} & \sci{0.02954620} & \sci{0.00038791} & \decc{0.79814371} & \decc{0.16971319} & \sci{0.00121439} & \sci{0.00000430} & \decc{72.532443}\\
2  & \decc{0.84982443} & \decc{0.07263948} & \sci{0.03323397} & \sci{0.00018746} & \decc{0.79090221} & \decc{0.17597251} & \sci{0.00126281} & \sci{0.00000724} & \decc{73.382435}\\
3  & \decc{0.85751889} & \decc{0.07147789} & \sci{0.03016186} & \sci{0.00043289} & \decc{0.80057997} & \decc{0.17277908} & \sci{0.00122130} & \sci{0.00001317} & \decc{73.974760}\\
4  & \decc{0.86133814} & \decc{0.06943808} & \sci{0.02939858} & \sci{0.00036715} & \decc{0.80480512} & \decc{0.17080515} & \sci{0.00118212} & \sci{0.00000644} & \decc{73.574846}\\
5  & \decc{0.84940674} & \decc{0.07389642} & \sci{0.03073328} & \sci{0.00006405} & \decc{0.78042756} & \decc{0.17778180} & \sci{0.00135231} & \sci{0.00000840} & \decc{73.610566}\\
6  & \decc{0.86491833} & \decc{0.06884264} & \sci{0.02817768} & \sci{0.00109931} & \decc{0.79651536} & \decc{0.17259996} & \sci{0.00123840} & \sci{0.00000550} & \decc{74.776071}\\
7  & \decc{0.84479024} & \decc{0.07522555} & \sci{0.03321152} & \sci{0.00033187} & \decc{0.77714655} & \decc{0.17881705} & \sci{0.00136428} & \sci{0.00000609} & \decc{72.969865}\\
8  & \decc{0.86130151} & \decc{0.07095884} & \sci{0.02822826} & \sci{0.00094709} & \decc{0.79340480} & \decc{0.17427275} & \sci{0.00126414} & \sci{0.00000870} & \decc{73.525416}\\
9  & \decc{0.85635680} & \decc{0.07117056} & \sci{0.03063645} & \sci{0.00018628} & \decc{0.79156576} & \decc{0.17140193} & \sci{0.00129500} & \sci{0.00000226} & \decc{72.807923}\\
10 & \decc{0.85291453} & \decc{0.07316158} & \sci{0.03156879} & \sci{0.00017458} & \decc{0.77797059} & \decc{0.17690013} & \sci{0.00135505} & \sci{0.00000484} & \decc{73.168425}\\
\bottomrule
\end{tabular}
\end{table}

\begin{table}
\caption{Performance for the recovery of the model parameters $\vect{\theta}_{\text{dyn}}$ (left block) and covariance entries (right block) using a DNN across different random seeds. We report values for CDET, MdAPE, CMSE, and SQB. The training time is reported in seconds.} \label{t:dnn_seed_inits_supp}
\centering\adjtab
\begin{tabular}{rccccccccc}\toprule
&
\multicolumn{4}{c}{\textbf{parameter recovery}}
&
\multicolumn{4}{c}{\textbf{covariance entries}}
\\
\midrule
\textbf{id}
& \textbf{CDET} & \textbf{MdAPE} & \textbf{CMSE} & \textbf{SQB}
& \textbf{CDET} & \textbf{MdAPE} & \textbf{CMSE} & \textbf{SQB} & \textbf{time}
\\
\midrule
1  & \decc{0.77429643} & \decc{0.10465082} & \sci{0.07316643} & \sci{0.00064851} & \decs{0.62366736} & \decc{0.23029476} & \sci{0.00219708} & \sci{0.00002238} & \decc{47.30905900} \\
2  & \decc{0.79425307} & \decc{0.10207084} & \sci{0.06750974} & \sci{0.00022292} & \decs{0.65380618} & \decc{0.22594168} & \sci{0.00204072} & \sci{0.00001420} & \decc{50.21563900} \\
3  & \decc{0.77101114} & \decc{0.10355067} & \sci{0.07173137} & \sci{0.00115763} & \decs{0.64612204} & \decc{0.22654514} & \sci{0.00206216} & \sci{0.00003052} & \decc{48.24864000} \\
4  & \decc{0.78682716} & \decc{0.10347767} & \sci{0.06923674} & \sci{0.00043100} & \decs{0.64930001} & \decc{0.22289624} & \sci{0.00204874} & \sci{0.00002919} & \decc{49.18488400} \\
5  & \decc{0.77979862} & \decc{0.10437568} & \sci{0.06959951} & \sci{0.00494260} & \decs{0.63981368} & \decc{0.23169672} & \sci{0.00212889} & \sci{0.00004747} & \decc{49.74031600} \\
6  & \decc{0.78473497} & \decc{0.10618298} & \sci{0.06872220} & \sci{0.00129582} & \decs{0.64031446} & \decc{0.23171644} & \sci{0.00210541} & \sci{0.00002960} & \decc{49.20512200} \\
7  & \decc{0.78709973} & \decc{0.10211766} & \sci{0.06769572} & \sci{0.00002371} & \decs{0.64773904} & \decc{0.22657106} & \sci{0.00207223} & \sci{0.00002224} & \decc{46.98229800} \\
8  & \decc{0.77816604} & \decc{0.10202335} & \sci{0.06761681} & \sci{0.00185316} & \decs{0.65157269} & \decc{0.22400179} & \sci{0.00201093} & \sci{0.00001798} & \decc{47.42461700} \\
9  & \decc{0.78261366} & \decc{0.10685442} & \sci{0.06935019} & \sci{0.00015708} & \decs{0.64394879} & \decc{0.23105876} & \sci{0.00212079} & \sci{0.00001519} & \decc{47.39379200} \\
10 & \decc{0.78798554} & \decc{0.10202425} & \sci{0.06645635} & \sci{0.00138135} & \decs{0.64719754} & \decc{0.22689586} & \sci{0.00206336} & \sci{0.00002463} & \decc{48.73358600} \\
\bottomrule
\end{tabular}
\end{table}

\subsubsection{Sensitivity to Initialization --- Cross Validation}\label{s:nn_initialization_supp_cc}

We additionally evaluate the quality of the selected NN models and the impact of different weight initializations through $k$-fold cross-validation with $k = 6$.

\ipoint{Purpose:} To assess generalizability of the considered NN models independent of the dataset.

\ipoint{Setup:} The dataset is partitioned into six equally sized subsets; in each fold, five subsets are used for training and one for testing. Evaluation metrics are recorded on the held-out test subset for each fold. These experiments follow the same configuration described in \Cref{s:covariance_results}, using \tsfc\ as NN inputs, with observations generated from a combined noise model and using the full set of available samples. The setup and hyperparameters are held constant across runs, except for the random seed controlling initialization. We train for 100 epochs per fold to standardize across experiments.

\ipoint{Results:} The aggregated results are presented in  \Cref{t:cross-validation-para} and \Cref{t:cross-validation-covariance}, respectively.

\ipoint{Observations:} Similar to the results reported in the prior section, we do not observe a significant deterioration with respect to the initialization. Our cross validation results indicate that we have good generalization performance.

\begin{table}
\caption{6-fold cross-validation with (noisy) training and testing sets of sizes $n_{\text{train}} = 8653$ and $n_{\text{test}} = 1731$, respectively; random seeds for initializing NN weights vary per row. We consider a CNN with 5 layers, with a filter setting of $n_f = 8$. This table summarizes inference for all parameters, $\vect{\theta}_{\text{dyn}} = (\theta_0,\theta_1,\theta_2) \in \mathbb{R}^3$ and $\vect{\theta}_{\text{noise}} = (\beta, \rho,\sigma) \in \mathbb{R}^3$.}
\centering\adjtab
\label{t:cross-validation-para}
\begin{tabular}{rcccccccc}
\toprule
 & \multicolumn{2}{c}{\textbf{CDET}} & \multicolumn{2}{c}{\textbf{MdAPE}} & \multicolumn{2}{c}{\textbf{SQB}} & \multicolumn{2}{c}{\textbf{CMSE}} \\
\cmidrule(lr){2-3} \cmidrule(lr){4-5} \cmidrule(lr){6-7} \cmidrule(lr){8-9}
\textbf{id} & \textbf{mean} & \textbf{std} & \textbf{mean} & \textbf{std} & \textbf{mean} & \textbf{std} & \textbf{mean} & \textbf{std} \\
\midrule
1  & \decc{0.839190} & \decc{0.005064} & \decc{0.077335} & \decc{0.001678} & \sci{0.001050} & \sci{0.001038} & \sci{0.033857} & \sci{0.001068} \\
2  & \decc{0.834112} & \decc{0.006916} & \decc{0.077845} & \decc{0.003586} & \sci{0.000537} & \sci{0.000428} & \sci{0.034706} & \sci{0.001764} \\
3  & \decc{0.840920} & \decc{0.008994} & \decc{0.077972} & \decc{0.003224} & \sci{0.000937} & \sci{0.001087} & \sci{0.034750} & \sci{0.001088} \\
4  & \decc{0.837276} & \decc{0.004638} & \decc{0.079842} & \decc{0.002204} & \sci{0.000609} & \sci{0.000864} & \sci{0.034134} & \sci{0.001109} \\
5  & \decc{0.837643} & \decc{0.006804} & \decc{0.078250} & \decc{0.001553} & \sci{0.000207} & \sci{0.000137} & \sci{0.033598} & \sci{0.002170} \\
6  & \decc{0.837265} & \decc{0.007344} & \decc{0.077701} & \decc{0.001892} & \sci{0.000334} & \sci{0.000148} & \sci{0.033988} & \sci{0.001660} \\
7  & \decc{0.833065} & \decc{0.006504} & \decc{0.078845} & \decc{0.001394} & \sci{0.000906} & \sci{0.000777} & \sci{0.035515} & \sci{0.001172} \\
8  & \decc{0.836654} & \decc{0.004688} & \decc{0.078370} & \decc{0.001530} & \sci{0.000523} & \sci{0.000614} & \sci{0.034131} & \sci{0.001389} \\
9  & \decc{0.839720} & \decc{0.005849} & \decc{0.077301} & \decc{0.002416} & \sci{0.000510} & \sci{0.000426} & \sci{0.033978} & \sci{0.001754} \\
10 & \decc{0.838279} & \decc{0.005268} & \decc{0.077799} & \decc{0.001137} & \sci{0.000449} & \sci{0.000764} & \sci{0.034235} & \sci{0.001439} \\
\bottomrule
\end{tabular}
\end{table}

\begin{table}
\caption{6-fold cross-validation with (noisy) training and testing sets of sizes $n_{\text{train}} = 8653$ and $n_{\text{test}} = 1731$, respectively. The random seeds used for initializing NN weights vary per row. We consider a CNN with 5 layers and filter setting of $n_f = 8$. This table summarizes inference for all entries of the covariance matrix.}
\label{t:cross-validation-covariance}
\centering\adjtab
\begin{tabular}{rcccccccc}
\toprule
 & \multicolumn{2}{c}{\textbf{CDET}} & \multicolumn{2}{c}{\textbf{MdAPE}} & \multicolumn{2}{c}{\textbf{SQB}} & \multicolumn{2}{c}{\textbf{CMSE}} \\
\cmidrule(lr){2-3} \cmidrule(lr){4-5} \cmidrule(lr){6-7} \cmidrule(lr){8-9}
\textbf{id} & \textbf{mean} & \textbf{std} & \textbf{mean} & \textbf{std} & \textbf{mean} & \textbf{std} & \textbf{mean} & \textbf{std} \\
\midrule
1  & \decc{0.745472} & \decc{0.002547} & \decc{0.193965} & \decc{0.001977} & \sci{0.000011} & \sci{0.000007} & \sci{0.001620} & \sci{0.000034} \\
2  & \decc{0.745470} & \decc{0.006676} & \decc{0.190170} & \decc{0.003990} & \sci{0.000007} & \sci{0.000004} & \sci{0.001620} & \sci{0.000048} \\
3  & \decc{0.743441} & \decc{0.014386} & \decc{0.193887} & \decc{0.005397} & \sci{0.000011} & \sci{0.000011} & \sci{0.001630} & \sci{0.000061} \\
4  & \decc{0.733735} & \decc{0.005012} & \decc{0.197524} & \decc{0.001453} & \sci{0.000014} & \sci{0.000013} & \sci{0.001686} & \sci{0.000029} \\
5  & \decc{0.744516} & \decc{0.008346} & \decc{0.194047} & \decc{0.001899} & \sci{0.000006} & \sci{0.000003} & \sci{0.001638} & \sci{0.000056} \\
6  & \decc{0.743291} & \decc{0.006043} & \decc{0.193121} & \decc{0.005572} & \sci{0.000006} & \sci{0.000005} & \sci{0.001629} & \sci{0.000064} \\
7  & \decc{0.733419} & \decc{0.012136} & \decc{0.198348} & \decc{0.005115} & \sci{0.000015} & \sci{0.000012} & \sci{0.001696} & \sci{0.000090} \\
8  & \decc{0.732213} & \decc{0.009682} & \decc{0.198078} & \decc{0.004337} & \sci{0.000011} & \sci{0.000011} & \sci{0.001684} & \sci{0.000064} \\
9  & \decc{0.743753} & \decc{0.007760} & \decc{0.195118} & \decc{0.003848} & \sci{0.000008} & \sci{0.000008} & \sci{0.001628} & \sci{0.000047} \\
10 & \decc{0.741600} & \decc{0.002316} & \decc{0.195018} & \decc{0.003616} & \sci{0.000005} & \sci{0.000004} & \sci{0.001643} & \sci{0.000014} \\
\bottomrule
\end{tabular}
\end{table}

\subsubsection{Performance as a Function of the Noise Perturbation}\label{s:noise_sensitivity_supp}

We study the prediction performance as a function of increasing noise perturbations. The results reported in this section complement those reported in \Cref{s:noise_sensitivity}.

\ipoint{Purpose:} We explore the performance of the proposed framework as a function of increasing noise perturbations. We expect the prediction accuracy to deteriorate as the noise level increases. We note that these results are backed into those reported in prior experiments. We attempt to disentangle the errors and report prediction performance as a function of increasing noise.

\begin{figure}
\centering
\includegraphics[width=0.9\textwidth]{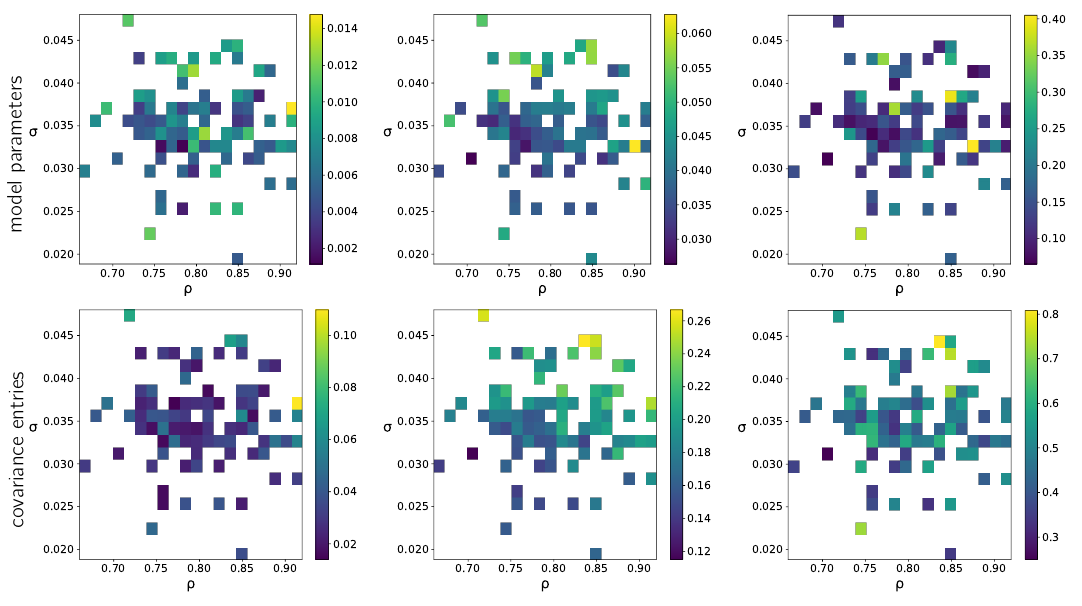}\quad
\caption{2D maps of the relative $\ell^2$-error between NN predictions and true values. We bin the errors as a function of $\rho$ and $\sigma$, respectively. The amplitude of the error is shown in color. The top row corresponds to the errors for the prediction of the parameters of the dynamical system. The bottom row shows the errors for the prediction of the covariance entries. We show (left column to right column) the minimum, mean, and maximum error.}\label{fig:hm_add_viz_2D}
\end{figure}

\begin{figure}
\centering
\includegraphics[width=0.9\textwidth]{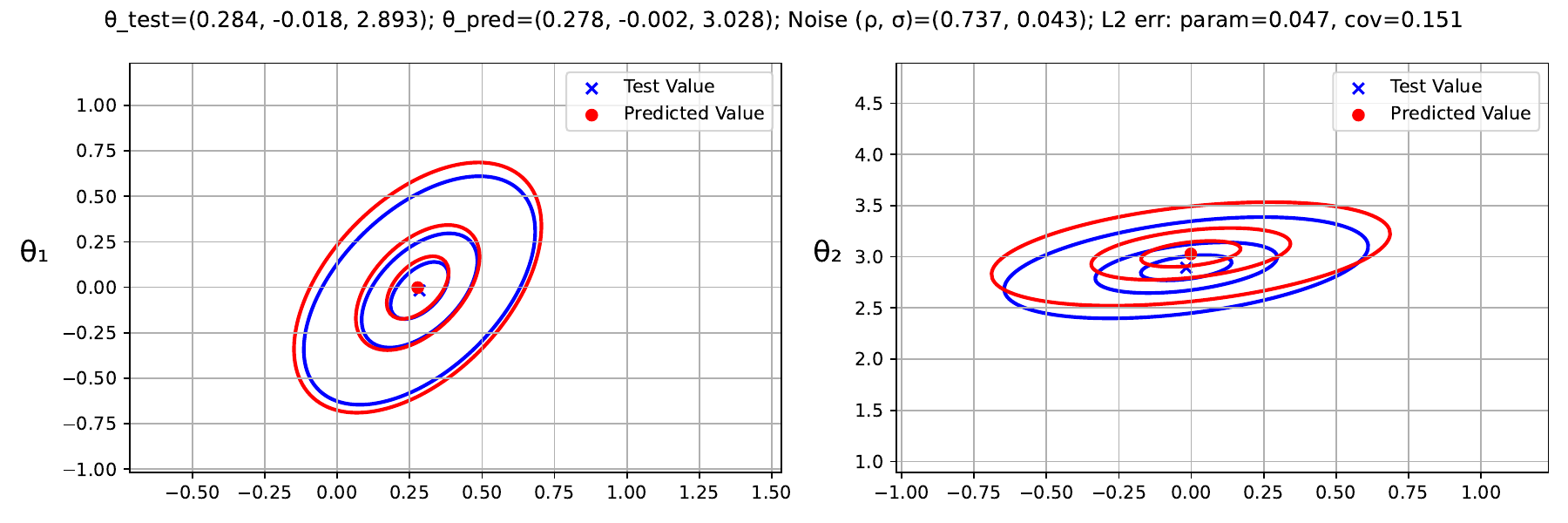}
\includegraphics[width=0.9\textwidth]{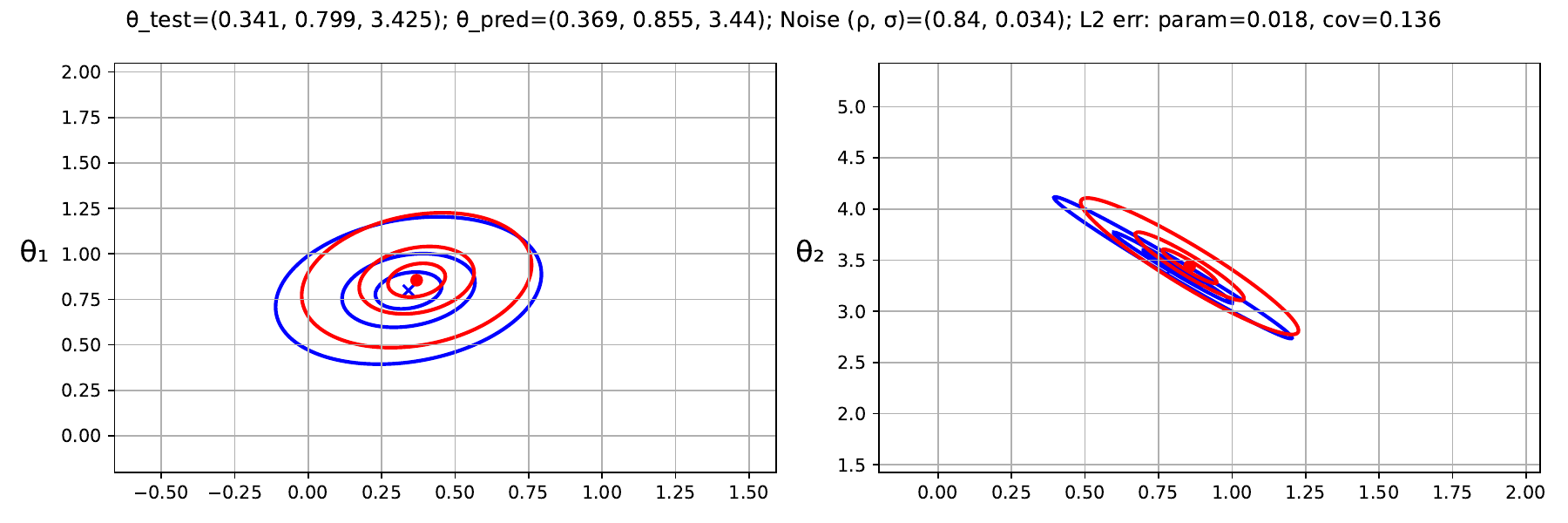}
\includegraphics[width=0.9\textwidth]{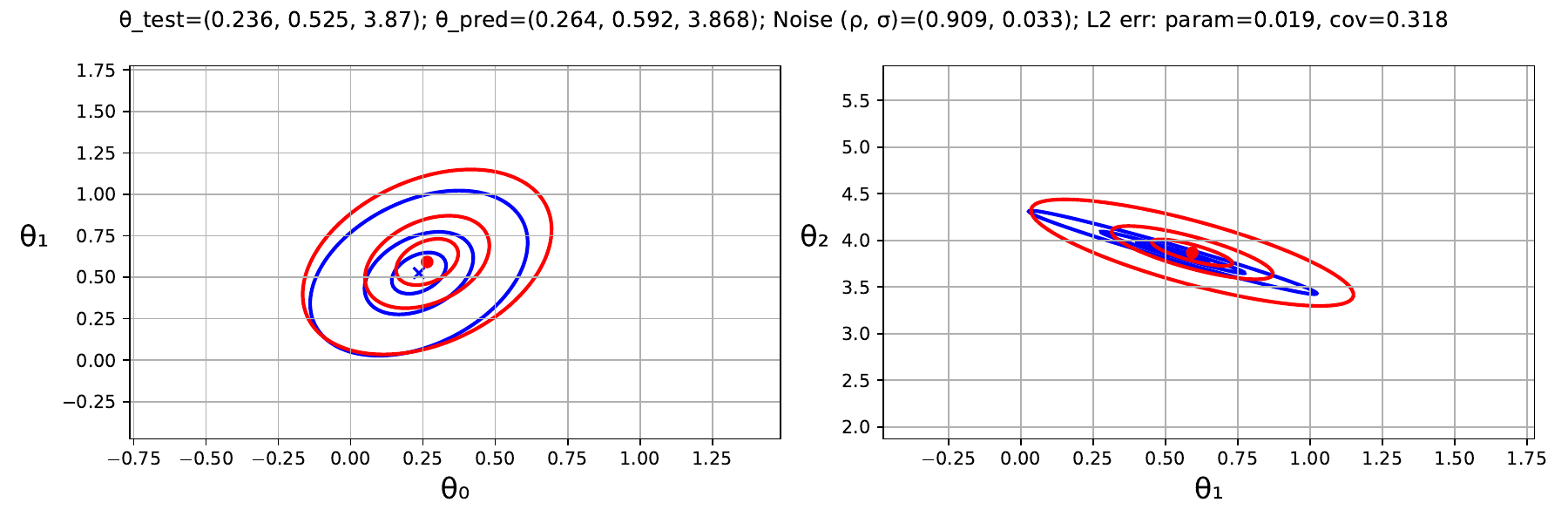}
\caption{2D projections of covariance pairs for increasing noise levels (from top to bottom). We compare the NN prediction to the test data. The predictions are based on observations generated under the additive noise model. Each subplot shows 2D uncertainty ellipses computed from the underlying covariance matrices. Each row corresponds to the 10th, 50th, and 90th percentile noise levels (based on the $\ell$-norm of the noise vector $\vect{\theta}_{\text{noise}} = (\sigma,\rho)$). The left column corresponds to the projection into the $\theta_1$-$\theta_0$ plane. The right column corresponds to the projection onto the $\theta_1$-$\theta_2$ plane. In the title of each pair of the plots (i.e., for each row), we include values for the predicted and true parameters of the dynamical system, the noise parameters, and the relative $\ell^2$-norm for the prediction of the parameters of the dynamical system and the covariance matrix.}\label{fig:cov_ellipse_viz_add}
\end{figure}

\begin{figure}
\centering
\includegraphics[width=0.9\textwidth]{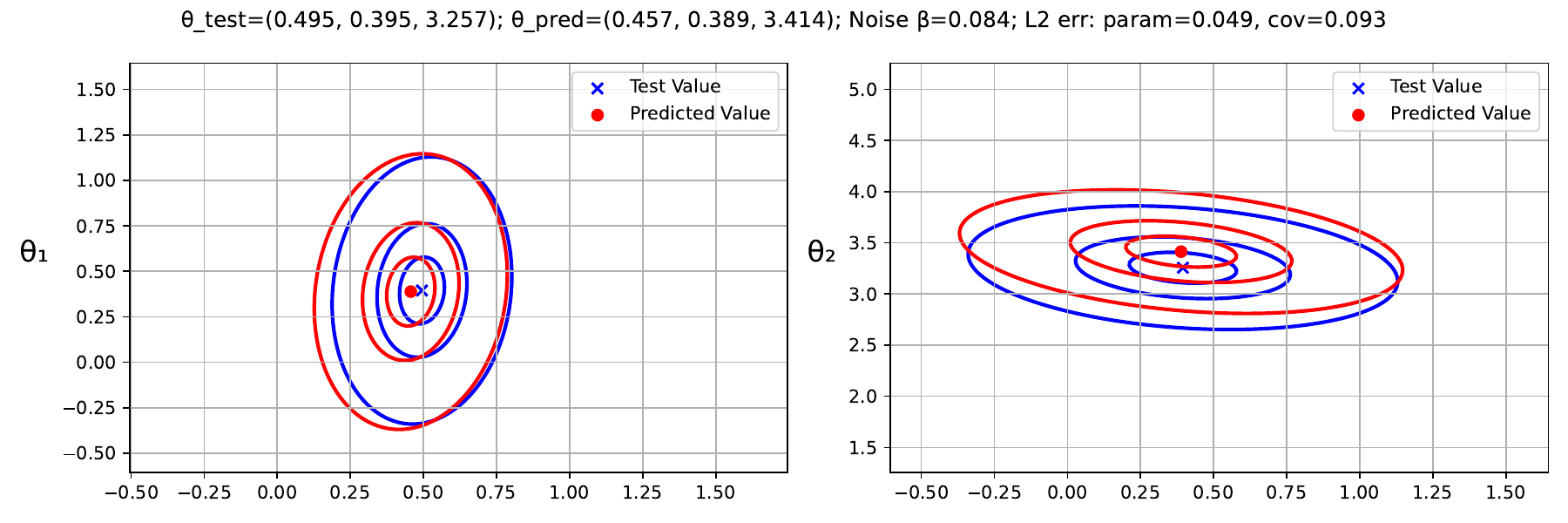}
\includegraphics[width=0.9\textwidth]{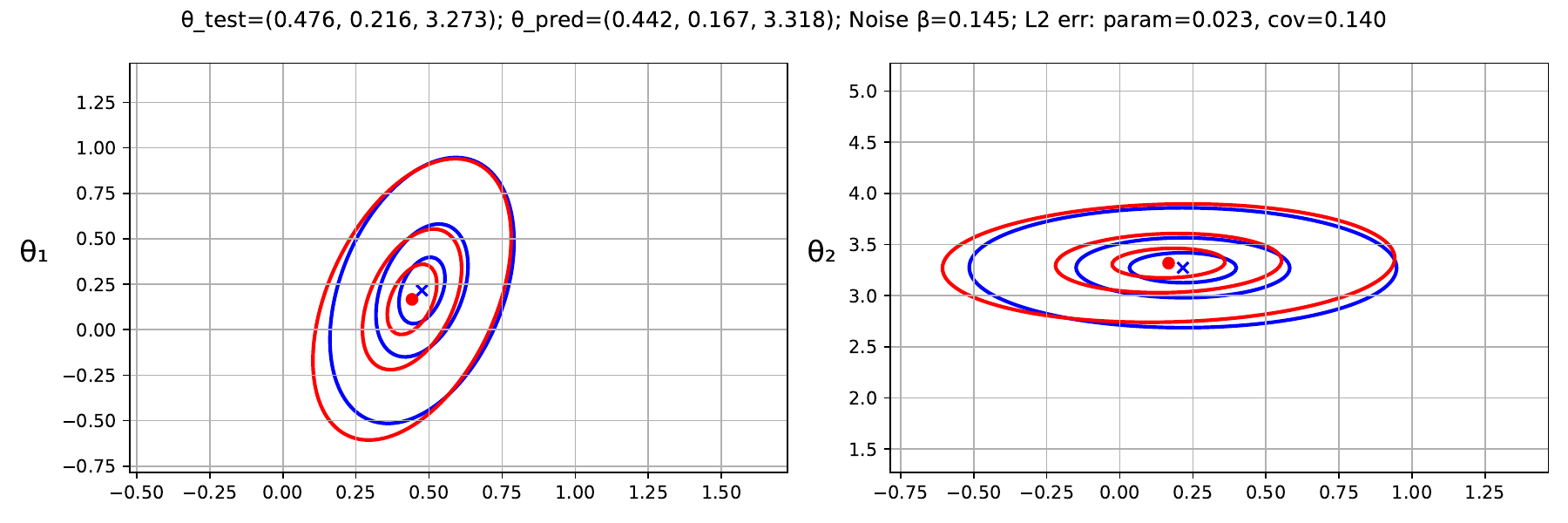}
\includegraphics[width=0.9\textwidth]{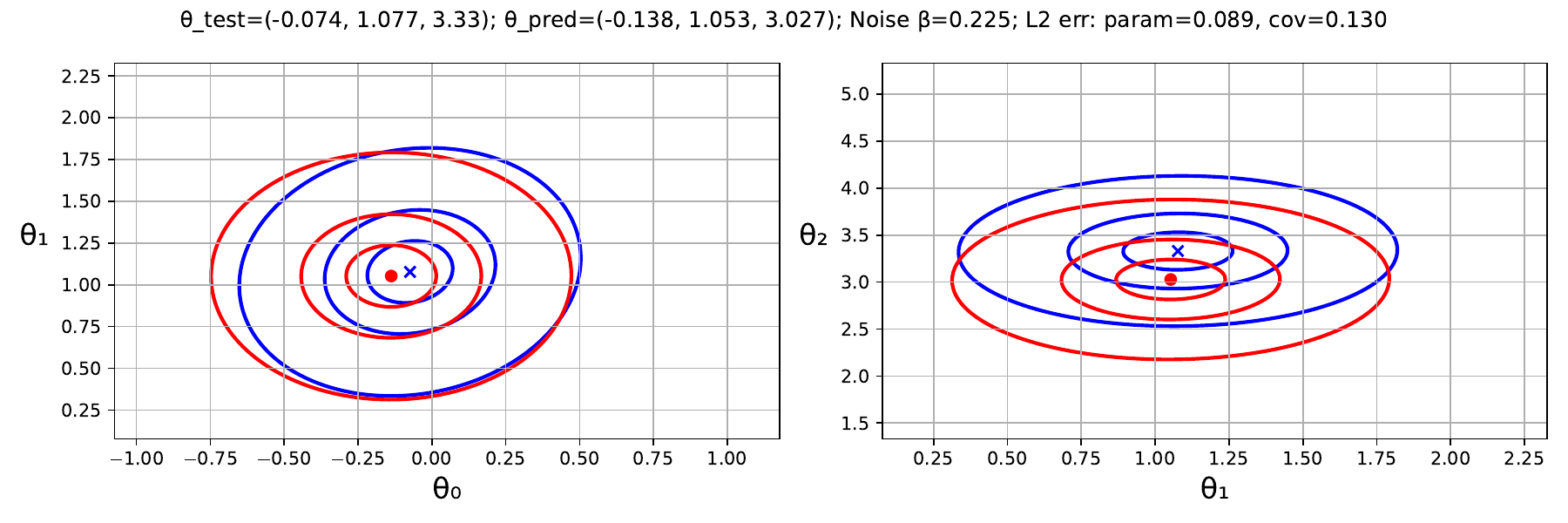}
\caption{2D projections of covariance pairs for increasing noise levels (from top to bottom). We compare the NN prediction to the test data. The predictions are based on observations generated under the intrinsic noise model. Each subplot shows 2D uncertainty ellipses computed from the underlying covariance matrices. Each row corresponds to the 10th, 50th, and 90th percentile noise levels (for the absolute value of the noise parameter $\vect{\theta}_{\text{noise}} = \beta$). The left column corresponds to the projection into the $\theta_1$-$\theta_0$ plane. The right column corresponds to the projection onto the $\theta_1$-$\theta_2$ plane. In the title of each pair of the plots (i.e., for each row), we include values for the predicted and true parameters of the dynamical system, the noise parameters, and the relative $\ell^2$-norm for the prediction of the parameters of the dynamical system and the covariance matrix.}\label{fig:cov_ellipse_viz_sde}
\end{figure}

\ipoint{Results:} We provide 2D maps for the errors in the prediction of the parameters of the dynamical system and the covariance entries for the additive noise model in \Cref{fig:hm_add_viz_2D}. We report errors as a function of the noise parameters $\vect{\theta}_{\text{noise}} = (\sigma,\rho)$. We report results for the prediction of the parameters of the dynamical system (top row) and the covariance entries (bottom row). We include (left column to right column) the smallest error, the mean error, and the largest error across multiple datasets that share the same noise perturbations. We also include some visualization for the prediction of covariance matrices to provide some intuition. Results for increasing noise levels for the additive noise model are shown in \Cref{fig:cov_ellipse_viz_add}. Similar results are shown in \Cref{fig:cov_ellipse_viz_sde} for the intrinsic noise model.

\ipoint{Observations:} The results included in this section again demonstrate that our framework remains relatively stable as the noise level increases. Compared to the box-whisker plots for the additive noise included in the main manuscript, the 2D maps shown in \Cref{fig:hm_add_viz_2D} seem to hint at a trend that the error increases as both noise parameters increases (especially for the mean values shown in the middle columns). However, overall the trend is not as striking as for the intrinsic noise model (see \Cref{s:noise_sensitivity} in the main manuscript).

\end{appendix}

\end{document}